\documentclass[12pt,a4paper]{article}
\usepackage[all]{xy}
\usepackage{amscd,amsmath,amssymb,amsfonts}

\numberwithin{equation}{section}

\newtheorem{thm}{Theorem}[section]
\newtheorem{prop}[thm]{Proposition}
\newtheorem{lem}[thm]{Lemma}
\newtheorem{cor}[thm]{Corollary}

\newtheorem{conj}[thm]{Conjecture}
\newtheorem{claim}[thm]{Claim}
\newtheorem{defn}[thm]{Definition}
\newtheorem{rem}[thm]{Remark}
\def\Th#1.{\vskip 6pt \medbreak\noindent{\bf Theorem #1.}}

\newenvironment{pf}{\smallskip\noindent\emph{Proof.}}{Q.E.D.\bigbreak}
\newcommand{\bysame}{\mbox{\rule{3em}{.4pt}}\,}

\begin{document}

%


\def\Spec{{\mathrm{Spec}}}
\def\Pic{{\mathrm{Pic}}}
\def\Ext{{\mathrm{Ext}}}
\def\NS{{\mathrm{NS}}}
\def\CH{{\mathrm{CH}}}
\def\deg{{\mathrm{deg}}}
\def\dim{{\mathrm{dim}}}
\def\codim{{\mathrm{codim}}}
\def\Coker{{\mathrm{Coker}}}
\def\ker{{\mathrm{ker}}}
\def\Image{{\mathrm{Image}}}
\def\Aut{{\mathrm{Aut}}}
\def\Hom{{\mathrm{Hom}}}
\def\Proj{{\mathrm{Proj}}}
\def\Sym{{\mathrm{Sym}}}
\def\Image{{\mathrm{Image}}}
\def\Gal{{\mathrm{Gal}}}
\def\GL{{\mathrm{GL}}}
\def\SL{{\mathrm{SL}}}
\def\End{{\mathrm{End}}}
\def\P{{\mathbb P}}
\def\C{{\mathbb C}}
\def\R{{\mathbb R}}
\def\Q{{\mathbb Q}}
\def\Z{{\mathbb Z}}
\def\F{{\mathbb F}}
\def\l{\ell}
\def\lra{\longrightarrow}
\def\ra{\rightarrow}
\def\hra{\hookrightarrow}
\def\ot{\otimes}
\def\op{\oplus}
\def\vg{\varGamma}
\def\O{{\cal{O}}}
\def\ol#1{\overline{#1}}
\def\wt#1{\widetilde{#1}}
\def\us#1#2{\underset{#1}{#2}}
\def\os#1#2{\overset{#1}{#2}}
\def\lim#1{\us{#1}{\varinjlim}}
\def\plim#1{\us{#1}{\varprojlim}}

%


\def\Gr{{\mathrm{Gr}}}
\def\rank{{\mathrm{rank}}}
\def\dlog{{\mathrm{dlog}}}
\def\Res{{\mathrm{Res}}}
\def\tor{{\mathrm{tor}}}
\def\reg{{\mathrm{reg}}}

\def\K{{\cal{K}}}
\def\E{{\cal{E}}}
\def\G{{\mathbb{G}}}
\def\dR{{\mathrm{dR}}}
\def\DR{{\mathrm{DR}}}
\def\Zar{{\mathrm{Zar}}}
\def\gen{{\mathrm{gen}}}
\def\NS{{\mathrm{NS}}}
\def\et{{\text{\'et}}}
\def\crys{{\mathrm{crys}}}
\def\Gr{{\mathrm{Gr}}}
\def\dego{\deg=0}
\def\Qb{\bar{\Q}}
\def\ord{{\mathrm{ord}}}
\def\tr{{\mathrm{tr}}}
\def\Tr{{\mathrm{Tr}}}
\def\Frob{{\mathrm{Frob}}}
\def\D{\Phi}
\def\kp{\phi_{\underline{\mu}}}
\def\cK{{\cal{K}}}
\def\cR{{\cal{R}}}

\title
{On dlog image of
$K_2$ of elliptic surface minus singular fibers}
\author{Masanori Asakura}

\date\empty

\maketitle
\begin{abstract}
Let $\pi:X\to C$ be an elliptic surface 
over an algebraically closed field of characteristic zero and
$D_i=\pi^{-1}(P_i)$ the singular fibers. 
Put $U=X-\cup_i D_i$. Our objective in this paper
is the image of the dlog map
$\vg(U,\K_2)\to \vg(X,\Omega^2_X(\log \sum D_i))$.
In particular, we give an upper bound of the rank of
the dlog image, which is computable in many cases.
This also allows us to construct
indecomposable parts of Bloch's higher Chow groups
$\CH^2(X,1)$ in special examples.
\end{abstract}
\section{Introduction}\label{intrsect}
Let $U$ be a nonsingular variety
over a field of characteristic zero.
Let
$X$ be a smooth compactification of $U$ such that
$D:=X-U$ is a normal crossing divisor.
Then there is the dlog map
$$
\dlog:\vg(U,\K_2)\lra \vg(X,\Omega^2_X(\log D))
$$
from the $K$-cohomology group to the space of the
algebraic differential 2-forms on $X$ with logarithmic poles
along $D$. When the base field is $\C$,
the dlog image is contained in
the Betti cohomology group $H^2_B(U,\Q(2))$:
\begin{equation}\label{beilinsonconjecture}
\dlog\vg(U,\K_2)\ot_\Z\Q\subseteq \vg(X,\Omega^2_X(\log D))\cap
H^2_B(U,\Q(2)).
\end{equation} 
In the paper \cite{bei}
Beilinson conjectured that the equality always holds true
in \eqref{beilinsonconjecture}
(see also \cite{jan} 5.18, 5.20).
However very little is known about it.

\smallskip

In the present paper we study the dlog image 
in a different way from the above
when $U$ is the complement of the singular fibers in an elliptic
surface.
The main results are Theorems A and B.

Theorem A concerns $K_2$ of Tate curves.
Let $K$ be a finite unramified extension of $\Q_p$
and $R$ its integer ring.
Let $E_q$ be the Tate curve over the ring $R((q)):=R[[q]][q^{-1}]$ 
of formal power series
with coefficient in $R$.
Then Theorem A gives a criterion for 2-forms on $E_q$
to be contained
in dlog image of $K_2(E_q)$.
It also gives a bounding space 
for $K_2$ of
elliptic surface over $R$.
Namely let $\pi_R:X_R\to C_R$ be an elliptic surface over $R$
which has $R$-rational multiplicative fibers $D_R$
(see \S \ref{setupsect} for the precise definition).
Put $X_K:=X_R\times_RK$, $X_{\ol{K}}:=X_R\times_R\ol{K}$ etc.
Let $U_{\ol{K}}$ be the complement of the all singular fibers.
Then we introduce a $\Z_p$-submodule $\Phi(X_R,D_R)_{\Z_p}\subset
\vg(X_R,\Omega^2_{X_R/R}(\log D_R))$
and show 
\begin{equation}\label{beilinsonconjecture2}
\dlog\vg(U_{\ol{K}},\K_2)\ot_\Z
\Q_p\subseteq\Phi(X_R,D_R)_{\Z_p}\ot_{\Z_p}\Q_p
\end{equation}
under a mild condition on $p$ (Theorem \ref{key2}).

Theorem B concerns an upper bound of the rank of
the dlog image of $K_2$ of elliptic surface minus 
singular fibers.
Due to Theorem A
 we get $\rank_\Z\dlog\vg(U_{\ol{K}},\K_2)\leq \rank_{\Z_p}
\Phi(X_R,D_R)_{\Z_p}$.
However this is not enough to compute an explicit bound
as it is usually difficult to compute the $\Z_p$-rank
of $\Phi(X_R,D_R)_{\Z_p}$. Theorem B enables us
to replace $\Phi(X_R,D_R)_{\Z_p}$ with a finite dimensional
$\F_p$-space
$\Phi(X_R,D_R)_{\F_p}$, 
which will be introduced in \S \ref{giving}.
Theorems A and B imply that
\begin{equation}\label{beilinsonconjecture3}
\rank_\Z\dlog\vg(U_{\ol{K}},\K_2)\leq \dim_{\F_p}
\Phi(X_R,D_R)_{\F_p}
\end{equation}
under some conditions on $p$ and $X_R$
(Theorem \ref{key3}).


One can compare \eqref{beilinsonconjecture2} or
\eqref{beilinsonconjecture3} with
\eqref{beilinsonconjecture}.
I expect that the equality also holds in \eqref{beilinsonconjecture2}
(Conjecture \ref{conj1}).
If it is true,
then the equality in \eqref{beilinsonconjecture3} also holds
(see Conjecture \ref{conj2} and 
the remark after it).
Conjecture \ref{conj1} is true in the modular case (
Theorem \ref{modularthm}).
However I have no idea how to attack
this problem in general.

\smallskip

Our $\Phi(X_R,D_R)_{\Z_p}$ and $\Phi(X_R,D_R)_{\F_p}$
are defined from the Fourier
expansions of 2-forms at the neighborhoods of multiplicative fibers.
They are completely different from the right hand side
of \eqref{beilinsonconjecture}.
One of
the advantage is that $\Phi(X_R,D_R)_{\F_p}$ is
computable in many examples.
In fact we will give the following examples in \S \ref{expmsect}.
\begin{cor}[Theorem \ref{two}]\label{dmain}
Let $\pi:X\to\P^1$ be the minimal elliptic surface over 
$\C$ 
such that the general fiber $\pi^{-1}(t)$ is the elliptic curve
defined by $Y^2=X^3+X^2+t^n$.
Let $U$ be the complement of the all singular fibers in $X$.
Then we have 
$$
\rank~\dlog\vg(U,\K_2)=2 \quad\text{for }n=2,3,5,7,11,13,17,19,23,29.
$$
The dlog image is generated by
$$
\dlog\left\{
\frac{Y-X}{Y+X},-\frac{t^n}{X^3}
\right\},\quad
\dlog\left\{
\frac{iY-(X+2/3)}{iY+(X+2/3)},-\frac{t^n+4/27}{(X+2/3)^3}
\right\}.
$$
\end{cor}
Since $\pi:X\to\P^1$ is defined over $\ol{\Q}$,
it seems difficult to obtain the above result
without using $\D(X_R,D_R)_{\F_p}$.
For example
I do not know whether $\dim_\Q F^2\cap H^2_B(U,\Q(2))=2$.

Computations of the rank of dlog image can be applied to
the constructions of the indecomposable parts of the Adams weight piece
$K_1(X)^{(2)}$ of $X$ (which is isomorphic
to Bloch's higher Chow group $\CH^2(X,1)\ot\Q$)
in special examples.
In fact using Corollary \ref{dmain} together with
Stiller's 
computations on N\'eron-Severi groups (\cite{stiller1}, \cite{stiller2}),
we can obtain the following:
\begin{cor}[Theorem \ref{dectwo}]\label{dmaintwo}
Let $\pi:X\to\P^1$ be as above and
$D_i$ $(1\leq i \leq n+1)$ the multiplicative fibers. 
Let $K_1^{\mathrm{ind}}(X)^{(2)}$ denotes the indecomposable $K_1$
$$
K_1^{\mathrm{ind}}(X)^{(2)}\os{\mathrm{def}}{=}
K_1(X)^{(2)}/(\C^*\ot \NS(V))$$
where $NS(X)$ denotes the N\'eron-Severi group.
Then we have
$$
\dim~\Image(\bigoplus_{i=1}^{n+1}K'_1(D_i)\lra
K_1^{\mathrm{ind}}(X)^{(2)})= n-1
\quad\text{for }n=7,11,13,17,19,23,29.
$$
\end{cor}
I do not know how to show the non-vanishing of
the regulator image of $\op_i K'_1(D_i)$ 
in the 
indecomposable part
$H^3_D(X,\Q(2))/(\C^*\ot \NS(X))$
of the Deligne-Beilinson cohomology group.

\medskip

\noindent{\it Acknowledgements}.
I would like to express special thanks to
Dr. Seidai Yasuda who pointed out that 
Theorem A in an earlier version can be improved by using
Kato's explicit reciprocity law.
Before his comments, I used the Artin-Hasse formula for the proof and
was dissatisfied with a half-finished statement
(Remark \ref{yasuda}).
I also thank Professor Kazuya Kato who informed me his papers
\cite{kato1} and \cite{kato2} with many valuable suggestions.
Most parts of this paper was written during my stay
at the University of Chicago in September 2005 supported by
JSPS Postdoctoral Fellowships for Research Abroad.
I thank for their hospitality,
especially to Professor Spencer Bloch.

\section{Preliminaries}\label{definitionsect}
For an abelian group $M$,
we denote by $M[n]$ (resp. $M/n$)
the kernel (resp. cokernel) of multiplication by $n$.
$M_\tor$ denotes the torsion subgroup of $M$.

\subsection{Algebraic $K$-theory and $K$-cohomology}
Let $X$ be a separated noetherian scheme.
Let $P(X)$ be the exact category of locally free sheaves, and
$BQP(X)$ the simplicial set attached to $P(X)$ by Quillen (\cite{Q},
\cite{Sr}).
{\it Quillen's higher $K$-groups} of $X$ are defined as the homotopy
groups of $BQP(X)$:
$$
K_i(X)\os{\mathrm{def}}{=}\pi_{i+1}BQP(X), \quad i\geq 0.
$$
We refer \cite{Sr} for the general properties of higher $K$-theory
such as,
products, localization exact sequences, norm maps (also called transfer
maps) etc.

Let ${\cal K}_i$ denote the Zariski
sheaf on $X$ associated to
the presheaf
$$
U\longmapsto K_i(U) \quad (U\subset X).
$$
The Zariski cohomology groups $H^j_\Zar(X,{\cal K}_i)$ are
called the {\it $K$-cohomology groups} (\cite{Sr} \S 5 etc.).
Assume that $X$ is a regular irreducible scheme
of dimension $d$.
We denote by $X^j$ the set of points of height $j$.
We write by $\kappa(x)$ the residue field of a point
$x\in X$
and $i_x:\{x\}\to X$ the inclusion.
Then we have the complex
\begin{equation}\label{gerstencomplex}
0\to {\cal K}_i \to K_i(\kappa(\eta)) \to
\bigoplus_{x\in X^1} i_{x*}K_{i-1}(\kappa(x)) \to \cdots
\to \bigoplus_{x\in X^{d}} i_{x*}K_{i-d}(\kappa(x))
\to 0
\end{equation}
of Zariski sheaves where $\eta$ is the generic point.
\begin{thm}[Gersten conjecture]\label{gersten}
The complex \eqref{gerstencomplex} is exact
in either of the following cases.
\begin{enumerate}
\item
$i\geq 0$ and
$X$ is a nonsingular variety over a field,
\item
$i=2$ and $X$ is smooth over a Dedekind domain.
\end{enumerate}
\end{thm}
\begin{pf}
The former is due to
Quillen \cite{Q} Thm. 5.11, and
the latter is due to Bloch \cite{blochgersten}.
\end{pf}
Suppose that $X$ is either of the cases in Theorem \ref{gersten}.
Then \eqref{gerstencomplex} gives the 
flasque resolution of the sheaf
${\cal K}_i$ so that we have the isomorphism
\begin{equation}\label{coisom}
H^j_\Zar(X,{\cal K}_i)\cong
\frac{\ker (\bigoplus_{x\in X^j} K_{i-j}(\kappa(x)) \to
\bigoplus_{x\in X^{j+1}} K_{i-j-1}(\kappa(x)) )}{\Image
(\bigoplus_{x\in X^{j-1}} K_{i-j+1}(\kappa(x)) \to
\bigoplus_{x\in X^{j}} K_{i-j}(\kappa(x)) )}
\end{equation}
Hereafter,
we often use the identification \eqref{coisom}.
In particular, we identify $\vg(X,{\cal K}_2)
=H^0_\Zar(X,{\cal K}_2)$ with the kernel of
the {\it tame symbol}
\begin{equation}\label{tame}
\tau=\bigoplus_{x\in X^1}\tau_x: K_2^M(\kappa(\eta)) \lra
\bigoplus_{x\in X^1} \kappa(x)^*
,\quad
\{f,g\}\mapsto\sum_{x\in X^1}(-1)^{{\mathrm{ord}}_x(f)
{\ord}_x(g)}\frac{f^{{\ord}_x(g)}}
{g^{{\ord}_x(f)}}.
\end{equation}
Here $K_2^M$ denotes Milnor's $K_2$.

Let $X$ be a regular scheme of dimension $d$.
There is the natural map
$K_2(X)\to \ker~\tau$, which
is surjective if tensoring with $\Z[1/(d+1)!]$
(\cite{soulecanada} Th\'eor\`eme 4 iv)).
In particular, if $X$ is either of the cases in Theorem \ref{gersten},
we have a natural surjection
\begin{equation}\label{souleger1}
K_2(X)\ot\Z[\frac{1}{(d+1)!}]
\lra\vg(X,\K_2)\ot\Z[\frac{1}{(d+1)!}].
\end{equation}

\subsection{Regulator and dlog maps on $\vg(X,\K_2)$}
Let $W_k$ be a nonsingular variety over a field $k$.
There are the {\it regulator maps} 
(also called the {\it Chern class maps})
\begin{equation}\label{dlogR0}
c_B:\vg(W_\C,\K_2)\lra H^2_B(W_\C,\Z(2))
\end{equation}
to the Betti cohomology group when $k=\C$,
\begin{equation}\label{dlogR1}
c_\et:\vg(W_k,\K_2)\lra H^2_\et(W_k,\Z/n(2))
\end{equation}
to the \'etale cohomology group when $n$ is invertible
in $k$,
and
\begin{equation}\label{dlogR2}
c_\dR:\vg(W_k,\K_2)\lra H^2_\dR(W_k/k)
\end{equation}
to the de Rham cohomology group
when $k$ is of characteristic zero
(cf. \cite{suslin} \S 23. See also \cite{gillet}, \cite{schneider}
for the general Chern class maps on Quillen's $K$-groups).

\medskip

The following theorem
will play an essential role in the proof of Theorem B.
\begin{thm}[Suslin's exact sequence ; \cite{suslin} Corollary 23.4]
\label{universalc0}
Suppose that $n$ is invertible in $k$.
There is the natural exact sequence
\begin{equation}
0\lra \vg(W_k,{\cal K}_2)/n
\os{c_\et}{\lra} H^2_\et(W_k,\Z/n(2)) \lra
H^1_\Zar(W_k,{\cal K}_2)[n]\lra 0.
\end{equation}
\end{thm}

Suppose that $k$ is of characteristic zero.
The de Rham regulator $c_\dR$ factors
through the Hodge filtration $F^2H^2_\dR(W_k/k)
=\vg(\ol{W}_k,\Omega_{\ol{W}_k}^2(\log D_k))$ where
$\ol{W}_k\supset W_k$ is a smooth compactification such that
$D_k:=\ol{W}_k-W_k$ is a normal crossing divisor.
Thus it gives rise to the {\it dlog map}
\begin{equation}\label{dlogR2d}
\dlog:\vg(W_k,\K_2)\lra \vg(\ol{W}_k,\Omega_{\ol{W}_k}^2(\log D_k)),
\end{equation}
which is written as $\sum\{f,g\}\mapsto \sum df/f.dg/g$ 
under the identification \eqref{coisom}.

Let $A$ be a discrete valuation ring with a uniformizer
$\pi$.
Suppose that the quotient field $K$ is of characteristic zero.
Let $W_A$ be a smooth scheme over $A$ and
$\ol{W}_A\supset W_A$ smooth over $A$ such that $D_A=\sum D_{i,A}
:=\ol{W}_A-W_A$ is a {\it relative} normal crossing divisor over $A$,
which means that all $D_{i,A}$ are smooth over $A$ and
any intersection locus of $D_A$ are also smooth.
Put $W_K:=W_A\times_AK$.
Then the dlog map \eqref{dlogR2d} induces
\begin{equation}\label{dlogR}
\dlog:\vg(W_K,\K_2)
\lra \vg(\ol{W}_A,\Omega^2_{\ol{W}_A/A}(\log D_A)).
\end{equation}
We can see it in the following way.
Let
$$
\Omega^j_{\ol{W}_A/A}(\log)
:=\lim{E}\vg(\ol{W}_A-E_{\mathrm{sing}},
\Omega^j_{\ol{W}_A/A}(\log E))\subset \Omega_{K(W)/K}^j
$$
where $E$ runs over all divisors on $\ol{W}_A$ such that
each irreducible component is flat over $A$ and
$E_{\mathrm{sing}}$ denotes the singular locus.
We claim that the image of the dlog map
$$
K_2^M(K(W))\lra \Omega_{K(W)/K}^2,\quad
\{f,g\}\longmapsto \frac{df}{f}\frac{dg}{g}
$$
is contained in $\Omega^2_{\ol{W}_A/A}(\log)$.
In fact let $f\in K(W)^*$. Since $d(cf)/(cf)=df/f$
for $c\in K^*$ and $\ol{W}_A$ is smooth over $A$, 
we may replace $f$ with $\pi^mf$ 
so that all irreducible components of the divisor of $f$ 
are flat over $A$.
This shows $df/f\in 
\Omega^1_{\ol{W}_A/A}(\log )$ for all $f\in K(W)^*$
and hence $df/f\cdot dg/g\in
\Omega^2_{\ol{W}_A/A}(\log )$.
Let
$$
{\mathrm{Res}}_x:\Omega^j_{\ol{W}_A/A}(\log)
\lra \Omega^{j-1}_{\kappa(x)/K}
$$
be the residue map
at $x$. Then we have a commutative diagram
\begin{equation}\label{blochger1}
\begin{CD}
K_2^M(K(W))@>{\tau}>>
\bigoplus_{x\in W_K^1}\kappa(x)^*\\
@V{\mathrm{dlog}}VV@VV{\mathrm{dlog}}V\\
\Omega^2_{\ol{W}_A/A}(\log)
@>{\op{\mathrm{Res}}_x}>> 
\bigoplus_{x\in W_K^1}\Omega^{1}_{\kappa(x)/K}.
\end{CD}
\end{equation}
The kernel of the bottom arrow is
\begin{align*}
\lim{E}\vg(\ol{W}_A-E_{\mathrm{sing}},
\Omega^2_{\ol{W}_A/A}(\log D_A))&\subset
\lim{\codim Z\geq 2}\vg(\ol{W}_A-Z,
\Omega^2_{\ol{W}_A/A}(\log D_A))\\
&=
\vg(\ol{W}_A,
\Omega^2_{\ol{W}_A/A}(\log D_A)).
\end{align*}
The last equality follows from the fact that 
$\Omega^2_{\ol{W}_A/A}(\log D_A)$
is locally free of finite rank.
Thus the diagram \eqref{blochger1} together with
the identification \eqref{coisom} gives rise to
\eqref{dlogR}.

\subsection{Elliptic surface and
multiplicative fiber}\label{setupsect}
Let $A$ be a field or a discrete valuation ring.
In this paper,
we mean by an {\it elliptic surface over $A$}
a projective flat morphism $\pi_A:X_A\to C_A$ of $A$-schemes
such that
\begin{enumerate}
\renewcommand{\theenumi}{(\roman{enumi})}
\item
$X_A$ and $C_A$ are projective and smooth schemes
over $A$ of relative dimension $2$ and $1$ respectively,
\item
The general fiber of $\pi_A$ is an elliptic curve,
\item
$\pi_A$ has a section $e_A:C_A\to X_A$.
\end{enumerate}
Let $K$ be the quotient field of $A$ and $\ol{K}$ the algebraic
closure.
Put $X_K:=X_A\times_AK$ and $X_{\ol{K}}:=X_A\times_A\ol{K}$
etc.
If all fibers of $\pi_K:X_K\to C_K$
are free from $(-1)$-curves,
we call it {\it minimal}.
It is well-known that there is a unique minimal
elliptic surface $\pi'_K:V_K\to C_K$ 
over $K$
with a surjective
map $\rho_K:X_K\to V_K$ 
such that $\pi'_K\rho_K=\pi_K$:
$$
\xymatrix{
X_K\ar[rd]_{\pi_K}\ar[rr]^{\rho_K}
& &V_K\ar[ld]^{\pi'_K}\\
&C_K 
}
$$
The classification of the singular fibers of
minimal elliptic surface is well-known thanks to
works of Kodaira and N\'eron.
I refer the reader to the book \cite{silverman} IV \S 8.
Note that if $\pi_K$ has a section, then there is 
no multiple fiber (often denoted by $\empty_nI_m$).

\medskip

Let $\tilde{C}:=\P^1_\Z\times \Z/m$ be the disjoint union of
copies of $\P^1_\Z$, indexed by $\Z/m$.
Attach the point $0$ of $i$-th $\P^1_\Z$ to the point
$\infty$ of $(i+1)$-th $\P^1_\Z$. Then we have a connected 
proper curve
$C_\Z$ over $\Spec\Z$ with the normalization $\tilde{C}\to C_\Z$.
We call $C_\Z\times S$ the {\it standard N\'eron polygon} (or
{\it standard $m$-gon}) over a scheme $S$
(cf. \cite{DeRa} II. 1.1).

Let $D_{\ol{K}}$ be a singular fiber of $\pi_{\ol{K}}$.
If there is a closed subscheme 
$D^\dag_{\ol{K}}\subset D_{\ol{K}}$ which is isomorphic to
a standard $m$-gon, then we call $D_{\ol{K}}$ a 
{\it multiplicative fiber} or type $I_m$.
This is equivalent to say that $\rho_{\ol{K}}(D_{\ol{K}})$
is a standard $m$-gon.

\medskip

If the elliptic surface $\pi_A:X_A\to C_A$ satisfies
the following condition, we say that it has {\it $A$-rational}
multiplicative fibers:
\begin{description}\item[(Rat)]
There is a closed subscheme $\Sigma_A\subset C_A$ which is a disjoint
union of finite copies of $\Spec A$ such that
\begin{enumerate}
\renewcommand{\theenumi}{(\roman{enumi})}
\item
a singular fiber $\pi_{\ol{K}}^{-1}(P)$ over $P\in C_A(\ol{K})$
is multiplicative
if and only if $P\in \Sigma_A(\ol{K})$,
\item
for each $A$-rational point $P$ of $\Sigma_A$, 
there is a closed subscheme
$D^\dag_A\subset  \pi^{-1}_A(P)$ which is isomorphic to
a standard N\'eron polygon over $A$.
\end{enumerate}
\end{description}
When $A=\ol{K}$, the above is automatically satisfied.

\medskip

\noindent\underline{Notation}.
Let $\pi_A:X_A\to C_A$ be an elliptic surface over a 
discrete valuation ring or a field $A$
satisfying
{\bf (Rat)}.
Let $\Sigma_A=\{P_1,\cdots,P_s\}$ with $P_i\cong \Spec A$.
We put $S_A:=C_A-\Sigma_A$,
$D_{i,A}:=\pi_A^{-1}(P_{i,A})$,
$D_A:=\sum D_{i,A}$ and $U_A:=X_A-D_A$:
$$
\begin{CD}
U_A@>>> X_A\\
@V{\pi_A}VV@VV{\pi_A}V\\
S_A@>>>C_A.
\end{CD}
$$
We put by $r_i\geq 1$ the number of the irreducible components
of the standard N\'eron polygon $D_{i,A}^\dag\subset D_{i,A}$
(i.e. $D_{i,A}^\dag$ is a standard $r_i$-gon).
Equivalently, $r_i$ is the pole order of the functional $j$-invariant
of $X_{\ol{K}}\to C_{\ol{K}}$ at $P_i$.
Moreover
let $S^0_A\subset S_A$ be an arbitrary
open subscheme such that
$T_A:=(S_A-S_A^0)_{\mathrm{red}}$ is flat over $A$
where `${\mathrm{red}}$' denotes the reduced subscheme.
We
put $C^0_A:=C_A-T_A$, $Y_A:=\pi_A^{-1}(T_A)_{\mathrm{red}}$,
$U^0_A:=X_A-D_A-Y_A$
and $X_A^0:=X_A-Y_A$:
$$
\begin{CD}
U^0_A@>>>X^0_A\\
@V{\pi_A}VV@VV{\pi_A}V\\
S^0_A@>>>C_A^0.
\end{CD}
$$

\subsection{Boundary maps}\label{boundarysect}
Let $\pi_A:X_A\to C_A$ be an elliptic surface over a 
discrete valuation ring or a field $A$
satisfying
{\bf (Rat)}.
For each multiplicative fiber $D_{i,A}$,
we fix an irreducible component $Z_{i,A}\subset D^\dag_{i,A}$
of the standard N\'eron polygon
and a singular point $Q_{i,A}$ of $D^\dag_{i,A}$ such that
$Q_{i,A}\in Z_{i,A}$.
Put $Z_{i,A}^*=Z_{i,A}\cap D^{\dag,\mathrm{reg}}_{i,A}$ 
where $D^{\dag,\mathrm{reg}}_{i,A}$ denotes
the regular locus. 
Then
$Z_{i,A}^*$ is isomorphic to ${\mathbb G}_{m,A}$.
Let $\tau_i:\vg(U^0_A,\K_2)\to \vg(Z^*_{i,A},\K_1)$
be the tame symbol \eqref{tame} at $Z^*_{i,A}$
and
$\ord_i:\vg(Z^*_{i,A},\K_1)
\to \Z$ the map of order at $Q_{i,A}$.
Put $\partial_i:=\ord_i\cdot \tau_i$.
We call
the following map
the {\it boundary map}
in algebraic $K$-theory:
\begin{equation}\label{phikthe0}
\partial=\bigoplus_{i=1}^s\partial_i:
\vg(U^0_A,\K_2)
\lra \Z^{\op s}.
\end{equation}
It is easy to see that $\partial_i$ does not depend on
the choices of $Z_{i,A}$ nor $Q_{i,A}$ up to sign
(cf. \cite{scholl} 1.5 or \S \ref{modelsect}
below).
Since our objective is the image of the boundary,
the sign ambiguity does not matter.
We write the composition $K_2(U^0_A)\to\vg(U_A^0,\K_2)
\os{\partial}{\to} \Z^{\op s}$ by the same notation $\partial$.

There are the corresponding maps (which we also call 
the {\it boundary maps})
\begin{equation}\label{phidr01}
\partial_{B}:H^2_B(U^0_\C,\Z(2))\lra \Z^{\op s}
\end{equation}
on the Betti cohomology when $A=\C$
and
\begin{equation}\label{phidr02}
\partial_{\et}:H^2_\et(U^0_{\ol{k}},\Z/n(2))\lra \Z/n^{\op s}
\end{equation}
on the \'etale cohomology when $A=\ol{k}$ 
is an algebraically
closed field in which $n$ is invertible. 
They are
compatible with $\partial$ under the regulator maps 
$c_B$ and $c_\et$.
The map $\partial_B$ is defined as the composition of
the following maps
$$
H^2_B(U^0_\C,\Z(2))
\lra
\bigoplus_{i=1}^sH^1_B(Z^*_{i,\C},\Z (1))
\os{\cong}{\lra}
\bigoplus_{i=1}^s\Z.
$$
The definition of $\partial_\et$
is similar.
Note that $\partial_B$ (or $\partial_\et$) is also defined
as the composition of
$H^2_B(U^0_\C,\Z(2))\to H^3_{D+Y,B}(X_\C,\Z(2))$ and the 
Poincare-Lefschetz duality
$H^3_{D,B}(X_\C,\Z(2))\cong H_1(D_\C,\Z)=\Z^{\op s}$. 
In particular, $\partial_B$ is a homomorphism of mixed Hodge
structure (we put the trivial Hodge structure of the right hand side of
\eqref{phidr01}), and $\partial_\et$ is compatible in
$G_k$-action.

\medskip

Let $K$ be the quotient field of $A$.
Suppose that the characteristic of $K$ is zero.
Put $X_K:=X_A\times_AK$ etc.
Then the dlog map \eqref{dlogR} induces
\begin{equation}\label{dlog0}
\dlog:
\vg(U_K^0,\K_2)
\lra \vg(X_A,\Omega^2_{X_A/A}(\log D_A)).
\end{equation}
We can see it in the following way.
A priori the dlog image is contained in $
\vg(U_K^0,\Omega^2_{U^0_K/K})$.
It follows from \eqref{dlogR} that we have
$\dlog
\vg(U_K^0,\K_2)
\subset \vg(X^0_A,\Omega^2_{X_A/A}(\log D_A)).
$
Moreover we claim that 
$\dlog
\vg(U_K^0,\K_2)
\subset \vg(X_K,\Omega^2_{X_K/K}(\log D_K)).
$
To see this, we may replace $K$ with $\C$.
Then it is enough to show
$F^2H_\dR^2(U^0_\C/\C)\cap H_B^2(U_\C^0,\Q(2))
\subset F^2H_\dR^2(U_\C/\C)=
\vg(X_\C,\Omega^2_{X_\C/\C}(\log D_\C))$.
However it follows from the exact sequence
$$
H^2_{Y,B}(U_\C,\Q)\lra
H^2_B(U_\C,\Q)\lra 
H^2_B(U^0_\C,\Q)\lra H^3_{Y,B}(U_\C,\Q)
$$
of mixed Hodge structures and the fact that
$H^2_{Y,B}(U_\C,\Q)\cong\Q(-1)^{\op}$
and $H^3_{Y,B}(U_\C,\Q)$ is of type $(1,2)+(2,1)$ 
(cf. \cite{d2} III).
Thus we get
\begin{align*}
\dlog:
\vg(U_K^0,\K_2)
\lra &\vg(X^0_A,\Omega^2_{X_A/A}(\log D_A))
\cap \vg(X_K,\Omega^2_{X_K/K}(\log D_K))\\
&=\vg(X_A,\Omega^2_{X_A/A}(\log D_A))
\end{align*}
where the equality follows from the fact that
$\Omega^2_{X_A/A}(\log D_A)$ is a locally free sheaf.

The Poincare residue gives rise to the boundary map
\begin{equation}\label{phidr0}
\partial_{\DR}:
\vg(X_A,\Omega^2_{X_A/A}(\log D_A))
\lra A^{\op s}
\end{equation}
on the de Rham cohomology group.
It is 
compatible with $\partial$ under the dlog map \eqref{dlog0}:
\begin{equation}\label{boundarydr}
\begin{CD}
\vg(X_A,\Omega^2_{X_A/A}(\log D_A))
@>{\partial_{\DR}}>>A^{\op s}\\
@A{\dlog}AA@AAA\\
\vg(U^0_K,\K_2)
@>{\partial}>>\Z^{\op s}.
\end{CD}
\end{equation}
\begin{lem}\label{d=b}
$\dlog\vg(U^0_A,\K_2)\cong\partial\vg(U^0_A,\K_2)$.
In particular, $0\leq\rank_\Z\dlog\vg(U^0_K,\K_2)\leq s$. 
\end{lem}
\begin{pf}
It is enough to show that $\partial_\dR$ is injective on 
the image $\dlog\vg(U_A^0,\K_2)$.
To do this, we may replace $A$ with $\ol{K}$.
Moreover it follows from a standard argument 
that we may assume $\ol{K}=\C$.
$\partial_\dR$ is compatible with $\partial_B$
and the dlog image is contained in 
$H^2_B(U^0_\C,\Z(2))/(\text{torsion})$.
Therefore it is enough to show that the following map
\begin{equation}\label{parinj0}
F^2H^2_\dR(U^0_\C/\C)\cap 
\left(H^2_B(U^0_\C,\Z(2))/\text{torsion}\right)
\os{\partial_B}{\lra} \Z^{\op s}
\end{equation}
is injective. Since
$\partial_B$
is a homomorphism of mixed Hodge
structures, the injectivity of \eqref{parinj0} 
follows from the Hodge symmetry.
\end{pf}
\begin{lem}\label{d=bmore}
\begin{enumerate}
\renewcommand{\labelenumi}{$(\theenumi)$}
\item
$\partial\vg(U^0_A,\K_2)=\partial\vg(U^0_K,\K_2)$.
\item 
Suppose $A=K$ (i.e. $A$ is a field). Then we have
$\partial\vg(U_K,\K_2)=\partial\vg(U^0_K,\K_2)$
and
$\partial\vg(U^0_K,\K_2)\ot\Q=\partial\vg(U^0_{\ol{K}},\K_2)\ot\Q$.
\end{enumerate}
\end{lem}
\begin{pf}
(2).
Let $k$ be the residue field of $K$ and $U^0_k$ and $S_k$
the special fibers over $\Spec k$.
It follows from a commutative diagram
$$
\begin{CD}
@.0 @.0\\
@.@VVV@VVV\\
@.k(U^0_k)^*@>>>\bigoplus_{x\in (U^0_k)^1}\Z\\
@.@VVV@VVV\\
K_2^M(K(U^0))@>>>\bigoplus_{x\in (U^0_A)^1}\kappa(x)^*@>>>
\bigoplus_{x\in (U^0_A)^2}\Z\\
@V{=}VV@VVV@VVV\\
K_2^M(K(U^0))@>>>\bigoplus_{x\in (U^0_K)^1}\kappa(x)^*
@>>>
\bigoplus_{x\in (U^0_K)^2}\Z\\
@.@VVV@VVV\\
@.0 @.0
\end{CD}
$$
and the Gersten conjecture (Theorem \ref{gersten})
that we have an exact sequence
\begin{equation}\label{vgexseq}
\vg(U^0_A,\K_2)\lra \vg(U^0_K,\K_2)\lra \vg(U^0_k,\K_1).
\end{equation}
Since $\vg(U^0_k,\K_1)\cong\vg(S^0_k,\K_1)$,
we have a commutative diagram
$$
\begin{CD}
\vg(S^0_A,\K_2)@>>> \vg(S^0_K,\K_2)@>>> \vg(S^0_k,\K_1)\\
@V{\pi_A^*}VV@V{\pi_K^*}VV@V{\pi_k^*}V{\cong}V\\
\vg(U^0_A,\K_2)@>>> \vg(U^0_K,\K_2)@>>> \vg(U^0_k,\K_1)\\
@V{e_A^*}VV@V{e_K^*}VV@V{e_k^*}V{\cong}V\\
\vg(S^0_A,\K_2)@>>> \vg(S^0_K,\K_2)@>>> \vg(S^0_k,\K_1)
\end{CD}
$$
with exact rows.
A diagram chase yields that the map
$$
\vg(U^0_A,\K_2)\op \pi_K^*\vg(S^0_K,\K_2)\lra
\vg(U^0_K,\K_2)
$$
is surjective.
It is clear from the definition that 
$\partial\pi_K^*\vg(S^0_K,\K_2)=0$.
Thus we have
$\partial\vg(U^0_A,\K_2)=\partial\vg(U^0_K,\K_2)$.

\medskip

(2).
We first show that $\partial\vg(U_K,\K_2)=\partial\vg(U^0_K,\K_2)$.
It follows from a commutative diagram
$$
\begin{CD}
@.0 @.0\\
@.@VVV@VVV\\
@.\bigoplus_{\eta\in Y_K^0}\kappa(\eta)^*@>{c}>>
\bigoplus_{x\in Y_K^1}\Z\\
@.@VVV@VVV\\
K_2^M(K(U))@>>>\bigoplus_{x\in U_K^1}\kappa(x)^*@>>>
\bigoplus_{x\in U_K^2}\Z\\
@V{=}VV@VVV@VVV\\
K_2^M(K(U))@>>>\bigoplus_{x\in (U_K^0)^1}\kappa(x)^*
@>>>
\bigoplus_{x\in (U_K^0)^2}\Z\\
@.@VVV@VVV\\
@.0 @.0
\end{CD}
$$
and the Gersten conjecture (Theorem \ref{gersten})
that we have an exact sequence
\begin{equation}\label{vgexseq=0}
\vg(U_K,\K_2)\lra \vg(U_K^0,\K_2)\lra \ker~c.
\end{equation}
We claim $\ker~c\cong \vg(T_K,\K_1)$.
In fact, let $\mu:\tilde{Y}_K\to Y_K$ be the normalization.
Then we have a commutative diagram
$$
\begin{CD}
@. 0\\
@.@VVV\\
@. \ker~\mu_*\\
@.@VVV\\
\bigoplus_{\eta\in \tilde{Y}_K^0}\kappa(\eta)^*
@>{\tilde{c}}>>\bigoplus_{x\in \tilde{Y}_K^1}\Z\\
@|@VV{\mu_*}V\\
\bigoplus_{\eta\in Y_K^0}\kappa(\eta)^*
@>{c}>>\bigoplus_{x\in Y_K^1}\Z\\
@.@VVV\\
@. 0.
\end{CD}
$$
Since $\tilde{Y}_K$ is regular, we have
$\ker~\tilde{c}\cong\vg(T_K,\K_1)$ and 
$\Coker~\tilde{c}\cong\Pic(\tilde{Y}_K)$.
Therefore it is enough to show that the map 
$\ker\mu_*\to \Pic(\tilde{Y}_K)$ is injective.
To do this, we may replace $K$ with $\ol{K}$.
Let $\Pic(\tilde{Y}_{\ol{K}})\to \Z^{\op}$ be the degree map.
Then the kernel of the composition
$\ker\mu_*\to \Pic(\tilde{Y}_{\ol{K}})\to\Z^{\op}$
is isomorphic to $H_1(\vg(Y_{\ol{K}}),\Z)$
where $\vg(Y_{\ol{K}})$ denotes the dual graph 
(cf. \cite{DeRa} I 3.5).
Since all fibers over $T_{\ol{K}}$ are not
multiplicative, it is zero.

Now we have a commutative diagram
$$
\begin{CD}
\vg(S_K,\K_2)@>>> \vg(S^0_K,\K_2)@>>> \vg(T_K,\K_1)\\
@V{\pi_K^*}VV@V{\pi_K^*}VV@V{\pi_K^*}V{\cong}V\\
\vg(U_K,\K_2)@>>> \vg(U^0_K,\K_2)@>>> \ker~c\\
@V{e_K^*}VV@V{e_K^*}VV@V{e_K^*}V{\cong}V\\
\vg(S_K,\K_2)@>>> \vg(S^0_K,\K_2)@>>> \vg(T_K,\K_1)
\end{CD}
$$
with exact rows. A diagram chase yields that the map
$$
\vg(U_K,\K_2)\op\pi_K^*\vg(S_K^0,\K_2)\lra \vg(U_K^0,\K_2)
$$
is surjective. Since $\partial\pi_K^*\vg(S_K^0,\K_2)=0$,
we have
$\partial\vg(U_K,\K_2)=\partial\vg(U^0_K,\K_2)$.

\medskip

Next we show that 
$\partial\vg(U^0_K,\K_2)\ot\Q=\partial\vg(U^0_{\ol{K}},\K_2)\ot\Q$.
There is a finite extension $L/K$ such that 
$\partial\vg(U^0_L,\K_2)=\partial\vg(U^0_{\ol{K}},\K_2)$.
We want to show $\partial\vg(U^0_K,\K_2)\ot\Q\supset
\partial\vg(U^0_L,\K_2)\ot\Q$.
Recall that the elliptic surface $X_K\to C_K$ has $K$-rational multiplicative
fibers (i.e. {\bf (Rat)}).
Therefore the norm map $N_{L/K}:\vg(U^0_L,\K_2)\to 
\vg(U^0_K,\K_2) $ on $K$-cohomology
induces a commutative diagram
$$
\begin{CD}
\vg(U^0_L,\K_2)@>{\partial}>>\Z^{\op s}\\
@V{N_{L/K}}VV@VV{\text{mult. by }[L:K]}V\\
\vg(U^0_K,\K_2)@>{\partial}>>\Z^{\op s}.
\end{CD}
$$
Thus the assertion follows.
\end{pf}
Due to Lemma \ref{d=bmore} (2), the quotient
$\partial\vg(U_{\ol{K}},\K_2)/\partial\vg(U_K,\K_2)$
is finite. We will later discuss when it is $p$-torsion free
in case that $K$ is a local field (Theorem B (1)).

\section{Theorem A and $\Phi(X_R,D_R)_{\Z_p}$}

\subsection{Tate curves over complete rings}
Let $A$ be a noetherian local ring which is complete
with respect to an ideal $qA$.
Then the {\it Tate curve} $E_q=E_{q,A[q^{-1}]}$
is defined to be the proper smooth scheme over $A[q^{-1}]$
defined by the equation
\begin{equation}\label{tate0}
y^2+xy=x^3+a_4(q)x+a_6(q)
\end{equation}
with
\begin{equation}\label{tate1}
a_4(q)=-5\sum_{n=1}^\infty \frac{n^3q^n}{1-q^n},\quad
a_6(q)=-\sum_{n=1}^\infty \frac{(5n^3+7n^5)q^n}{12(1-q^n)}.
\end{equation}
Let $O$ be the infinity point of $E_q$.
The series
\begin{equation}\label{xexp}
x(u)=\sum_{n\in\Z}\frac{q^nu}{(1-q^nu)^2}
-2\sum_{n\geq1}\frac{q^n}{(1-q^n)^2}
\end{equation}
\begin{equation}\label{yexp}
y(u)=\sum_{n\in\Z}\frac{(q^nu)^2}{(1-q^nu)^3}
+\sum_{n\geq1}\frac{q^n}{(1-q^n)^2}
\end{equation}
converge for all $u\in A[q^{-1}]^*-q^\Z$. They induce a 
homomorphism
\begin{equation}\label{thehom}
A[q^{-1}]^*/q^\Z\lra E_q(A[q^{-1}]),
\quad
u\longmapsto
\begin{cases}
(x(u),y(u)) & \text{if }u\not\in q^\Z\\
O & \text{if }u\in q^\Z,
\end{cases}
\end{equation}
which is bijective if $A$ is a complete discrete valuation ring.
We often use the uniformization \eqref{thehom}.
In particular, we have
\begin{equation}\label{canonicalinvariant}
\frac{dx}{2y+x}=\frac{du}{u}
\end{equation}
which we call the canonical invariant 1-form.
\begin{defn}[Theta function]\label{thetadefn}
$$
\theta(u)=\theta(u,q)
\os{{\text{def}}}{=}(1-u)\prod_{n=1}^\infty(1-q^nu)(1-q^nu^{-1})
$$
\end{defn}
$\theta(u)$ converges for all $u\in A[q^{-1}]^*$ and satisfy
\begin{equation}\label{eqA}
\theta(qu)=\theta(u^{-1})=-u^{-1}\theta(u).
\end{equation}
In particular a function
\begin{equation}\label{11}
f(u)=c\prod_i\frac{\theta(\alpha_iu)}{\theta(\beta_iu)}
\end{equation}
is $q$-periodic if $\prod_i{\alpha_i}/{\beta_i}=1$.
Suppose that $A=\cR$ is a complete discrete valuation ring.
We denote by $\cK$ its quotient field.
Then, for any rational function $f(u)$ on 
$E_{q,\ol{\cK}}$,
one can find $c,~\alpha_i,~\beta_i\in \ol{\cK}^*$ such that $f(u)$
is given as in \eqref{11}. Thus we have the one-one correspondence
\begin{equation}\label{cor}
\ol{\cK}(E_q)^*
\os{1:1}{\longleftrightarrow}
\left\{
c\prod_i\frac{\theta(\alpha_iu)}{\theta(\beta_iu)}
~;~
c,~\alpha_i,~\beta_i\in \ol{\cK}^*
\text{ with }\prod_i{\alpha_i}/{\beta_i}=1
\right\},
\end{equation}
by which we often identify the both sides.

\subsection{Statement of Theorem A}\label{seconleysect}
Let $p\geq 3$ be a prime number, and $r\geq1$ an integer.
Let $R_0$ be the integer ring of a finite unramified extension
$K_0$ of $\Q_p$ of degree $d$. 
Since $K_0/\Q_p$ is unramified, 
there is a cyclotomic basis
$\underline{\mu}=\{\zeta_1,\cdots,\zeta_d\}$
of $R_i$, namely $\zeta_i$ are roots of unity
such that $R_0=\bigoplus_{i=1}^d\Z_p\zeta_i$.
Letting $q_0$ be an indeterminate
we put
$q:=q_0^r$.
Let  
$E_q=E_{q,R_0((q_0))}$ be the Tate curve
over $R_0((q_0)):=R_0[[q_0]][q_0^{-1}]$.

The dlog map
$\K_2\to\Omega^2_{E_q/R_0}$ of Zariski
sheaves on $E_q$ gives rise to
\begin{equation}\label{dlogtate0}
\dlog:\vg(E_q,\K_2)\lra
\vg(E_q,\Omega^2_{E_q/R_0}).
\end{equation}
We also write by $\dlog$ the composite map
$$
K_2(E_q)\lra 
\vg(E_q,\K_2)\os{\dlog}{\lra}
\vg(E_{q},\Omega^2_{E_{q}/R_0})
\lra
R_0((q_0))\frac{dq_0}{q_0}
\frac{du}{u},
$$
where the last map is the natural map.
Define a map
\begin{equation}\label{naturalmap00}
\kp:R_0((q_0))\frac{dq_0}{q_0}\frac{du}{u}
\lra \prod_{k\geq 1}(\Z_p/k^2\Z_p)^{\op d}
\end{equation}
in the following way. (Note that the inner component of
the right hand side is zero
unless $p\vert k$.)
Express
\begin{align}
f(q_0)\frac{dq_0}{q_0}\frac{du}{u}
&=\left(\sum_{j\in\Z}c_jq_0^j\right)
\frac{dq_0}{q_0}\frac{du}{u}\\
&=\left(\sum_{j\leq 0}c_jq_0^j+
\sum_{k=1}^\infty\sum_{i=1}^d
a_k^{(i)}\frac{\zeta_iq_0^k}{1-\zeta_i q_0^k}\right)
\frac{dq_0}{q_0}\frac{du}{u}\label{exp00}
\end{align}
with $a^{(i)}_k\in \Z_p$, which are uniquely determined.
Then we define
$$
\kp:f(q_0)\frac{dq_0}{q_0}\frac{du}{u}
\longmapsto(\bar{a}_k^{(1)},
\cdots,\bar{a}_k^{(d)})_{k\geq1}
$$
where bars denote modulo $k^2\Z_p$.

\Th A.
{\it 
Suppose that $p$ is prime to $2r$.
Then the composition
\begin{equation}\label{aseq00}
K_2(E_q)\os{\dlog}{\lra}
R_0((q_0))\frac{dq_0}{q_0}\frac{du}{u}\os{\kp}{\lra}
\prod_{k\geq 1}(\Z_p/k^2\Z_p)^{\op d}
\end{equation}
is zero. }

\vskip 6pt

\begin{rem}\label{remafter}
We can replace $K_2(E_q)$
with $\vg(E_q,\K_2)$ in Theorem A if
$p \geq 5$.
In fact the map
$K_2(E_q)\ot\Z[1/6]\lra 
\vg(E_q,\K_2)\ot\Z[1/6]
$
is surjective as $\dim E_q=2$
(\cite{soulecanada} Th\'eor\`eme 4 iv)).
\end{rem}

Theorem A is the most important result in this paper.
It gives a criterion for 2-forms to be contained
in the dlog image.

We will prove Theorem A in \S \ref{tausect}
and \S \ref{pf2sect}.
The outline is as follows. We first construct the following
commutative diagram
$$
\xymatrix{
& K_2(E_q)\ar[d]^{\tau_\infty}\ar[rd]^{\tau_\infty^\et}
\ar[ld]_{\dlog} & \\
R_0((q_0))\frac{dq_0}{q_0}\frac{du}{u}
& R_0((q_0))^*\ot\Z[r^{-1}] \ar[l]^{\iota}\ar[r]_h & 
\left(R_0((q_0))[p^{-1}]\right)^*/p^\nu
\\
}$$
for $\nu\geq1$
where $h$ is the natural map from the ring homomorphism
$R_0((q_0))\to R_0((q_0))[p^{-1}]$ and
$
\iota$ is given by
$h\mapsto \frac{dh}{h}\frac{du}{u}.
$
The map $\tau^\et_\infty$ is defined from the \'etale regulator
and the weight exact sequence on the \'etale 
$H^1$ of Tate curve (cf. \eqref{wtseq}).
The map
$\tau_\infty$ is defined in the way of algebraic 
$K$-theory,
which turns out to be an analogue of the Beilinson regulator,
though our construction is completely different from Beilinson's one
(and it is rather simple).
Although the three maps $\dlog$,
$\tau_\infty$ and $\tau^\et_\infty$
are defined in different ways,
it turns out that they are compatible.
In particular we can reduce the problem to the characterization of
the image of $\tau_\infty^\et$.
It follows from the definition that the image of 
$\tau^\et_\infty$
is contained in the kernel of the map
$x\mapsto x\cup q$ of cup-product in \'etale cohomology groups
of
$R_0((q_0))[p^{-1}]$.
Thus
we get a vanishing 
\begin{equation}\label{vankato}
\{\tau_\infty(\xi),q\}=0
\quad\text{ in }
K^M_2(R_0 \langle q_0 \rangle)/p^\nu
\end{equation} for all $\xi\in K_2(E_q)$
where $R_0 \langle q_0 \rangle$ denotes the $p$-adic completion
of $R_0((q_0))$ (cf. \S \ref{fielddefsect}).
Finally we apply Kato's explicit reciprocity law to \eqref{vankato}
and obtain the desired result.

\begin{rem}\label{yasuda}
 In an earlier version,
I proved Theorem A only in the case of
$R_0=\Z_p$ and only for $k\leq p(2p-2)$.
The outline of the proof was the same, but I used the Artin-Hasse
formula instead of Kato's explicit reciprocity law.
Dr. Yasuda informed me Kato's work and pointed out
that we can get more general statement.
\end{rem}

The map $\kp$ depends on the basis $\underline{\mu}$.
However we have
\begin{lem}\label{dependrem}
The kernel of $\kp$ does not depend on the choice
of $\underline{\mu}$.
\end{lem}
\begin{pf}
Let $f(q_0)\in R_0[[q_0]]$ be expressed as in \eqref{exp00}.
Then 
\begin{equation}\label{dep1}
a_k^{(i)}\in k\Z_p\text{ for }\forall k,~i
\Longleftrightarrow
f(q_0)=q_0\frac{d}{dq_0}\log h(q_0),
\quad \exists h(q_0)\in R_0((q_0))^*
\end{equation}
where $\log(1-x)=-\sum_{i\geq 1}x^i/i$.
Note that $h(q_0)$ is uniquely determined up to constant.
Let $\sigma:R_0\to R_0$ be the Frobenius automorphism and 
$\varphi:R_0((q_0))\to R_0((q_0))$ the endomorphism
such that $\varphi(aq_0^i)=\sigma(a)q_0^{ip}$ for $a\in R_0$
(which is also
called the Frobenius).
Put
$l_\varphi(g)=p^{-1}\log(\varphi(g)/g^p)$ for $g\in R_0((q_0))^*$
(cf. \S \ref{step2}).
Let $\zeta\in R_0$ be an arbitrary root of unity.
Since $\sigma(\zeta)=\zeta^p$, we have
$$
l_\varphi(1-\zeta q_0^j)=\sum_{m\geq 1,(m,p)=1} 
\frac{(\zeta q_0^j)^m}{m},
\quad j\geq 1.
$$
Therefore if we express
$$
h(q_0)=c(-q_0)^m\prod_{k=1}^\infty\prod_{i=1}^d
(1-\zeta_i q_0^k)^{b^{(i)}_k}
$$
with $c\in R_0^*$, $m\in\Z$ and $-kb^{(i)}_k=a^{(i)}_k$,
then we have
$$
l_\varphi(h)=\sum_{k=1}^\infty\sum_{i=1}^d\sum_{m\geq 1,(m,p)=1} 
b_k^{(i)}\frac{(\zeta_i q_0^k)^m}{m}.
$$
One can easily show that 
\begin{align}
a_k^{(i)}\in k^2\Z_p\text{ for }\forall k,~i
&\Longleftrightarrow
b_k^{(i)}\in k\Z_p\text{ for }\forall k,~i\\
&\Longleftrightarrow\label{dep2}
l_\varphi(h)=q_0\frac{d}{dq_0}g(q_0),
\quad \exists g(q_0)\in R_0((q_0))
\end{align}
by the induction on $k$.
Since the right hand side of \eqref{dep1} and \eqref{dep2}
do not depend on $\underline{\mu}$, we see that so does not
the kernel of $\kp$.
\end{pf}

\subsection{$\Phi(X_R,D_R)_{\Z_p}$ : a bounding space of dlog image}
Let $p\geq 3$ be a prime number.
Let $K$ be a finite unramified extension of $\Q_p$
and $R$ its integer ring.
Let $\pi_R:X_R\to C_R$ be an elliptic surface over $R$
satisfying
{\bf (Rat)}.
We use the notations in \S \ref{setupsect}.
Assume that
$p$ is prime to $r_1\cdots r_s$.
Let $t_i\in \O_{C_R}$ be the uniformizer of 
$P_{i,R}$ for $1\leq i \leq s$.
Let
$
\iota_i:\Spec R((t_i))\lra S_R
$
be the punctured neighborhood of $P_{i,R}$
and $X_i$ the fiber product: 
$$
\begin{CD}
X_i@>>>U_R\\
@VVV@VV{\pi_R}V\\
\Spec R((t_i))@>{\iota_i}>> S_R.
\end{CD}
$$
Then $X_i$ is isomorphic to the Tate curve
(e.g. \cite{DeRa} VII Cor. 2.6).
More precisely 
let $q\in R((t_i))$ be the unique power series such that
$\ord_{t_i}(q)=r_i$ and 
\begin{equation}\label{jtate}
j(X_i)=\frac{1}{q}+744+196884q+\cdots.
\end{equation}
Then there is an isomorphism
\begin{equation}\label{tateisom0}
X_i\cong E_{q,R((t_i))}
\end{equation}
of $R((t_i))$-schemes (e.g. \cite{DeRa} VII. 2.6). 
The isomorphism \eqref{tateisom0} is unique
up to the translation and the involution $u\mapsto u^{-1}$.
Put $a_i:=t_i^{r_i}/j(X_i)\vert_{t_i=0}=t_i^{r_i}/q\vert_{t_i=0}\in R^*$.
Let $K_i=K(a_i^{1/r_i})$ and $R_i$ its integer ring.
Since $p$ is prime to $r_i$, $K_i$ is also unramified over $\Q_p$.
There is $q_i\in R_i((t_i))$ such that $q_i^{r_i}=q$
and we have $R_i((t_i))= R_i((q_i))$.
Put $\vg(X_R,\Omega^2_{X_R/R}(\log D_R))_{\Z_p}:=
\partial_\dR^{-1}(\Z_p^{\op s})$ where 
$\partial_\dR$ is the boundary map \eqref{phidr0}.
We have the composition of natural maps
\begin{multline}\label{drlocal1}
\vg(X_R,\Omega^2_{X_R/R}(\log D_R))_{\Z_p}\lra
\vg(X_i,\Omega^2_{X_i/R})\cong
\vg(E_{q,R((t_i))},\Omega^2_{E_{q,R((t_i))}/R})\\
\lra
R_i((q_i))\frac{dq_i}{q_i}\frac{du}{u}
\end{multline}
where the isomorphism in the middle is due to \eqref{tateisom0}.
Let $d_i=[K_i:\Q_p]$ and fix
a cyclotomic basis $\underline{\mu_i}$
of $R_i$.
There are the corresponding maps on $K$-cohomology to \eqref{drlocal1}
which give
rise to a commutative diagram
\begin{equation}\label{ccdia}
\begin{CD}
\vg(U_R,\K_2)\ot\Z_p@>{\dlog}>>
\vg(X_R,\Omega^2_{X_R/R}(\log D_R))_{\Z_p}\\
@VVV@VV{\eqref{drlocal1}}V\\
\vg(E_{q,R_i((q_i))},\K_2)\ot\Z_p@>{\dlog}>>
R_i((q_i))\frac{dq_i}{q_i}\frac{du}{u}\\
@.@VV{\phi_{\underline{\mu_i}}}V\\
@.\prod_{k\geq1}(\Z_p/k^2\Z_p)^{\op d_i}.
\end{CD}
\end{equation}
Let $\phi_i$ be the composition of the right vertical arrows.
\begin{defn}
$\Phi(X_R,D_R)_{\Z_p}:=\ker(\op\phi_i)=\bigcap_{i=1}^s\ker~\phi_i
\subset \vg(X_R,\Omega^2_{X_R/R}(\log D_R))_{\Z_p}$.
\end{defn}
It follows from Lemma \ref{dependrem} that
$\Phi(X_R,D_R)_{\Z_p}$ does not depend on the choice of the 
cyclotomic basis
$\underline{\mu_i}$ and hence depends only on $(X_R,D_R)$.
By Theorem A and Remark \ref{remafter} we have
$\phi_{\underline{\mu_i}}\dlog=0$ in \eqref{ccdia}
if $p\geq 5$ and
$(\phi_{\underline{\mu_i}}\dlog)\ot\Q=0$ if $p\geq3$.
Thus we have $\dlog\vg(U_R,\K_2)\ot\Z_p\subseteq 
\Phi(X_R,D_R)_{\Z_p}$ if $p\geq 5$
and $\dlog\vg(U_R,\K_2)\ot\Q_p\subseteq 
\Phi(X_R,D_R)_{\Z_p}\ot_{\Z_p}\Q_p$ if $p\geq 3$.

\medskip

Summarizing above together with Lemmas \ref{d=b} and \ref{d=bmore},
we have
\begin{thm}\label{key2}
Let $R$ be the ring of integer in an unramified extension $K$
of $\Q_p$, and
$\pi_R:X_R\to C_R$ an elliptic surface over $R$ satisfying
{\bf (Rat)}. 
Then we have
\begin{equation}\label{conj1seq}
\dlog\vg(U^0_K,\K_2)\ot_{\Z}\Z_p
\subseteq 
\Phi(X_R,D_R)_{\Z_p}
\quad\text{ if }
p\not\vert 6r_1\cdots r_s,
\end{equation}
\begin{equation}\label{conj1seq2}
\dlog\vg(U^0_{\ol{K}},\K_2)\ot_{\Z}\Q_p
\subseteq 
\Phi(X_R,D_R)_{\Z_p}\ot_{\Z_p}\Q_p
\quad\text{ if }
p\not\vert 2r_1\cdots r_s.
\end{equation}
\end{thm}
Although the sequence \eqref{aseq00} is far from exact,
I expect that
the equality in \eqref{conj1seq2} holds, 
namely the dlog image is characterized
by the numerical conditions on the Fourier coefficients:
\begin{conj}[characterization of dlog image]\label{conj1}
$$\dlog\vg(U^0_{\ol{K}},\K_2)\ot_\Z\Q_p=
\Phi(X_R,D_R)_{\Z_p}\ot_{\Z_p}\Q_p.$$
\end{conj}
We will see in \S \ref{modularsect} that it is true
for modular elliptic surfaces.



\def\rr{\cR\langle u \rangle}
\def\kk{\cK\langle u \rangle}
\def\rrd{\cR\langle u \rangle^{\flat}}
\def\kkd{\cK\langle u \rangle^{\flat}}
\section{Construction of $\langle~,~\rangle$ and $\tau_\infty$}
\label{tausect}
In this section we construct maps
\begin{equation}\label{first1}
\langle~,~\rangle:K_2^M(\kkd)\lra \cK^*,
\end{equation}
\begin{equation}\label{second2}
\tau_\infty:K_2(E_{q,A[q^{-1}]})\ot\Z[\frac{1}{r}]\lra 
A[q^{-1}]^*\ot\Z[\frac{1}{r}],
\end{equation}
where $\cK$ is a complete
discrete valuation field
(see \S \ref{fielddefsect} below for the definition
of $\kkd$), 
$(A,q,r)$ is as in the beginning of \S \ref{constrsect},
and
$E_{q,A[q^{-1}]}$ is the Tate curve over $A[q^{-1}]$.
Both maps play important roles in the proof of 
Theorem A.
I keep in mind the cases of $\cK=L$ a finite extension
of $\Q_p$ or $\cK=k((q_0))$ 
for \eqref{first1}
and the case that $A=R_0((q))$ or $R_0$ for \eqref{second2}.
\subsection{Construction of the pairing
$\langle~,~\rangle:K_2^M(\kkd)\to \cK^*$}
\label{sect2} 
Let $\cR$ be a complete
discrete valuation ring with a uniformizer $\pi_{\cK}$
and $\cK$ the quotient field of $\cR$. 
Let
$\ord_\cK:\cK\to \Z\cup\{+\infty\}$ be the
valuation such that $\ord_\cK(\pi_\cK)=1$. By convention,
$\ord_\cK(0)=+\infty$. This map is uniquely extended on $\ol{\cK}$ the
algebraic closure of $\cK$, which we also write by $\ord_\cK$. 
We put $\cR_n=\cR/\pi_\cK^n$,
$U_n=1+\pi_\cK^nR$ ($n\geq 1$) and $k=\cR/\pi_\cK\cR$ the residue field.

\subsubsection{$\kk$ and $\kkd$}\label{fielddefsect}
Letting $u$ be an indeterminate, we define a ring
$\rr$ to be the $\pi_{\cK}$-adic completion of 
$\cR((u)):=\cR[[u]][u^{-1}]$:
$$
\rr \os{\hbox{def}}{=} \plim{i}
\left(\cR/\pi_\cK^i((u))\right).
$$
$\rr$ is the subring of the ring
$\cR[[u,u^{-1}]]$ of formal Laurent series which is generated by
$$
\sum_{i=-\infty}^{+\infty}a_iu^i \quad (a_i\in \cR)
$$
such that $\ord_\cK(a_i)\to +\infty$ as $i\to -\infty$.
Since $\cR((u))$ is a PID, $\rr$ is a discrete valuation ring with a uniformizer $\pi_\cK$.
We denote by $\kk$ the fractional field :
$
\kk=\cK\ot_\cR \rr=\rr[\pi_\cK^{-1}].
$
The valuation on $\kk$ is given by
$$
\ord_{\kk}\left(\sum_{i=-\infty}^{+\infty}a_iu^i\right)
=\mathrm{min}(\ord_\cK(a_i);i\in\Z).
$$
Any $f\in \kk$ has the following expression
$$
f=a_0u^n\prod_{i=1}^\infty(1-a_{-i}u^{-i})(1-a_{i}u^{i}),
$$
with
\begin{equation}\label{exp1}
\begin{cases}
a_0\in \cK^*& \\
a_i\in \cR& i\not=0 \\
\ord_\cK(a_{-i})>0& i>0\\
\ord_\cK(a_{-i})\to
+\infty&\text{as }i\to +\infty.
\end{cases}
\end{equation}
Note $\ord_{\kk}(f)=\ord_\cK(a_0)$.
The residue field of $\kk$ is the field $k((u))$ 
of formal power series
with coefficient in $k$.
Note that
the field $\kk$ contains $\cK(u)$, however neither $\kk\subset \cK((u))$
nor $\kk\supset \cK((u))$.

\begin{lem}\label{flat}
Let $f\in \kk$.
The following are equivalent.
\begin{description}
\item[(1)]
For any $n\geq 1$, there are
$f_n\in \cK(u)$ and $\phi_n\in \pi_\cK^n\rr$ such that
$f=f_n(1+\phi_n)$.
\item[(2)]
For any $n\geq 1$, there are
$f'_n\in \cK(u)$ and $\phi'_n\in \pi_\cK^n\rr$ such that
$f=f'_n+\phi'_n$.
\end{description}
\end{lem}
\begin{pf}
Assume that (1) holds.
Put $N=\ord_{\kk}(f)$. We have $N=\ord_{\kk}(f_n)$ for all $n\geq 1$.
Let $n_1\geq
\mathrm{max}(n-N,1)$.
Put $f'_n=f_{n_1}$ and $\phi'_n=f_{n_1}\phi_{n_1}$.
Since $\ord_{\kk}(\phi_n)\geq n_1+N$, we have (2).
Conversely, assume that (2) holds. Then for all $n>N$ we have
$N=\ord_{\kk}(f'_n)$. Put $f_n=f'_{n+N}$ and $\phi_n=\phi'_{n+N}/f'_{n+N}$.
Since $\ord_{\kk}(\phi_n)\geq (n+N)-N=n$, we have (1).
\end{pf}

\begin{defn}
We define
$$
\kkd\subset\kk
$$
as the subset of the elements such that the equivalent conditions
in Lemma \ref{flat} hold.
Moreover, we define
$$
\rrd\os{\hbox{def}}{=}
\{f\in\kkd;\ord_{\kk}(f)\geq 0\}.
$$
\end{defn}
\begin{lem}\label{fielddef}
$\kkd$ is a complete discrete valuation field with $\pi_\cK$ a uniformizer,
and $\rrd$
is its valuation ring. The residue field is $k(u)$.
\end{lem}
\begin{pf}
If $f\in \kkd$ and $g\in (\kkd)^*$ then
$fg^{-1}\in\kkd$
by Lemma \ref{flat} (1) and
$f-g\in \kkd$ by Lemma \ref{flat} (2).
Therefore $\kkd$ is a field. The map $\ord_{\kk}$ gives
the discrete valuation on $\kkd$.
Then $\rrd$ is its valuation ring by definition.
It follows from Lemma \ref{flat} (2) that
any Cauchy sequence in $\rrd$ converges to an element in $\rrd$.
Thus $\rrd$ is complete. The last assertion is clear.
\end{pf}

\noindent{\bf Example.}
Let $q\in \cK^*$ such that $\ord_\cK(q)>0$. The typical example
of elements of $\kkd$ is the theta function
$\theta(u,q)=(1-u)\prod_{n\geq 1}(1-q^nu)(1-q^nu^{-1})$.
On the other hand, let $f=a_0+a_1u+a_2u^2+\cdots\in \cR[[u]]$
such that $f~\mathrm{mod}(\pi_\cK)\not\in k(u)$.
Then $f$ is contained in $\rr$ but not in $\rrd$.

\subsubsection{$\langle~,~\rangle_n$ 
and $\langle~,~\rangle$}\label{defsect0}
For $\alpha\in \ol{\cK}$
let $\tau_\alpha$ denote the tame symbol at $u=\alpha$:
$$
\tau_\alpha:K_2^M(\ol{\cK}(u))\lra \ol{\cK}^*,\quad
\{f,g\}\longmapsto
(-1)^{\ord_\alpha(f)\ord_\alpha(g)}
\left(
\frac{f^{\ord_\alpha(g)}}{g^{\ord_\alpha(f)}}
\right)\vert_{u=\alpha}.
$$
Define a pairing
\begin{equation}\label{pairing0}
(\kkd)^*\times (\kkd)^*\lra \cK^*/U_n, \quad (f,g)\longmapsto
\langle f,g\rangle_n
\os{\hbox{def}}{=}
\sum_{\ord_\cK(\alpha)>0}\tau_\alpha\{f_n,g_n\}
\end{equation}
where $f=f_n(1+\phi_n)$ and $g=g_n(1+\psi_n)$ are
as in Lemma \ref{flat} (1) and $\alpha$ runs over all
$\alpha\in \ol{\cK}$ such that $\ord_\cK(\alpha)>0$ (including $\alpha=0$).
Since the choices of $f_n$ or $g_n$ are not unique, we need to show
that the above pairing does not depend on them:

\begin{lem}\label{welldeflem}
The pairing \eqref{pairing0} is well-defined.
\end{lem}
\begin{pf}
We want to show that if $h\in \cK(u)^*\cap (1+\pi_\cK^n\rr)$
and $h'\in \cK(u)^*$ then we have
\begin{equation}\label{well0}
\sum_{\ord_\cK(\alpha)>0}\tau_\alpha\{h,h'\}\in U_n.
\end{equation}
More generally, we show
\begin{equation}\label{well1}
\sum_{\ord_\cK(\alpha)>0}\tau_\alpha\{h,a\}\in \ol{U}_n:=1+\pi_\cK^n\ol{\cR}
\end{equation}
\begin{equation}\label{well2}
\sum_{\ord_\cK(\alpha)>0}\tau_\alpha\{h,b-u\}\in \ol{U}_n
\end{equation}
for $a,b\in \ol{\cK}$
where $\ol{\cR}$ denotes the integer ring of $\ol{\cK}$.

\smallskip

For any $h\in \cK(u)^*$ there is a decomposition
\begin{equation}\label{decomp}
h=cu^n
\prod_{\ord_\cK(\alpha_i)>0}(1-\alpha_iu^{-1})^{r_i}
\prod_{\ord_\cK(\beta_j)\geq0}(1-\beta_ju)^{s_j}
\end{equation}
with $n,r_i,s_j\in \Z$ and
$c,~\alpha_i,~\beta_j\in \ol{\cK}^*$.
Since $h\in \cK(u)^*$ and $\cK$ is complete, both of
$$
\prod_{\ord_\cK(\alpha_i)>0}(1-\alpha_iu^{-1})^{r_i},
\quad
\prod_{\ord_\cK(\beta_j)\geq0}(1-\beta_ju)^{s_j}
$$
are in $\cK(u)^*$
and hence $c\in \cK^*$.
Assume $h\in 1+\pi_\cK^n\rr$. Then reduction mod $\pi_\cK$ shows
$n=0$. Moreover since
\begin{align*}
c\prod_{\ord_\cK(\alpha_i)>0}(1-\alpha_iu^{-1})^{r_i}
&=c(1+a_1u^{-1}+a_2u^{-2}+a_3u^{-3}+\cdots)\\
&\equiv \prod_{\ord_\cK(\beta_j)\geq0}(1-\beta_ju)^{-s_j}\mod \pi_\cK^n\rr\\
&=1+b_1u+b_2u^2+b_3u^3+\cdots,
\end{align*}
we have $c-1\equiv a_i\equiv b_j\equiv 0$ mod $\pi_\cK^n$.
Therefore it is enough to show \eqref{well1} and \eqref{well2}
in the following cases:

\begin{description}
\item[Case (i)]
$h=c$ where $c\equiv 1$ mod $\pi_\cK^n$,
\item[Case (ii)]
$h=\prod_{\ord_\cK(\alpha_i)>0}(1-\alpha_iu^{-1})^{r_i}
=1+a_1u^{-1}+a_2u^{-2}+\cdots$
with $\forall a_i\equiv 0$ mod $\pi_\cK^n$,
\item[Case (iii)]
$h=\prod_{\ord_\cK(\beta_j)\geq 0}(1-\beta_ju)^{s_j}
=1+b_1u+b_2u^2+\cdots$
with $\forall b_i\equiv 0$ mod $\pi_\cK^n$.
\end{description}

\noindent\underline{Proof : Case (i)}.
Easy.

\medskip

\noindent\underline{Proof : Case (ii)}.
\eqref{well1} is clear. We show \eqref{well2}.
If $\ord_\cK(b)\leq 0$, then
$$
\sum_{\ord_\cK(\alpha)>0}\tau_\alpha\{h,b-u\}
=\prod_i\left(\frac{b}{b-\alpha_i}\right)^{r_i}=
(1+a_1b^{-1}+a_2b^{-2}+\cdots)^{-1}\equiv 1\mod \pi_\cK^n.
$$
If $\ord_\cK(b)> 0$ and $b\not=0$, then
\begin{align*}
\sum_{\ord_\cK(\alpha)>0}\tau_\alpha\{h,b-u\}
&=
\sum_{\ord_\cK(\alpha)>0}\tau_\alpha(\{\prod_{\alpha_i\not=b}
(1-\alpha_iu^{-1})^{r_i},b-u\}+\{(1-bu^{-1})^r,b-u\})\\
&=\prod_{\alpha_i\not=b}
\left(\frac{b}{b-\alpha_i}\right)^{r_i}
\left(1-\alpha_ib^{-1}\right)^{r_i}\cdot \frac{b^r}{b^r}\\
&=1.
\end{align*}
If $b=0$, then
\begin{align*}
\sum_{\ord_\cK(\alpha)>0}\tau_\alpha\{h,-u\}
&=
\sum_{\ord_\cK(\alpha)>0}\tau_\alpha\{\prod_{\alpha_i\not=0}
(1-\alpha_iu^{-1})^{r_i},-u\}\\
&=\prod_{\alpha_i\not=0}
\left(\frac{-\alpha_i}{-\alpha_i}\right)^{r_i}\\
&=1.
\end{align*}

\medskip

\noindent\underline{Proof : Case (iii)}.
\eqref{well1} is clear. \eqref{well2} is also clear if $\ord_\cK(b)\leq 0$.
Suppose $\ord_\cK(b)>0$. Then
$$
\sum_{\ord_\cK(\alpha)>0}\tau_\alpha\{h,b-u\}
=\prod_{\ord_\cK(\beta_j)\geq 0}(1-\beta_jb)^{s_j}
=1+b_1b+b_2b^2+\cdots
\equiv 1 \mod \pi_\cK^n.
$$
This completes the proof.
\end{pf}
\begin{lem}\label{functlem}
\begin{description}
\item[(1)]
$\langle ~,~\rangle_n$ is biadditive.
\item[(2)]
Let $f,g\in (\kkd)^*$ and
$f\in 1+\pi_\cK^n\rr$. Then $\langle f,g\rangle_n=1$.
\item[(3)]
For $m\geq n$ the pairings $\langle ~,~\rangle_m$
and $\langle ~,~\rangle_n$ are compatible:
$$
\begin{CD}
(\kkd)^*\times (\kkd)^*@>{\langle ~,~\rangle_m}>> \cK^*/U_m\\
@|@VVV\\
(\kkd)^*\times (\kkd)^*@>{\langle ~,~\rangle_n}>> \cK^*/U_n
\end{CD}
$$
\item[(4)]
$\langle f,1-f\rangle_n=1$ for any $f\in (\kkd)^*$ such that $f\not=1$.
\end{description}
\end{lem}
\begin{pf}
(1), (2) and (3) are obvious. To show (4), let $f=f_m(1+\phi_m)$
($m\geq n$) be the decomposition
as in Lemma \ref{flat} (1). Put $N=\ord_{\kk}(f)=\ord_{\kk}(f_m)$ and
$M=\ord_{\kk}(1-f)$.
Then for $m> M-N$ the order of
$
1-f_m=1-f+f_m\phi_m
$ is equal to $M$. Therefore we have
$$
1-f=1-f_m-f_m\phi_m=(1-f_m)\left(
1-\frac{f_m\phi_m}{1-f_m}\right)=(1-f_m)(1-\psi_m)
$$
with $\ord_{\kk}(\psi_m)\geq m+N-M$. Take $m\geq \mathrm{max}(n,n+M-N)$.
Thus we have
$$
\langle f,1-f\rangle_n=\sum_{\ord_\cK(\alpha)>0}\tau_\alpha\{f_m,1-f_m\}=1.
$$
This completes the proof.
\end{pf}

Due to Lemma \ref{functlem}
the pairing $\langle ~,~\rangle_n$
induces a map on Milnor's $K_2$:
\begin{equation}\label{pairing01}
\langle ~,~\rangle_n:
K_2^M(\kkd)\lra \cK^*/U_n, \quad \{f,g\}\longmapsto
\langle f,g\rangle_n.
\end{equation}
Due to Lemma \ref{functlem} (3) it
forms the projective system.
We define
\begin{equation}\label{pairing02}
\langle ~,~\rangle=\plim{n}\langle ~,~\rangle_n
:K_2^M(\kkd)\lra \plim{n}\cK^*/U_n=\cK^*.
\end{equation}

\begin{lem}[Explicit description
of the pairing]\label{cclem}
Let $f,g\in \kkd$ be expressed as
$$
f=a_0u^n\prod_{i=1}^\infty(1-a_{-i}u^{-i})(1-a_{i}u^{i}),
\quad
g=b_0u^m\prod_{i=1}^\infty(1-b_{-i}u^{-i})(1-b_{i}u^{i}),
$$
where $a_i$ and $b_i$ are as in \eqref{exp1}.
Then we have
\begin{equation}\label{cc}
\langle f,g\rangle
=(-1)^{nm}\frac{a_0^m}{b_0^n}\prod_{i,j\geq 1}
\frac{(1-a_i^{j/(i,j)}b^{i/(i,j)}_{-j})^{(i,j)}}
{(1-a_{-i}^{j/(i,j)}b^{i/(i,j)}_{j})^{(i,j)}}
\end{equation}
where $(i,j)$ denotes the greatest common divisor of $i$ and $j$.
\end{lem}
\begin{pf}
Note that $\langle~,~\rangle_n$ annihilates the subgroup
$\{1+\pi_\cK^n\rr\cap\kkd,-\}$ by Lemma \ref{functlem} (2).
Now \eqref{cc} is straightforward from the definition.
\end{pf}
Let $\rrd_n=\rrd/\pi_\cK^n\rrd\hra \rr/\pi_\cK^n\rr=\cR_n((u))$
be the subring.
Then due to Lemma \ref{functlem} (2) (and Lemma \ref{cclem}),
our $\langle~,~\rangle$ also induces the following
map
\begin{equation}\label{ccflat}
\langle~,~\rangle:K_2^M(\rrd_n)\lra \cR_n^*.
\end{equation}
\begin{rem}
We will see that the pairing \eqref{ccflat} can be extended
on $K_2^M(\cR_n((u)))$ and the same formula as \eqref{cc} holds
(\S \ref{step3}).
However I do not know how to show it directly from \eqref{cc}.
\end{rem}


\subsubsection{Definition of
$\tau_\infty:K_2(E_{q,\cK})\to \cK^*$}\label{useful}
Let $E_q=E_{q,\cK}$ 
be the Tate curve over $\cK$.
The uniformization \eqref{thehom}
induces a natural inclusion
$\cK(E_q)\hookrightarrow \cK\langle u\rangle^\flat$
as $\theta(u)\in \cK\langle u\rangle^\flat$.
We define $\tau_\infty:K_2(E_{q,\cK})\to \cK^*$
as the composition of the maps
$$
K_2(E_{q,\cK})\lra 
K_2^M(\cK(E_q))\lra K_2^M(\cK\langle u\rangle^\flat
)\os{\langle~,~\rangle}{\lra} \cK^*.
$$
This gives the definition of \eqref{second2} when $A$ is a complete
discrete valuation ring.

\begin{rem}[Beilinson's regulator]\label{beirem2}
Suppose $\cK=\C((q))$ and $E_q=\C^*/q^\Z$ with $\vert q\vert\ll1$. 
Then one can easily see that
our $\tau_\infty$ is equal to Beilinson's regulator:
$$
\tau_{\infty}(\xi)
=\exp\left(\sum\frac{1}{2\pi i}\int_{\delta}\log f\frac{dg}{g}
-\log g(o)\frac{df}{f}\right),
\quad \xi=\sum\{f,g\}\in K_2(E_q)
$$
where $o$ is an initial point and $\delta\in\pi_1(E_q,o)$ 
is the vanishing cycle.
\end{rem}

\subsection{Construction of 
$\tau_\infty$ for arbitrary $A$}\label{constrsect}
Let $A$ be a regular local ring which is complete with
respect to the maximal ideal $\frak m$.
Let $q\in \frak m$ and $q\not=0$. 
We assume that $\O:=A/\sqrt{qA}$ is also regular.
Since $A$ is a UFD, the prime ideal $\sqrt{qA}$ is generated
by a prime element $\pi$. Let $r\geq 1$ be the integer such that
$\pi^rA=qA$.
Let $\hat{A}$ be the completion of $A_{(\pi)}$ with respect to $\pi$.
Then $\hat{A}$ is a complete discrete valuation ring with
the uniformizer $\pi$.
We have already constructed the map
\begin{equation}\label{kt0}
\tau_\infty:K_2(E_{q,\hat{A}[q^{-1}]})\lra \hat{A}[q^{-1}]^*
\end{equation}
in \S \ref{useful}.
There is the natural map $K_2(E_{q,A[q^{-1}]})
\to K_2(E_{q,\hat{A}[q^{-1}]})$.
In this section we show 
$\tau_\infty(K_2(E_{q,A[q^{-1}]})\ot\Z[1/r])
\subset A[q^{-1}]^*\ot\Z[1/r]$, which allows us to 
define
\begin{equation}\label{kt0a}
\tau_\infty:K_2(E_{q,A[q^{-1}]})\ot\Z[\frac{1}{r}]
\lra A[q^{-1}]^*\ot\Z[\frac{1}{r}].
\end{equation}
(Unfortunately I could not remove ``$\ot\Z[1/r]$".
However it does not matter for our purpose because we only handle
$K_2(-)\ot\Z_p$ and the prime number 
$p$ is supposed to be prime to $r$.)

\subsubsection{Regular model 
${\cal E}_{q,A}$
}\label{modelsect}
Let $X$ be the proper scheme over $A$ defined
by the equation \eqref{tate0}. 
Then $X\times_A\Spec A[q^{-1}]$ is isomorphic to $E_{q,A[q^{-1}]}$,
and $X$ is regular except at the locus $x=y=\pi=0$.
Taking $(r-1)$-times blowing ups along the locus, 
we get a regular scheme
$\E_{q,A}$ over $A$. It is proper over $A$.
Moreover 
$\E_{q,A}\times_A\Spec A[q^{-1}]$ is isomorphic to $E_{q,A[q^{-1}]}$
and the special fiber $D_\O:=\E_{q,A}\times_A\Spec\O$ 
is the standard N\'eron $r$-gon over $\O$. 
Note that
$\E_{q,A}$ is the generalized elliptic curve in the sense of
\cite{DeRa} II 1.12 (see also VII. 1 for rigid geometric 
construction of $\E_{q,A}$).

\medskip

Quillen's localization exact
sequence (cf. \cite{Sr} Prop.(5.15))
yields an exact sequence
\begin{equation}\label{locz1}
K_2(\E_{q,A})\lra K_2(E_{q,A[q^{-1}]}) 
\lra K'_1(D_\O) \os{i_*}{\lra} K_1(\E_{q,A})
\end{equation}
where $i:D_\O\hra \E_{q,A}$.
We first calculate $K'_1(D_\O)$.
Let $D^{\mathrm{reg}}_\O$ be the regular locus of $D_\O$.
It is the disjoint union of
$r$-copies of $\G_{m,\O}$.
Then we have an exact sequence
\begin{equation}\label{locz2}
K_2(\G_{m,\O})^{\op r}\os{\delta}{\lra} 
K_1(\O)^{\op r}\lra K'_1(D_\O) \lra K_1(\G_{m,\O})^{\op r}
\os{\delta'}{\lra}\Z^{\op r}.
\end{equation}
Note $K_1(\G_{m,\O})\cong \Z\op \O^*$ and $K_2(\G_{m,\O})\cong
K_2(\O)\op\O^*$ (loc.cit. Cor.(5.5)).
Then $\delta(K_2(\O)^{\op r})=0$ and the induced map
$$
\left(K_2(\G_{m,\O})/K_2(\O)\right)^{\op r}
\cong (\O^*)^{\op r}\lra K_1(\O)^{\op r}\cong (\O^*)^{\op r}
$$
is given by the matrix
\begin{equation}\label{smat}
\begin{pmatrix}
1&&&&-1\\
-1&1&&&\\
&-1&&&\\
&&\ddots\\
&&&1\\
&&&-1&1
\end{pmatrix}
\end{equation}
Therefore the cokernel of $\delta$ is isomorphic to $\O^*$.
Similarly, $\delta'(\O^{\op r})=0$ and the induced map
$
K_1(\G_{m,\O})/K_1(\O)^{\op r}=\Z^{\op r}\to\Z^{\op r}
$
is also given by \eqref{smat}. Thus we have $\ker\delta'\cong
\Z\op\O^{*\op r}$:
\begin{equation}\label{smat0}
0\lra \O^*\lra K'_1(D_\O)\lra \Z\op \O^{*\op r}\lra 0.
\end{equation}
Since the first map in \eqref{smat0} has the splitting by
the map $K'_1(D_\O)\to K_1(\O)$ induced from the structure
morphism, we have
\begin{equation}\label{smat00}
K'_1(D_\O)\cong 
\Z\op \O^{*\op r}\op\O^*.
\end{equation}
This also implies
\begin{equation}\label{smat1}
\bigoplus_{k=1}^r
K_1(D_\O^{(k)})\lra K'_1(D_\O)\lra \Z\lra 0
\end{equation}
where $D_\O^{(k)}$ are the irreducible components of $D_\O$.

Next we study the kernel of $i_*$.
Let us consider the composition of the following 4 maps
\begin{equation}\label{smat2}
\bigoplus_{k=1}^r K_1(D_\O^{(k)})\lra K'_1(D_\O)
\os{i_*}{\lra} K_1(\E_{q,A})
\os{i^*}{\lra} K_1(D_\O)
\lra 
\bigoplus_{k=1}^r K_1(D_\O^{(k)}).
\end{equation}
Note $D_\O^{(k)}\cong \P^1_\O$ and $K_1(\P^1_\O)\cong
\O^{*\op 2}$.
The composition \eqref{smat2} is given by the matrix
\begin{equation}\label{smat3}
\begin{pmatrix}
-2&1&&&1\\
1&-2&&&\\
&1&&&\\
&&\ddots&1\\
&&&-2&1\\
1&&&1&-2
\end{pmatrix}
\end{equation}
Therefore we have
$$
\text{kernel of \eqref{smat2}}=\{(\eta_1,\cdots,\eta_r)\in 
\bigoplus_{k=1}^r K_1(D_\O^{(k)})~\vert~r\eta_1=\cdots=r\eta_r\}.
$$
This shows 
\begin{equation}\label{smat4}
\ker~i_*\subset \Z\op M\op \O^*,\quad
M:=\{(c_1,\cdots,c_r)\in \O^{*\op r}~\vert~c_1^r=\cdots=c_r^r\}.
\end{equation}
under the identification \eqref{smat00}.

We see that the inclusion in \eqref{smat4}
is equal if we invert $r$.
Let us see the last component $\O^*$.
Let $p:\Spec A[q^{-1}]\to E_{q,A[q^{-1}]}$ be a rational point.
Then $p_*K_2(A[q^{-1}])\subset K_2(E_{q,A[q^{-1}]})$ 
is onto the component $\O^*$.
Next, we see the second component.
Let $f_E:E_{q,A[q^{-1}]}\to\Spec A[q^{-1}]$ be the
structure morphism. Then $f_E^*K_2(A[q^{-1}])
\subset K_2(E_{q,A[q^{-1}]})$ is onto the diagonal component of
$\O^{*\op r}$, which is equal to $M$ if we invert $r$.

Finally we see the first component $\Z$.
Note that the composition $K_2(E_{q,A[q^{-1}]})\to K'_1(D_\O)
\to \Z$ is the boundary map $\partial$
in \S \ref{boundarysect}.
We show that the image of $\partial$ contains $r\Z$.
Since $A$ is complete with respect to $q$,
there is the homomorphism $\Z[[q]]\to A$ and it
gives rise to a commutative diagram
$$
\begin{CD}
K_2(E_{q,\Z((q))}) @>{\partial}>>\Z\\
@VVV@VV{\text{mult. by }r}V\\
K_2(E_{q,A[q^{-1}]}) @>{\partial}>>\Z.
\end{CD}
$$
Therefore it is enough to show that the top arrow is surjective.
Let $0<a<b<c$ be integers.
Put $q_0=q^{1/c}$ and consider the embedding
$\Z((q))\hra \Z((q_0))$.
Let
$$
f(u):=\frac{\theta(q_0^au)^{c}}{\theta(u)^{c-a}\theta(qu)^a}
=(-u)^a\left(\frac{\theta(q_0^au)}{\theta(u)}
\right)^{c}
$$
and
$$
g(u):=\frac{\theta(q_0^bu)^{c}}{\theta(u)^{c-b}\theta(qu)^b}
=(-u)^b\left(\frac{\theta(q_0^b u)}{\theta(u)}
\right)^{c}
$$
the rational functions on $E_{q,\Z((q_0))}$.
More precisely,
let $Q_i$ be the $\Z((q_0))$-rational point of 
$E_{q,\Z((q_0))}$ corresponding to $u=q_0^{-i}$.
Then $f$ (resp. $g$) has a zero at $Q_a$ (resp. $Q_b$)
and a pole at $Q_0$ (resp. $Q_0$).
Put
\begin{equation}\label{ts0}
\xi^M_{(a,b,c)}:=\left\{
\frac{f(u)}{f(q_0^{-b})},
\frac{g(u)}{g(q_0^{-a})}
\right\}\in K_2^M(\O(E_{q,\Z((q_0))}-Q_*)),
\end{equation}
where $Q_*=Q_a+Q_b+Q_0$.
The Milnor $K_2$-symbol $\xi^M_{(a,b,c)}$ defines a 
Quillen $K_2$-symbol
$\xi^Q_{(a,b,c)}\in K_2(E_{q,\Z((q_0))}-Q_*)$.
Recall the localization exact sequence
$$
K_2(E_{q,\Z((q_0))})\lra
K_2(E_{q,\Z((q_0))}-Q_*)\os{\tau}{\lra}
K_1(Q_*)=K_1(\Z((q_0)))^{\op3}
$$
where $\tau$ is the tame symbol.
It is easy to see $\tau(\xi^Q_{(a,b,c)})=0$.
Therefore there is a lifting
$\xi'_{(a,b,c)}\in 
K_2(E_{q,\Z((q_0))})$.
We put 
\begin{equation}\label{eissymbol}
\xi_{(a,b,c)}:=N(\xi'_{(a,b,c)})\in K_2(E_{q,\Z((q))})
\end{equation}
where $N:K_2(E_{q,\Z((q_0))})\to K_2(E_{q,\Z((q))})$ 
denotes the norm map for $\Z((q_0))/\Z((q))$.
A direct calculation yields that 
$
\partial(\xi_{(a,b,c)}')
=a(b-a)(b-c)
$ (cf. \cite{tate2} Rem.5.5).
Since the diagram
$$
\begin{CD}
K_2(E_{q,\Z((q_0))}) @>{\partial}>>\Z\\
@V{N}VV@VV{=}V\\
K_2(E_{q,\Z((q))}) @>{\partial}>>\Z
\end{CD}
$$
is commutative, we have
$
\partial(\xi_{(a,b,c)})
=a(b-a)(b-c)
$.
Choose $(a,b,c)=(1,2,3)$. Thus we have the surjectivity of
$\partial:K_2(E_{q,\Z((q))}) \to\Z$.
\begin{prop}
Under the isomorphism
$K'_1(D_\O)\cong \Z\op\O^{*\op r}\op\O^*$ 
$$
\ker~i_*\subset \Z\op M\op \O^*
$$
where 
$M:=\{(c_1,\cdots,c_r)\in \O^{*\op r}~\vert~c_1^r=\cdots=c_r^r\}$.
The inclusion is equal if we tensor $\Z[1/r]$.
\end{prop}
More precisely, we can see from the above that
$(\Z\op M\op \O^*)/\ker~i_*$ is a 
subquotient of $(\Z/r\Z)^{\op r+1}$.
\begin{cor}\label{6rsurj}
Let $p:\Spec A[q^{-1}]\to E_{q,A[q^{-1}]}$ be a rational point
and $f_E:E_{q,A[q^{-1}]}\to\Spec A[q^{-1}]$ the
structure morphism.
Let $\xi_{(a,b,c)}$ be as in \eqref{eissymbol} and also
think of it being a $K_2$-symbol of $E_{q,A[q^{-1}]}$
by the natural homomorphism $\Z[[q]]\to A$.
Then 
\begin{equation}
\label{6rsurj1234}
K_2(\E_{q,A})\op 
\Z\cdot\xi_{(1,2,3)}\op
p_*K_2(A[q^{-1}])\op f^*_EK_2(A[q^{-1}])
\lra K_2(E_{q,A[q^{-1}]})
\end{equation}
is surjective tensoring with $\Z[1/r]$.
\end{cor}

\subsubsection{Boundary map in algebraic $K$-theory}\label{step3}
Let $S$ be a noetherian ring.
We recall the boundary map
\begin{equation}\label{ktb}
\partial:K_{i+1}(S((u)))\lra K_i(S)
\end{equation}
from K. Kato \cite{kato1} \S 2.
Let
\begin{align*}
{\cal P}_S&:\text{the category of projective 
$S$-modules of finite rank}\\
{\cal M}_S&:\text{the category of finite $S$-modules}\\
{\cal M}^f_S&:\text{the category of finite 
$S$-modules which have finite projective 
resolutions}.
\end{align*}
Moreover let ${\cal H}_{S[[u]]}$ be the category of finite
$S[[u]]$-modules $M$ such that
$M\ot_{S[[u]]}S((u))=0$ and there is a resolution
\begin{equation}\label{flatres}
0\lra E_0\lra E_1 \lra M\lra 0
\end{equation}
with $E_i\in {\cal P}_{S[[u]]}$.

By Quillen's resolution theorem
(\cite{Sr} Theorem (4.6))
$K_i(S)\os{\mathrm{def}}{=}K_i({\cal P}_S)\os{\cong}{\to} 
K_i({\cal M}_S^f)$.
By Quillen's localization theorem 
(\cite{Q2}, \cite{Sr} Theorem (9.1))
we have the exact sequence
\begin{equation}\label{qloc2}
\cdots\lra K_{i+1}(S[[u]])\lra K_{i+1}(S((u)))
\os{\partial'}{\lra} K_i({\cal H}_{S[[u]]})
\lra K_i(S[[u]])\lra \cdots.
\end{equation}
Let $M\in {\cal H}_{S[[u]]}$. Then we can think of $M$ 
being a finite $S$-module.
Since $S[[u]]$ is a flat $S$-module, \eqref{flatres} gives
a flat resolution of $S$-module. Therefore $M$ has a finite
projective resolution of $S$-modules. Thus we can think of
$M$ being an object of ${\cal M}^f_S$.
Let
\begin{equation}\label{iota0}
\iota:{\cal H}_{S[[u]]}
\lra {\cal M}^f_S
\end{equation} be the functor which sends $M$ to $M$
(as $S$-module).
It gives rise to a map $\iota_*:K_i({\cal H}_{S[[u]]})\to K_i(S)$.
Then we define $\partial=\iota_*\partial'$.

\medskip

Let us consider the Tate curve $E_{q,A[q^{-1}]}$.
Let $\E_{q,A}$ be as in \S \ref{modelsect}.
Put 
$A_n=A/q^nA$ for $n\geq 1$.
Then there are the morphisms of schemes
$$
\G_{m,A_n}\hra \E_{q,A}\times_{A} \Spec A_n
\lra \E_{q,A}\longleftarrow E_{q,A[q^{-1}]}.
$$
This give rise to 
$$
K_2(\E_{q,A})\lra
K_2(A_n[u,u^{-1}])\lra
K_2(A_n((u)))\os{\partial}{\lra}K_1(A_n)=A_n^*.
$$
Passing to the projective limit we have 
\begin{equation}\label{ktsymb}
\tau_\infty^\dag:K_2(\E_{q,A})\lra \plim{n}A_n^*=A^*.
\end{equation}

\begin{prop}\label{katocom0prop}
Suppose $A=\cR$ a complete discrete valuation ring.
Then the following diagram
\begin{equation}\label{katocom0}
\begin{CD}
K_2(\E_{q,\cR})@>{\tau_\infty^\dag}>>
\cR^*\\
@VVV@VV{\bigcap}V\\
K_2(E_{q,\cK})@>{\tau_\infty}>>
\cK^*\\
\end{CD}
\end{equation}
is commutative.
Here the bottom map $\tau_\infty$ is as in \S \ref{useful}.
\end{prop} 
\begin{cor}\label{katocom0propcor}
$\tau_\infty(K_2(\E_{q,A}))\subset A^*$.
\end{cor}
\begin{pf}
Apply Proposition \ref{katocom0prop} for $\cR=\hat{A}$.
\end{pf}

We prove Proposition \ref{katocom0prop}.
Note that the map 
$K_2(\cR_n[u,u^{-1}])\to
K_2(\cR_n((u)))$ factors through $K_2(\rrd_n)$
where $\rrd_n=\rrd/\pi_\cK^n\rrd$.
By a theorem of van der Kallen
$K_2(\rrd_n)$ is isomorphic to Milnor's
$K_2^M(\rrd_n)$ 
as the residue field
is infinite (\cite{vdK}, \cite{vdK2}).
Therefore, to show Proposition \ref{katocom0prop}
it is enough to show that the following diagram
\begin{equation}\label{katocom00}
\begin{CD}
K_2^M(\rrd_n)@>{\langle~,~\rangle}>>\cR_n\\
@VVV@VV{=}V\\
K_2(\cR_n((u)))@>{\partial}>>\cR_n
\end{CD}
\end{equation}
is commutative.
\begin{lem}\label{ktlem1}
Let $f\in \cR_n[[u]]^*$. Then $\partial\{f,u\}=f(0)$ where $\partial$
is as in \eqref{ktb}.
\end{lem}
\begin{pf}
There are commutative diagrams
$$
\begin{CD}
K_2(\cR_n((u)))@>{\partial}>> \cR_n^*\\
@AAA@AAA\\
K_2(\cR((u)))@>{\partial}>> \cR^*\\
@VVV@VVV\\
K_2(\cK((u)))@>{\partial}>> \cK^*.
\end{CD}
$$
Since $f$ comes from $\cR[[u]]^*$ it is enough to show
$\partial\{g,u\}=g(0)$ for any $g\in \cK[[u]]^*$.
This is well-known (e.g. \cite{Sr} Lemma (9.12)).
\end{pf}
\begin{lem}\label{ktlem2}
Let $f\in \cR_n[[u]]^*$. Let $a\in \cR_n$ be a nilpotent element
(i.e. $a\in \pi_\cK\cR_n$).
Then $\partial\{f,u-a\}=f(a)$.
\end{lem}
\begin{pf}
Since $a$ is nilpotent we have an isomorphism
\begin{equation}\label{iota1}
h:\cR_n[[u]]\lra \cR_n[[v]],\quad u\longmapsto v+a
\end{equation}
of $\cR_n$-algebras. The functor $\Phi_h:M\mapsto
M\ot_{\cR_n[[u]]}\cR_n[[v]]$
induces a commutative
diagram
$$
\xymatrix{
{\cal H}_{\cR_n[[u]]}\ar[rd]^{\iota}\ar[dd]_{\Phi_h}&\\
& {\cal M}^f_{\cR_n}\\
{\cal H}_{\cR_n[[v]]}\ar[ru]^{\iota}
&
}
$$where $\iota$ is as in \eqref{iota0}.
Thus we get a commutative diagram
\begin{equation}\label{iota2}
\xymatrix{
K_2(\cR_n((u)))\ar[r]^{\partial'}\ar[dd]_{h^*}&
K_2({\cal H}_{\cR_n[[u]]})\ar[rd]^{\iota_*}\ar[dd]_{\Phi_{h*}}&\\
&& \cR_n^*\\
K_2(\cR_n((u)))\ar[r]^{\partial'}&
K_2({\cal H}_{\cR_n[[v]]})\ar[ru]^{\iota_*}&
}
\end{equation}
and hence we have $\partial h^*=\partial$.
Now we show $\partial\{f,u-a\}=f(a)$.
Due to the commutative diagram \eqref{iota2}
we have
$\partial\{f,u-a\}=\partial\{f(v+a),v\}$.
By Lemma \ref{ktlem1} we have
$\partial\{f(v+a),v\}=f(a)$. This completes the proof.
\end{pf}
\begin{lem}\label{ktlem3}
Let $a\in \cR_n$ be a nilpotent element.
Then we have $\partial\{a-u,u\}=1$.
\end{lem}
\begin{pf}
Let $k=\cR/\pi_\cK$ be the residue field.
It follows from a commutative diagram
$$
\begin{CD}
K_2(\cR_n((u)))@>{\partial}>>
\cR_n^*\\
@VVV@VV{\mathrm{mod}~ q}V\\
K_2(k((u)))@>{\partial}>>
k^*
\end{CD}
$$
that we have $\partial\{a-u,u\}\equiv \partial\{-u,u\}=1$ mod $q$.
This means $
\partial\{a-u,u\}=1+\pi_\cK f
$
for some $f\in \cR_n$.
Therefore 
it is enough to show that there is $N\geq 1$ which is prime to 
$\mathrm{char}(k)$ such that $\partial\{a-u,u\}^N=1$.

Let $N\geq1$ be an integer which is prime to $\mathrm{char}(k)$
and such that $a^N=0$. Let $\zeta_N$ be a primitive $N$-th root
of unity.
We may suppose that $\cR$ contains $\zeta_N$
by replacing $\cR$ with $\cR[\zeta_N]$.
In the same way as
the diagram \eqref{iota2}, we can see that
the isomorphism
$$
v:\cR_n[[u]]\lra \cR_n[[u]],\quad u\longmapsto \zeta_N u
$$
induces a commutative diagram
\begin{equation}\label{iotav2}
\xymatrix{
K_2(\cR_n((u)))\ar[rd]^{\partial}\ar[dd]^{v^*}&\\
&\cR_n^*\\
K_2(\cR_n((u))).\ar[ru]^{\partial}&
}
\end{equation}
Thus we have
\begin{align*}
\partial\{a-u,u\}^{N^2}&=
\prod_{i=0}^{N-1}\partial\{a-\zeta_N^iu,\zeta^i_Nu\}^N\\
&=
\prod_{i=0}^{N-1}\partial\{a-\zeta_N^iu,u\}^N\\
&=
\partial\{a^N-u^N,u^N\}\\
&=
\partial\{-u^N,u^N\}\\
&=1.
\end{align*}
This completes the proof.
\end{pf}
Now we show that the diagram \eqref{katocom00} is commutative.
Due to Lemma \ref{cclem},
it is enough to show
\begin{equation}\label{katocom1}
\partial\{f,g\}
=(-1)^{nm}\frac{a_0^m}{b_0^n}\prod_{i,j\geq 1}
\frac{(1-a_i^{j/(i,j)}b^{i/(i,j)}_{-j})^{(i,j)}}
{(1-a_{-i}^{j/(i,j)}b^{i/(i,j)}_{j})^{(i,j)}}
\end{equation}
for
$$
f=a_0u^n\prod_{i=1}^\infty(1-a_{-i}u^{-i})(1-a_{i}u^{i}),
\quad
g=b_0u^m\prod_{i=1}^\infty(1-b_{-i}u^{-i})(1-b_{i}u^{i}),
$$
with
$$
\begin{cases}
a_i,b_i\in \cR_n^*& i\not=0 \\
a_{-i},b_{-i}\in \pi_\cK\cR_n& i>0\\
a_{-i}=b_{-i}=0& i\gg 0.
\end{cases}
$$
We can reduce it to check \eqref{katocom1}
in the following cases:
\begin{equation}\label{case1}
(f,g)=
\begin{cases}
(u,h)& h\in \cR_n[[u]]^*\\
(u,1-au^{-i})&a\in \pi_\cK\cR_n,~i>0\\
(1-au^{-i},h)&a\in \pi_\cK\cR_n,~i>0,~h
\in \cR_n[[u]]^*\\
(1-au^{-i},1-bu^{-j})&a,b\in \pi_\cK\cR_n,~i,j>0.
\end{cases}
\end{equation}
The first case immediately follows from Lemma \ref{ktlem1}.
To show the rest, we may replace
$\cR_n$ with $\cR[a^{1/m}]/(\pi_\cK^{n})$.
Thus we can assume that $a^{1/i},~b^{1/j}\in \cR_n$ and hence
may assume $i=j=1$.
Then the assertion follows from Lemmas \ref{ktlem1},
\ref{ktlem2} and \ref{ktlem3}.
This completes the proof of the commutativity of \eqref{katocom00}
and hence Proposition \ref{katocom0prop}.

\subsubsection{Proof of $\tau_\infty(
K_2(E_{q,A[q^{-1}]})\ot\Z[1/r])
\subset A[q^{-1}]^*\ot\Z[1/r]$.}
Due to Corollary \ref{6rsurj}, it is enough to show
the assertion for each component in \eqref{6rsurj1234}.
Corollary \ref{katocom0propcor} asserts 
$\tau_\infty(K_2(\E_{q,A}))
\subset A^*$. Moreover it is easy to see
$\tau_\infty(p_*K_2(A[q^{-1}]))=0$ and
$\tau_\infty(f_E^*K_2(A[q^{-1}]))=0$.
Therefore the rest of the proof is to show
$\tau_\infty(\xi_{(1,2,3)})\in A[q^{-1}]^*$.
However it follows from Lemma \ref{cclem} that we have
\begin{equation}\label{explicitxi}
\tau_\infty(\xi_{(a,b,c)})=
N_{A'/A}\tau_\infty(\xi^{\prime\prime}_{(a,b,c)})=
N_{A'/A}\left((-q_0)^{a(b-a)(b-c)}
\left(
\frac{M_b}{M_{b-a}M_a}\right)^c\right)
\end{equation}
where
$$
M_i:=\prod_{n=1}^\infty
\frac{(1-q_0^{nc-c+i})^{nc-c+i}}{(1-q_0^{nc-i})^{nc-i}},
\quad(1\leq i<c)
$$
and $N_{A'/A}$ is the norm map (cf. \cite{tate2} Cor.5.3).
Thus it is clearly contained in $A[q^{-1}]^*$.
This completes the proof.
\subsubsection{Functoriality of $\tau_\infty$}
Let $(A_i,q_i, \pi_i)$ ($i=1,2$)
be as in the beginning of \S \ref{constrsect}.
Let
$\phi:A_1\to A_2$ be a homomorphism of local rings
such that $\phi(q_1)=q_2$.
Let $q_iA_i=\pi_i^{r_i}A_i$ (then $r_1\vert r_2$).
\begin{prop}\label{comprop}
The following diagram
$$
\begin{CD}
K_2(E_{q_1,A_1[q^{-1}_1]})\ot\Z[1/r_1]
@>{\tau_\infty}>> A_1[q_1^{-1}]^*\ot\Z[1/r_1]\\
@V{\phi^*}VV@V{\phi^*}VV\\
K_2(E_{q_2,A_2[q_2^{-1}]})\ot\Z[1/r_2]
@>{\tau_\infty}>> A_2[q_2^{-1}]^*\ot\Z[1/r_2]
\end{CD}
$$
is commutative.
\end{prop}
\begin{pf}
Due to Corollary \ref{6rsurj}, it is enough to show
the assertion for $\xi_{(a,b,c)}$ 
and $K_2(\E_{q_1,A_1})$.
As for $\xi_{(a,b,c)}$, it 
directly follows from \eqref{explicitxi}.
As for $K_2(\E_{q_1,A_1})$, it follows from the functoriality
of $\tau_\infty^\dag$.
\end{pf}
\begin{rem}
Since $\langle~,~\rangle$ is clearly functorial,
Proposition \ref{comprop} is trivial when $A_i$ are complete
discrete valuation rings.
However the ring homomorphism $\phi$ does not necessarily
induce $\hat{A}_1\to\hat{A}_2$ in general.
Hence we can not reduce the proof to
the functoriality of $\langle~,~\rangle$.
\end{rem}

\subsection{Compatibility with dlog map}
\label{dlsect}
Let $\O$ be a regular local ring and $A=\O[[t]]$
the ring of formal power series with coefficients in $\O$.
Let $r\geq 1$ and $a\in A^*$ and put $q=at^r$.
We have constructed the map
$$
\tau_\infty:K_2(E_q)\lra \O((t))^*\ot\Z[1/r]
$$
for the Tate curve $E_q=E_{q,\O((t))}$ over
$\O((t))$.
\begin{prop}\label{xyprop}
Suppose that $r$ is invertible in the quotient field of $\O$.
Then the dlog map factors through $\tau_\infty$ as follows
$$\xymatrix{
&\O((t))^*\ot\Z[1/r]\ar[d]^{\iota}\\
K_2(E_q)\ot\Z[1/r]\ar[ur]^{\tau_\infty}\ar[r]_{\dlog}&
\O((t))\frac{dt}{t}\frac{du}{u}\ot\Z[1/r]
}$$
where the map $\iota$ is given as follows
$$
\iota:\O((t))^*
\lra
\O((t))\frac{dt}{t}\frac{du}{u},\quad
h\longmapsto \frac{dh}{h}\frac{du}{u}.
$$
\end{prop}
\begin{pf}
We may replace
$\O$ with its quotient field $F$.
Then $\cK=F((t))$ is a complete discrete valuation field
with the valuation ring $\cR=F[[t]]$.
It is enough to show that the following diagram
\begin{equation}\label{pcc}
\xymatrix{
&F((t))^*\ar[d]\\
K_2^M(\kkd)\ar[ur]^{\langle~,~\rangle}\ar[r]_{\dlog}&
F((t))\frac{dt}{t}\frac{du}{u}
}\end{equation}
is commutative where
the right vertical arrow is given by
$$
h\longmapsto \frac{dh}{h}\frac{du}{u}.
$$
Put
$$
\hat{\Omega}^1_{\cK/F}:=\left(
\plim{\nu}\Omega^1_{\cR/F}/t^\nu\Omega^1_{\cR/F}
\right)\ot_\cR\cK,
\quad
\hat{\Omega}^1_{\kkd/\cK}:=\left(
\plim{\nu}\Omega^1_{\rrd/\cR}/t^\nu\Omega^1_{\rrd/\cR}
\right)\ot_\cR\cK.
$$
Note$$
\frac{dx}{2y+x}=\frac{du}{u}\quad\text{ in }
\hat{\Omega}^1_{\kkd/\cK}.
$$
Let 
$$
\Res:\hat{\Omega}^1_{\kkd/\cK}\lra \cK
$$ be the residue map at $u=0$, namely if
we express $\omega=\sum_{n\in\Z}a_nu^ndu$ in the unique way then
$\Res(\omega)=a_{-1}$.
It is extended to a map 
\begin{equation}\label{res0}
\Omega^2_{\kkd/F}\lra
\hat{\Omega}^1_{\cK/F}\ot
\hat{\Omega}^1_{\kkd/\cK}
\os{\mathrm{id}\ot\Res}{\lra}
\hat{\Omega}^1_{\cK/F} .
\end{equation}
Then the commutativity of \eqref{pcc} is equivalent to 
that the following diagram
\begin{equation}\label{pcc1}
\begin{CD}
K_2^M(\kkd)@>{\langle~,~\rangle}>> \cK^*\\
@V{\dlog}VV@VV{\dlog}V\\
\Omega^2_{\kkd/\Q_p}@>{\eqref{res0}}>>
\hat{\Omega}^1_{\cK/\Q_p}.
\end{CD}
\end{equation}
is commutative.
However this follows from Lemma \ref{cclem}
and the fact that
$$
(u-\alpha)^{-1}=
\begin{cases}
u^{-1}\sum_{n=0}^\infty(\alpha u^{-1})^n & \ord_{\cK}
(\alpha)>0\\
-\alpha^{-1}\sum_{n=0}^\infty(\alpha^{-1} u)^n
& \ord_{\cK}(\alpha)\leq 0
\end{cases}
$$
in $\cK$.
This completes the proof of Proposition \ref{xyprop}.
\end{pf}

\subsection{Compatibility with \'etale regulator}
\label{etasect}
Let $A$ be a complete regular local ring and $q\in A$ a nonzero 
and not invertible element.
Let $n$ be a nonzero integer.
Let $f_E:E_{q,A[q^{-1},n^{-1}]}\to A[q^{-1},n^{-1}]$ 
be the Tate curve.
Then there is the weight exact sequence
\begin{equation}\label{wtseq}
0 \lra \Z/n(j)\lra
R^1f_{E*}\Z/n(j)\lra \Z/n(j-1)\lra 0
\end{equation}
of etale sheaves on $\Spec A[q^{-1},n^{-1}]$
(cf. \cite{DeRa} VII 1.13).
Recall the etale regulator map (also called the Chern class map)
\begin{equation}\label{chern}
c_{2,2}^\et:K_2(E_{q,A[q^{-1},n^{-1}]})
\lra H^2_\et(E_{q,A[q^{-1},n^{-1}]},\Z/n(2))
\end{equation}
to the \'etale cohomology group (\cite{gillet}, \cite{soule2}). 
The regulator \eqref{chern} together with the Leray spectral sequence
gives rise to a map
\begin{equation}\label{chern1}
\rho:K_2(E_{q,A[q^{-1},n^{-1}]})
\lra H^1_\et(A[q^{-1},n^{-1}],R^1f_{E*}\Z/n(2)).
\end{equation}
(cf. \cite{tate2} Lem.2.1).
We have from \eqref{wtseq} and \eqref{chern1} 
\begin{equation}\label{chern2}
\tau_\infty^\et:K_2(E_{q,A[q^{-1},n^{-1}]})
\lra H^1_\et(A[q^{-1},n^{-1}],\Z/n(1))\cong A[q^{-1},n^{-1}]^*/n,
\end{equation}
where ``$\cong$" follows from the fact that $A[q^{-1},n^{-1}]$
is a UFD.
\begin{prop}\label{etcomp}
Let $(A,q,\pi,r)$ be as in the beginning of \S \ref{constrsect}.
Suppose that $n$ is prime to $r$.
Then the following diagram
$$
\begin{CD}
K_2(E_{q,A[q^{-1}]})\ot\Z[1/r]
@>{\tau_\infty}>> A[q^{-1}]^*\ot\Z[1/r]\\
@VVV@VVV\\
K_2(E_{q,A[q^{-1},n^{-1}]})/n
@>{\tau^\et_\infty}>>A[q^{-1},n^{-1}]^*/n
\end{CD}
$$
is commutative
in the following cases.
\begin{enumerate}
\renewcommand{\theenumi}{(\roman{enumi})}
\item
$A$ is a complete discrete valuation ring.
\item
$A=R_0[[t]]$ where $R_0$ is an integer ring of a finite
extension $K_0$ of $\Q_p$ and $q=at^r$ with $a\in 
A^*$.
\end{enumerate}
\end{prop}

\subsubsection{Proof of (i)}
Let $A=\cR$ be a complete discrete valuation ring.
In this case we do not need ``$\ot\Z[1/r]$":
\begin{equation}\label{icase}
\begin{CD}
K_2(E_{q,\cK})
@>{\tau_\infty}>> \cK^*\\
@VVV@VVV\\
K_2(E_{q,\cK})/n
@>{\tau^\et_\infty}>>\cK^*/n
\end{CD}
\end{equation}
where $n$ is supposed to be invertible in $\cK$.

\medskip

To prove the commutativity of the diagram \eqref{icase}
we first extend $\tau_\infty^\et$
to a pairing $\langle~,~\rangle^\et$.
Let
\begin{equation}\label{etaud0}
K_2^M(\cK\langle u\rangle)\lra H^2_{\et}(\cK\langle u\rangle
,\Z/n(2))
\end{equation}
be the Galois symbol map.
Let $\left(K_2^M(
\cK\langle u\rangle
)\right)'$ be defined as
$$
0\lra \left(K_2^M(
\cK\langle u\rangle
)\right)'\lra
K_2^M(\cK\langle u\rangle)\lra H^2_{\et}(\cK\langle u\rangle
\ot_{\cK}\ol{\cK},\Z/n(2)).
$$
Then the Leray spectral sequence yields
\begin{equation}\label{etaud1}
\rho:
\left(K_2^M(
\cK\langle u\rangle
)\right)'
\lra H^1(\cK,H^1_\et(\cK\langle u\rangle
\ot_{\cK}\ol{\cK},\Z/n(2))).
\end{equation}
We define the natural map
\begin{equation}\label{etaud11}
H^1_\et(\cK\langle u\rangle
\ot_{\cK}\ol{\cK},\Z/n(j))
\cong(\cK\langle u\rangle
\ot_{\cK}\ol{\cK})^*/n\to \Z/n(j-1)
\end{equation} as
$$
f\longmapsto {\mathrm{Res}}\frac{df}{f}
$$
where $\mathrm{Res}$ denotes the residue map at $u=0$, namely if
we express $\omega=\sum_{n\in\Z}a_nu^ndu$ in the unique way then
${\mathrm{Res}}(\omega)=a_{-1}$.
The maps \eqref{etaud1} and \eqref{etaud11}
give rise to a pairing
\begin{equation}\label{surj1}
\langle~,~\rangle^\et
:\left(K_2^M(
\cK\langle u\rangle
)\right)'\lra \cK^*/n.
\end{equation}
(In \cite{tate2}, the pairing $\langle~,~\rangle^\et$ 
is written as $\hat{\tau}^\et_\infty$.)
By definition of $\cK\langle u\rangle^\flat$ there is a natural map
$K_2^M(\cK\langle u\rangle^\flat)\to 
\left(K_2^M(\cK\langle u\rangle)\right)'$.
We also write the 
composition $K_2^M(\cK\langle u\rangle^\flat)\to
\left(K_2^M(\cK\langle u\rangle)\right)'\to \cK^*/n$
by $\langle~,~\rangle^\et$.

It is clear from the construction that $\langle~,~\rangle^\et$
is compatible with $\tau_\infty^\et$.
On the other hand $\langle~,~\rangle$
is compatible with $\tau_\infty$.
Therefore it is enough to show
the following:
\begin{lem}\label{etcc}
The following diagram
$$
\begin{CD}
K_2^M(\cK\langle u\rangle^\flat)@>{\langle~,~\rangle}>> \cK^*\\
@VVV@VVV\\
K_2^M(\cK\langle u\rangle^\flat)/n
@>{\langle~,~\rangle^\et}>> \cK^*/n
\end{CD}
$$
is commutative.
\end{lem}
\begin{pf}
The pairing $\langle~,~\rangle$ is characterized by the 
following commutative
diagram
$$
\begin{CD}
K_2^M(\cK(u))@>{\sum_{\alpha}\tau_\alpha}>> \cK^*\\
@VVV@|\\
K_2^M(\cK\langle u\rangle^\flat)@>{\langle~,~\rangle}>> \cK^*.
\end{CD}
$$
The same commutative diagram holds for $\langle~,~\rangle^\et$
(\cite{tate2} Thm.4.4).
Thus the assertion follows.
\end{pf}
This completes the proof of Proposition \ref{etcomp} (i).

\subsubsection{Proof of (ii)}
We can reduce the proof of (ii) to the commutativity of
\eqref{icase}.
In fact, we may assume that $n=\l^\nu$ with $\l$ a prime number
(possibly $\l=p$).
Let $\pi_0$ be a uniformizer of $R_0$ and
$s_m:R_0[[t]]\to R_0$ the homomorphism given by $t\mapsto\pi^m$.
Then the map
$$
\prod_{ m}s_m:
\plim{\nu}A[q^{-1}]^*/\l^\nu \lra
\prod_{m}\plim{\nu}K_0^*/\l^\nu
$$
is injective where $m$ runs over all positive integers
which are prime to $\l$. 
Due to the functoriality of $\tau_\infty$
which follows from Proposition \ref{comprop} and
the functoriality of $\tau_\infty^\et$
which follows from the functoriality of etale regulator,
we can reduce the proof to the case $A=R_0$.

This completes the proof of Proposition \ref{etcomp} (ii).




\section{Proof of Theorem A}\label{pf2sect}
We are now in a position to prove Theorem A.

\medskip

Let $K_0$, $R_0=\bigoplus_{i=1}^d\Z_p\zeta_i$, $d=[K_0:\Q_p]$ 
and $R_0((q_0))$ with
$q=q_0^r$ be as in 
\S \ref{seconleysect}.
Let $E_q=E_{q,R_0((q_0))}$ be the Tate curve
and
\begin{equation}\label{sk1}
\dlog:K_2(E_q)\ot\Z_p\lra
R_0((q_0))\frac{dq_0}{q_0}\frac{du}{u}
\end{equation}
the dlog map.
Since $p$ is prime to $r$, it follows from 
Proposition \ref{xyprop} that
it factors through
$\tau_\infty$:
\begin{equation}\label{sk2}
\xymatrix{
&R_0((q_0))^*\ot\Z_p\ar[d]^{\iota}\\
K_2(E_q)\ot\Z_p\ar[ur]^{\tau_\infty}\ar[r]_{\dlog}&
R_0((q_0))\frac{dt}{t}\frac{du}{u}
}
\end{equation}
where 
$$
\iota:R_0((q_0))^*\ot\Z_p
\lra
R_0((q_0))\frac{dq_0}{q_0}\frac{du}{u},\quad
h\ot a\longmapsto a\frac{dh}{h}\frac{du}{u}.
$$
Let $\xi\in K_2(E_q)$ be an arbitrary element.
Express $h(q_0):=\tau_\infty(\xi)$ in the following way.
\begin{equation}\label{sk3}
h(q_0)=c(-q_0)^m\prod_{k=1}^\infty\prod_{i=1}^d
(1-\zeta_iq_0^k)^{b_k^{(i)}}
\end{equation}
with $c\in R_0^*$, $m\in\Z$ and $b_k^{(i)}\in \Z_p$.
Then 
\begin{align*}
\dlog(\xi)&=q_0\frac{d\log h(q_0)}{dq_0}
\cdot\frac{dq_0}{q_0}\frac{du}{u}\\
&=\left(m+
\sum_{k=1}^\infty\sum_{i=1}^d
(-kb^{(i)}_k)\frac{\zeta_iq_0^k}{1-\zeta_iq_0^k}
\right)
\frac{dq_0}{q_0}\frac{du}{u}
\end{align*}
due to the commutativity of \eqref{sk2}.
Therefore Theorem A is equivalent to say 
$b^{(i)}_k\in k\Z_p$ for all $k$ and $i$.
We will show it in the following steps:
\begin{description}
\item[(a0)] $c^r=1$. 
\item[(a1)] Let $R_0 \langle q_0 \rangle$ be
the $p$-adic completion of $R_0((q_0))$. Then
$$\{h(q_0),q_0\}=0\quad
\text{in }K_2^M(R_0 \langle q_0 \rangle)/p^\nu,
\quad\nu\geq 1.$$
\item[(a2)] Apply Kato's explicit reciprocity law for
$K_2^M(R_0 \langle q_0 \rangle)/p^\nu$.
\end{description}

\subsection{Proof of (a0)}
Embed $R_0((q_0)) \hra K_0((q_0))$. Then $\cK=K_0((q_0))$ is a complete
discrete valuation field with the valuation ring $\cR=K_0[[q_0]]$. 
We constructed the map
$\tau_\infty:K_2(E_{q,\cK})\to \cK^*$ in \S
\ref{useful}, which is compatible with the
\'etale regulator (see \eqref{icase}).
We thus have a commutative diagram
\begin{equation}\label{ckvan}
\xymatrix{
K_2(E_{q,\cK}) \ar[d]_{\rho}\ar[dr]^{\tau_\infty^\et}
\ar[r]^{\tau_\infty}&\cK^*\ar[d]\\
H^1_\et(\cK,
R^1f_{E*}\Z/n(2)
)\ar[r]&\cK^*/n\ar[r]^{s_q\quad}&K^M_2(\cK)/n
}
\end{equation}
for $n\geq1$
where the bottom is the exact sequence arising from the weight
exact sequence \eqref{wtseq} and the isomorphisms
$H^1_\et(\cK,\Z/n(1))\cong \cK^*/n$
and 
$H^2_\et(\cK,\Z/n(2))\cong K_2^M(\cK)/n$ (\cite{ms}).
Since the extension datum of $\eqref{wtseq}$ is $q$,
$s_q$ is the map given by $x\mapsto\{x,q\}$.
Therefore we have
$$
\{h(q_0),q\}=r\{h(q_0),q_0\}=0 \quad \text{in }K_2^M(\cK)/n.
$$
Applying the tame symbol $K_2^M(\cK)/n\to K_0^*/n$
we have $c^r=1$ in $K_0^*/n$. Since $n\geq 1$ is arbitrary,
we have $c^r=1$ in $K_0$.
This completes the proof.
\subsection{Proof of (a1)}
Let $K_0 \langle q_0 \rangle
=R_0 \langle q_0 \rangle[p^{-1}]$ be the quotient field.
In the same way as \eqref{ckvan} we have a commutative diagram
\begin{equation}\label{ckvan1}
\xymatrix{
K_2(E_{q})\ot\Z_p \ar[d]_{\rho}\ar[dr]^{\tau_\infty^\et}
\ar[r]^{\tau_\infty}&R_0((q_0))^*\ot\Z_p\ar[d]\\
H^1_\et({K_0 \langle q_0 \rangle},
R^1f_{E*}\Z/p^\nu(2)
)\ar[r]&{K_0 \langle q_0 \rangle}^*/p^\nu\ar[r]^{s_q\quad}
&K^M_2({K_0 \langle q_0 \rangle})/p^\nu
}
\end{equation}
and hence
\begin{equation}\label{ckvan2}
\{h(q_0),q\}=r\{h(q_0),q_0\}=0 \quad 
\text{in }K_2^M({K_0 \langle q_0 \rangle})/p^\nu.
\end{equation}
\begin{claim}\label{ckvan3}
The map $K_2^M(R_0 \langle q_0 \rangle)/p^\nu\to 
K_2^M({K_0 \langle q_0 \rangle})/p^\nu$ is injective.
\end{claim}
\begin{pf}
By a theorem of van der Kallen
$K^M_2(R_0 \langle q_0 \rangle)=K_2(R_0 \langle q_0 \rangle)$
as the residue field
is infinite (\cite{vdK}, \cite{vdK2}).
Recall the localization exact sequence
\begin{equation}\label{ckvan4}
K_2(\F_q((q_0)))\lra
K_2(R_0 \langle q_0 \rangle)\lra
K_2({K_0 \langle q_0 \rangle})\lra
K_1(\F_q((q_0)))=\F_q((q_0))^*.
\end{equation}
The right arrow is surjective as $R_0 \langle q_0 \rangle^*
\to \F_q((q_0))^*$ is 
surjective. 
Separate \eqref{ckvan4} as follows:
$$
K_2(\F_q((q_0)))\lra
K_2(R_0 \langle q_0 \rangle)\lra M\lra 0,
$$
$$
0\lra M \lra K_2({K_0 \langle q_0 \rangle})\lra
\F_q((q_0))^*\lra 0.
$$
Since the multiplication by $p^\nu$ on $\F_q((q_0))^*$ is injective,
we have $M/p^\nu\hra K_2({K_0 \langle q_0 \rangle})/p^\nu$.
On the other hand, since $K_2(\F_q((q_0)))/p^\nu=
K_2^M(\F_q((q_0)))/p^\nu=0$ (\cite{kato2} Lemma 5.6),
we have $K_2(R_0 \langle q_0 \rangle)/p^\nu
\hra M/p^\nu$. Thus we get 
$K^M_2(R_0 \langle q_0 \rangle)/p^\nu=
K_2(R_0 \langle q_0 \rangle)/p^\nu\hra 
K_2({K_0 \langle q_0 \rangle})/p^\nu=
K_2^M({K_0 \langle q_0 \rangle})/p^\nu$ 
as required.
\end{pf}
Now {\bf (a1)} follows from
\eqref{ckvan2} together with Claim \ref{ckvan3}

\subsection{Proof of (a2) : Kato's explicit reciprocity law}
\label{step2}
Since we work on $R_0((q_0))\ot\Z_p$, we may neglect $c$
due to {\bf (a0)}.

\smallskip

Let $\varphi:R_0 \langle q_0 \rangle\to R_0 \langle q_0 \rangle$ 
be the Frobenius such that
$\varphi(a\zeta^i q_0^j)=a\zeta^{ip}q_0^{jp}$
for $a\in \Z_p$ and $i,~j\in\Z$.
Let
$l_\varphi(f):=p^{-1}\log(\varphi(f)/f^p)$.
Then Kato's explicit reciprocity law describes
the syntomic regulator
$$
\theta_\varphi:K_2^M(R_0 \langle q_0 \rangle)/p^\nu 
\lra\Omega^1_{R_0 \langle q_0 \rangle}
/d{R_0 \langle q_0 \rangle}\ot\Z/p^\nu,\quad \nu\geq1
$$
for $p\geq 3$
explicitly (\cite{kato2} Cor. 2.9)
as follows:
$$
\theta_\varphi(\{a,b\})=l_\varphi(a)\frac{1}{p}\varphi
\left(\frac{db}{b}\right)
-l_\varphi(b)\frac{da}{a}.
$$
Due to {\bf (a1)}, we have
\begin{equation}\label{ckvan5}
\theta_\varphi(\{h,q_0\})=l_\varphi(h)\frac{dq_0}{q_0}=0,
\end{equation}
and hence
\begin{equation}\label{ckvan5b}
l_\varphi(h)=q_0\frac{dg}{dq_0}, \quad
\exists g=\sum_kc_kq_0^k\in R_0\langle q_0\rangle
\end{equation}
in ${R_0 \langle q_0 \rangle}/p^\nu=R_0((q_0))/p^\nu$
for all $\nu\geq1$.
We have from \eqref{sk3} and \eqref{ckvan5b} that
\begin{align*}
l_\varphi(h)&=\sum_{k=1}^\infty\sum_{i=1}^d\sum_{m\geq 1,(m,p)=1} 
b^{(i)}_k\frac{(\zeta_iq_0^k)^m}{m}\\
&=\sum kc_kq_0^k 
\end{align*}
in ${R_0 \langle q_0 \rangle}/p^\nu=R_0((q_0))/p^\nu$
for all $\nu\geq1$.
One can easily show
$b^{(i)}_k\in k\Z_p$ for all $k$ and $i$ by the induction on $k$.
This completes the proof.


\section{Theorem B and $\Phi(X_R,D_R)_{\F_p}$}
We use the notations in \S \ref{setupsect}.

\Th B.{\it
Let $F$ be a field of characteristic zero and
and $\pi_F:X_F\to C_F$ an elliptic surface over $F$
satisfying
{\bf (Rat)}.
Let $p$ be a prime number.
Put $X_{\ol{F}}=X_F\times_F\ol{F}$ etc.
\begin{enumerate}
\renewcommand{\labelenumi}{$(\theenumi)$}
\item\label{key11}
Assume that for any finite $p$-torsion $G_F$-module $M$,
$H^i_\et(F, M)$ is finite for all $i\geq 0$
 where $G_F$
denotes the absolute Galois group of $F$
(e.g. $F$ is a local field).
Then we have
$$\partial\vg(U^0_F,\K_2)\ot\Z_p
=\partial\vg(U^0_{\ol{F}},\K_2)\ot\Z_p
\subset\Z_p^{\op s}$$
if $p$ is prime to $6r_1\cdots r_s$ and 
$H_\et^3(X_{\ol{F}},\Z_p)$ is torsion free.
\item\label{key12}
The quotient
$$
\Z_p^{\op s}/(\partial\vg(U^0_{\ol{F}},\K_2)\ot\Z_p)
$$
is torsion free if
$p$ is prime to $6r_1\cdots r_s$,
$H_\et^3(X_{\ol{F}},\Z_p)
$ is torsion free and
$$
H^2_\et(X_{\ol{F}},\Z/p(2))^{G_F}=0.
$$
\end{enumerate}
}

\medskip

Note that there is the natural isomorphism
$H_\et^3(X_{\ol{F}},\Z_p)\cong
H_B^3(X_\C,\Z)\ot\Z_p
$
for a fixed embedding $\ol{F}\hra\C$
due to a theorem of M. Artin and the universal coefficient
theorem.
Therefore $H_\et^3(X_{\ol{F}},\Z_p)$ is torsion free for almost
all $p$.

The proof of Theorem B will be given in 
\S \ref{suslinsecct}
where Suslin's exact sequence (Theorem \ref{universalc0})
plays an essential role
(for example Claim \ref{c2} below is one of the key step).

\subsection{$\Phi(X_R,D_R)_{\F_p}$ : 
an upper bound of the rank of the dlog image 
}
\label{giving}
Admitting Theorems A and B, we give an upper bound
of the rank of the dlog image 
of $K_2$ of elliptic surface minus singular fibers
over an algebraically closed field of characteristic zero.

\medskip

Let $R$ be the ring of integer in an unramified extension $K$
of $\Q_p$, and
$\pi_R:X_R\to C_R$ an elliptic surface over $R$ satisfying
{\bf (Rat)}.

\begin{defn}
Let $f:\Z_p^{\op s}\to\F_p^{\op s}$ be the reduction modulo $p$.
Then we put 
$$\Phi(X_R,D_R)_{\F_p}:=f(\partial_\dR(\Phi(X_R,D_R)_{\Z_p}))
\subset
\F_p^{\op s}.$$
\end{defn}
\begin{thm}\label{key3}
Put $X_K:=X_R\times_RK$, $X_{\ol{K}}:=X_R\times_R\ol{K}$ etc.
Assume the following conditions (1) and (2).
\begin{enumerate}
\renewcommand{\labelenumi}{$(\theenumi)$}
\item
$p\not|6r_1\cdots r_s$ and $H^3_\et(X_{\ol{K}},\Z_p)$ is 
torsion free.
\item
There is a finitely generated subfield $F\subset \ol{K}$
such that $\pi_{\ol{K}}:X_{\ol{K}}\to C_{\ol{K}}$ is 
defined over $F$, which we write $\pi_F:X_F\to C_F$,
and
it satisfies
{\bf (Rat)}
and 
$H^2_\et(X_{\ol{F}},\Z/p(2))^{G_F}=0$.
\end{enumerate}
Then we have
\begin{equation}\label{conj1seq3}
\dlog\vg(U^0_{\ol{K}},\K_2)\ot\Z_p
\subseteq 
\Phi(X_R,D_R)_{\Z_p},
\end{equation}
\begin{equation}\label{conj1seq4}
\rank_\Z\dlog\vg(U^0_{\ol{K}},\K_2)\leq
\dim_{\F_p}\Phi(X_R,D_R)_{\F_p}.
\end{equation}
\end{thm}
\begin{pf}
\eqref{conj1seq3} follows from Theorem \ref{key2}
\eqref{conj1seq} and Theorem B (1).
We show \eqref{conj1seq4}.
Applying $\partial_\dR$ on \eqref{conj1seq3} we have
\begin{equation}\label{a+b3}
\partial\vg(U^0_{\ol{K}},\K_2)\ot\Z_p
\subset\partial_\dR(\Phi(X_R,D_R)_{\Z_p}).
\end{equation}
Applying $f$ to \eqref{a+b3}, we have
\begin{equation}\label{a+b4}
\dim_{\F_p}f(\partial\vg(U^0_{\ol{K}},\K_2)\ot\Z_p)
\leq \dim_{\F_p}\Phi(X_R,D_R)_{\F_p}.
\end{equation}
Due to Theorem B (2), $\partial\vg(U^0_{\ol{K}},\K_2)\ot\Z_p$
is a direct summand of $\Z_p^{\op s}$.
Therefore
\begin{equation}\label{a+b5}
\rank_{\Z_p}\partial\vg(U^0_R,\K_2)\ot\Z_p
=
\dim_{\F_p}f\partial(\vg(U^0_R,\K_2)\ot\Z_p).
\end{equation}
Then \eqref{conj1seq4} follows from
\eqref{a+b4} and \eqref{a+b5}.
\end{pf}
$\dim_{\F_p}\Phi(X_R,D_R)_{\F_p}$
is the desired bound. It is computable in many cases.
See \S \ref{expmsect}.

\begin{conj}\label{conj2}
Under the assumptions in Theorem \ref{key3},
$$
\rank_\Z\dlog\vg(U^0_{\ol{K}},\K_2)=
\dim_{\F_p}\Phi(X_R,D_R)_{\F_p}.$$
\end{conj}
We show that Conjecture \ref{conj1} implies the above.
In fact, Conjecture \ref{conj1} implies
$\partial\vg(U^0_{\ol{K}},\K_2)\ot\Q_p=
\partial_{\dR}\Phi(X_R,D_R)_{\Z_p}\ot\Q_p$ and hence
$\partial_{\dR}\Phi(X_R,D_R)_{\Z_p}/
\partial\vg(U^0_{\ol{K}},\K_2)\ot\Z_p$
is finite. However due to Theorem B (2), it must be 
torsion free, which
means $\partial\vg(U^0_{\ol{K}},\K_2)\ot\Z_p=
\partial_{\dR}\Phi(X_R,D_R)_{\Z_p}$.
Due to Theorem B (2) we have
Conjecture \ref{conj2}.


\section{Proof of Theorem B}\label{suslinsecct}

\subsection{Some Computations of cohomology groups}\label{cohsect}
Let 
$\pi:V\to C$ be a minimal elliptic surface over $\C$
with a section $e:C\to V$.
We omit to write the subscription ``$\C$"
for the simplicity of the notations.
Let $T=\{P_i\}\subset C$ be the set of points such that
$\pi^{-1}(P)$ is a singular fiber if and only of $P\in T$.
Put $Y_i:=\pi^{-1}(P_i)$, $Y=\sum Y_i$, $S^0=C-T$ and $U^0=V-Y$.
We assume that there is at least one
multiplicative fiber.
$H^\bullet(V)$ (resp. $H_\bullet(V)$) denotes the Betti (singular)
cohomology group (resp. homology group).

\bigskip

Let $n\geq 1$ be an integer.
There is the long exact sequence
\begin{equation}\label{localexseqcohomology}
\cdots\to H^2_Y(V,\Z/n(j))\to H^2(V,\Z/n(j))\to
H^2(U^0,\Z/n(j))\to \cdots.
\end{equation}
Let $Y_i=\cup_k Y_i^{(k)}$ be the irreducible
decomposition.
There is the isomorphism
$H^2_Y(V,\Z/n(1))\cong \bigoplus_{i,k} \Z/n\cdot y_i^{(k)}$
in which a base $y_i^{(k)}$ corresponds to the 
component $Y^{(k)}_i$.
Denote by $[-]$ the cycle class in the cohomology group
of $V$.
Then the map $H^2_Y(V,\Z/n(1))\to H^2(V,\Z/n(1))$
is given by $y_i^{(k)}\mapsto [Y_i^{(k)}]$.
Let $H^2(V,\Z/n(1))\to 
\bigoplus_{i,k}H^2(Y_i^{(k)},\Z/n(1))$
and
$H^2(V,\Z/n(1))\to H^2(e(C),\Z/n(1))$
be the pull-back maps. Then the composition 
\begin{multline}\label{intersectionmap}
\bigoplus_{i,k} \Z/n\cdot y_i^{(k)}
=H^2_Y(V,\Z/n(1))
\lra H^2(V,\Z/n(1))\\
\lra\bigoplus_{i,k}H^2(Y_i^{(k)},\Z/n(1))\op H^2(e(C),\Z/n(1))\\
=\bigoplus_{i,k} \Z/n\cdot y_i^{(k)}\op \Z/n\cdot e
\end{multline}
is given as follows
\begin{equation}\label{splitting}
y_i^{(k)}\longmapsto \sum_{i',k'}(Y_i^{(k)},Y_{i'}^{(k')})
\cdot y_{i'}^{(k')}+(Y_i^{(k)},e(C))\cdot e
\end{equation}
where $(-,-)$ denotes the intersection numbers.
\begin{lem}\label{intlem}
Let $Y_i=\pi^{-1}(P_i)=\sum_{k\geq 1} m_i^{(k)}\cdot Y_i^{(k)}$ be the
scheme theoretic fiber. We put
$$
[Y_i]:=\sum_km_i^{(k)}\cdot y_i^{(k)}.
$$
Suppose that $n$ is prime to $6r_1\cdots r_s$.
Then the composition \eqref{intersectionmap} induces
an isomorphism
$$
\bigoplus_{i,k\geq1} \Z/n\cdot y_i^{(k)}/\langle [Y_i]-[Y_j];i<j\rangle
\os{\cong}{\lra}
\bigoplus_{i\geq1,k\geq2}\Z/n\cdot y_i^{(k)}\op \Z/n\cdot e
$$
where $\langle [Y_i]-[Y_j];i<j\rangle$ is the subgroup
generated by $[Y_i]-[Y_j]$ for all $i<j$. 
\end{lem}
\begin{pf}
Since the cardinality of the both sides is the same, it is enough
to show the surjectivity.
To do this, it is enough to see that the map
\begin{equation}\label{intersectionmatrix}
\bigoplus_{k\geq1} \Z/n\cdot y_i^{(k)}
\lra
\bigoplus_{k\geq2}\Z/n\cdot y_i^{(k)}\op \Z/n\cdot e
\end{equation}
is surjective (and hence bijective) for each $i$.
We can check it as the case may be.
For example, if $Y_i$ is of type $I_{r_i}$,
then the matrix of \eqref{intersectionmatrix} is given by
$$
\begin{pmatrix}
-2&1&&&&&\\
1&-2&1&&&&\\
&1&-2&&&&\\
&&&\ddots&&&&\\
&&&&-2&1&0\\
&&&&1&-2&0\\
1&&&&0&1&1
\end{pmatrix}.
$$
The determinant of it is equal to $(-1)^{r_i-1}r_i$.
It is prime to $n$. 
Therefore \eqref{intersectionmatrix} is bijective.
The proofs for the other types are similar.
\end{pf}

\begin{lem}\label{hodge2}
Let $V_t=\pi^{-1}(t)$ be a smooth fiber.
Then
$H_1(V_t,\Q)\to H_1(V,\Q)$ is zero.
\end{lem}
\begin{pf}
Recall the assumption that there is a multiplicative fiber,
say $Y_1$. Let $\sigma$
be the local monodromy around $P_1$.
Then the action of $\sigma$ on $H^1(V_t,\Z)$ is given by the matrix
$$
\begin{pmatrix}
1&r\\
0&1
\end{pmatrix}
,\quad r\geq 1.
$$
Therefore the image of $H^1(V,\Q)\to H^1(V_t,\Q)$ has dimension $\leq1$.
On the other hand 
it constitutes a sub Hodge structure of $H^1(V_t,\Q)$ (\cite{d2} II).
Thus it is zero.
\end{pf}
\begin{prop}\label{equiporp}
Suppose that
$n$ is prime to $6r_1\cdots r_s$.
Then the following are equivalent.
\begin{enumerate}
\renewcommand{\labelenumi}{$(\theenumi)$}
\item\label{lemc1}
$n$ is prime to $\sharp H_1(V,\Z)_{\mathrm{tor}}$.
\item\label{lemc2}
$H^1(V_t,\Z/n)^{\pi_1(S^0,t)}=0$.
\item\label{lemc3}
$H^1(S^0,\Z/n)\cong H^1(U^0,\Z/n)$
\item\label{lemc3x}
$H^1(C,\Z/n)\cong H^1(V,\Z/n)$
\item\label{lemc4}
$H_1(V_t,\Z/n)\to H_1(V,\Z/n)$ is zero.
\end{enumerate}
\end{prop}
\begin{pf}
Note that $H_1(V,\Z/n)=H_1(V,\Z)/n$ for any CW complex $V$.
By the universal coefficient theorem and the fact that
$H_1(S^0,\Z)$ and $H_1(C,\Z)$ are torsion free, we have
$$
\eqref{lemc3}\Longleftrightarrow
H_1(S^0,\Z/n)\cong H_1(U^0,\Z/n)
$$
$$
\eqref{lemc3x}\Longleftrightarrow
H_1(C,\Z/n)\cong H_1(V,\Z/n).
$$

\noindent\eqref{lemc2}$\Longleftrightarrow$\eqref{lemc3}.
Let 
$$
E_2^{pq}=H^p(S^0,R^q\pi_*\Z/n)\Longrightarrow H^{p+q}(U^0,\Z/n)
$$
be the Leray spectral sequence. Due to the section $e$,
the maps $E_2^{p0}\to E^p$ are injective. 
Therefore we have an exact sequence
$$
0\lra H^1(S^0,\Z/n)\lra H^1(U^0,\Z/n)\lra H^1(V_t,\Z/n)^{\pi_1(S^0,t)}
\lra 0
$$
and hence the desired equivalence.

\noindent\eqref{lemc3}$\Longleftrightarrow$\eqref{lemc3x}.
In the following commutative diagram
$$
\begin{CD}
0@>>>H^1(S^0,\Z/n)@>>> H^1(U^0,\Z/n)
@>>> H^1(V_t,\Z/n)^{\pi_1(S^0,t)}
@>>> 0\\
@.@A{f_1}AA@A{f_2}AA@A{f_3}AA\\
0@>>>H^1(C,\Z/n)@>>> H^1(V,\Z/n)
@>>> \Coker
@>>> 0
\end{CD}
$$
the rows are split exact sequences.
We want to show that $f_3$ is bijective.
To do this, 
it is enough to show that $\Coker~f_1 \to \Coker~f_2$ is surjective
as $f_1$ and $f_2$ are injective.
Let the notations be as in Lemma \ref{intlem}.
There are the isomorphisms
$$
\Coker~f_1\cong \ker(\bigoplus_i\Z/n\cdot p_i\lra H^2(C,\Z/n)),
$$
$$
\Coker~f_2\cong \ker(\bigoplus_{i,k}\Z/n\cdot 
y_i^{(k)}\lra H^2(V,\Z/n))
$$
where $p_i$ corresponds to $P_i$.
The map $\Coker~f_1 \to \Coker~f_2$
is given by
$$
p_i\longmapsto [Y_i]:=\sum_km_i^{(k)}\cdot y_i^{(k)}.
$$
By Lemma \ref{intlem}, $\Coker~f_2$ is generated by $[Y_i]-[Y_j]$
($i<j$).
Therefore $\Coker~f_1 \to \Coker~f_2$ is surjective.

\noindent\eqref{lemc1}$\Longleftrightarrow$\eqref{lemc3x}.
Due to Lemma \ref{hodge2} we have a split exact sequence
$$
0\lra H_1(V,\Z)_{\mathrm{tor}}\lra H_1(V,\Z)
\os{\pi_*}{\lra} H_1(C,\Z)\lra 0.
$$
By the universal coefficients theorem we have
$$
0\lra H^1(C,\Z/n)\lra H^1(V,\Z/n)
\lra \Hom(H_1(V,\Z)_{\mathrm{tor}},\Z/n)\lra 0
$$
and hence the desired equivalence.

\noindent\eqref{lemc1}$\Longrightarrow$\eqref{lemc4}.
It follows from \eqref{lemc1} that we have
$$
H_1(V,\Z/n)=H_1(V,\Z)/n=(H_1(V,\Z)/H_1(V,\Z)_{\mathrm{tor}})\ot\Z/n.
$$
Therefore it is enough to show that the map
$$H_1(V_t,\Z)\lra
H_1(V,\Z)/H_1(V,\Z)_{\mathrm{tor}}\hra H_1(V,\Q)
$$
is zero. This follows from Lemma \ref{hodge2}.

\noindent\eqref{lemc4}$\Longrightarrow$\eqref{lemc3x}.
We want to show that
the map $e_*:H_1(C,\Z/n)\to H_1(V,\Z/n)$ is surjective.
To do this, it is enough
to show that $H_1(S^0,\Z/n)\to H_1(V,\Z/n)$ is surjective.
There is the split exact sequence
$$
1\lra \pi_1(V_t)\lra \pi_1(U^0)\lra \pi_1(S^0)\lra 1.
$$
In particular we have that $\pi_1(V_t)\op\pi_1(S^0)^{\mathrm{ab}}
\to\pi_1(U^0)^{\mathrm{ab}}$ is surjective and hence
so is $H_1(V_t,\Z/n)\op H_1(S^0,\Z/n)\to H_1(U^0,\Z/n)$.
On the other hand since
$H_1(U^0,\Z/n)\to H_1(V,\Z/n)$ is surjective, 
we have the surjectivity of  
$H_1(V_t,\Z/n)\op H_1(S^0,\Z/n)\to H_1(V,\Z/n)$.
However by \eqref{lemc4}, $H_1(V_t,\Z/n)\to H_1(V,\Z/n)$
is zero. Thus we have the surjectivity of 
$H_1(S^0,\Z/n)\to H_1(V,\Z/n)$.
\end{pf}

\begin{lem}\label{simplyc}
If there is a singular fiber which is not a multiplicative one,
then the equivalent conditions in 
Proposition \ref{equiporp} are satisfied for all $n$ such that
$(n,6r_1\cdots r_s)=1$.
\end{lem}
\begin{pf}
We see that \eqref{lemc4} is satisfied.
Let $Y_j$ be the singular fiber which is not multiplicative.
Then it is simply connected.
Since $H_1(V_t,\Z/n)\to H_1(V,\Z/n)$ factors through $H_1(Y_j,\Z/n)$
it is zero.
\end{pf}
\begin{lem}\label{equiporpcor}
Assume that $n$ is prime to $6r_1\cdots r_s$ and satisfies
the equivalent conditions in Proposition \ref{equiporp}.
Then $n$ is prime to the cardinality of the torsion part of
$H^2(V,\Z)$.
\end{lem}
\begin{pf}
This follows from a commutative diagram
$$
\begin{CD}
0@>>> H^1(V,\Z)/n@>>> H^1(V,\Z/n)@>>> H^2(V,\Z)[n]@>>>0\\
@.@AAA@AAA\\
@.H^1(C,\Z)/n@>{\cong}>> H^1(C,\Z/n)
\end{CD}
$$
and \eqref{lemc3x}.
\end{pf}
We put
$$
C_{\Z/n}(j):=\Coker(H^2_Y(V,\Z/n(j))\lra H^2(V,\Z/n(j))),
$$
$$
I_{\Z/n}(j):=\Image(H^2_Y(V,\Z/n(j))\lra H^2(V,\Z/n(j))).
$$
From \eqref{localexseqcohomology} we have exact sequences
\begin{equation}\label{41}
0\lra I_{\Z/n}(j)\lra H^2(V,\Z/n(j))\lra C_{\Z/n}(j)\lra 0,
\end{equation}
\begin{equation}\label{42}
0\lra C_{\Z/n}(j)\lra H^2(U^0,\Z/n(j))\os{\delta}{\lra} 
H^3_Y(V,\Z/n(j))\lra H^3(V,\Z/n(j)).
\end{equation}
Note that $\delta$ in \eqref{42} gives rise to the boundary map
(cf. \S \ref{boundarysect})
\begin{equation}\label{43}
\partial_B: H^2(U^0,\Z/n(j))\lra 
H^3_Y(V,\Z/n(j))\os{\cong}{\lra}\Z/n(j-2)^{\op s}.
\end{equation}
\begin{prop}\label{cohprop}
Assume that $n$ is prime to $6r_1\cdots r_s$ and satisfies
the equivalent conditions in Proposition \ref{equiporp}.
\begin{enumerate}
\renewcommand{\labelenumi}{$(\theenumi)$}
\item\label{cohprop1}
The exact sequence \eqref{41} has a splitting induced from
\eqref{splitting}. 
In particular 
$C_{\Z/n}(j)$ is a direct summand of $H^2(V,\Z/n(j))$.
Moreover we have $I_{\Z/n}(j)\cong \Z/n(j-1)^\op$.
\item\label{cohprop2}
The boundary
map $\partial_B$ is surjective.
Hence we have an exact sequence
\begin{equation}\label{420}
0\lra C_{\Z/n}(j)\lra H^2(U^0,\Z/n(j))\os{\partial_B}{\lra} 
\Z/n(j-2)^{\op s}\lra0.
\end{equation}
\item\label{cohprop3}
$H^2(V,\Z/n)\cong (H^2(V,\Z)/H^2(V,\Z)_{\mathrm{tor}})\ot\Z/n$.
\end{enumerate}
\end{prop}
\begin{pf}
\noindent\eqref{cohprop1}.
This follows from Lemma \ref{intlem}.

\noindent\eqref{cohprop2}.
It is enough to show that the map 
$H^3_Y(V,\Z/n(j))\to H^3(V,\Z/n(j))$ in \eqref{42} is zero.
Due to the Poincare-Lefschetz duality we have
$$
\begin{CD}
H^3_Y(V,\Z/n(j))@>>> H^3(V,\Z/n(j))\\
@A{\cong}AA@AA{\cong}A\\
\bigoplus_{i=1}^s H_1(Y_i,\Z/n(j-2))@>>> H_1(V,\Z/n(j-2)).
\end{CD}
$$
We see that $
H_1(Y_i,\Z/n)\to H_1(V,\Z/n)$ is zero for each $i$.
Since there is a surjective map $H_1(V_t,\Z/n)\to H_1(Y_i,\Z/n)$,
it is enough to see that $H_1(V_t,\Z/n)\to H_1(V,\Z/n)$ is zero.
However this follows from Proposition \ref{equiporp} \eqref{lemc4}.

\noindent\eqref{cohprop3}.
This follows from Proposition \ref{equiporp} \eqref{lemc1},
Lemma \ref{equiporpcor} and the exact sequence
$$
\begin{CD}
H^2(V,\Z)@>{n}>> H^2(V,\Z)@>>> H^2(V,\Z/n) 
@>>> H^3(V,\Z)[n]@>>> 0\\
@. @. @. @|\\
@. @. @. H_1(V,\Z)[n].\\
\end{CD}
$$
\end{pf}


\def\y{{F}}

\subsection{Proof of Theorem B}\label{proof1sect}
Due to Lemma \ref{d=bmore} the assertions in Theorem B
do not depend on
the choice of $U_F^0$.
Thus we may assume that
$U^0_\y=X_\y-$(all singular fibers).
Let $\pi'_\y:V_\y\to C_\y$ be the minimal elliptic surface
associated to $X_\y$.
Then $H^2_\et(X_{\ol{\y}},\Z/n)\cong 
H^2_\et(V_{\ol{\y}},\Z/n)\op\Z/n^{\op}$
and $H^i_\et(X_{\ol{\y}},\Z/n)\cong 
H^i_\et(V_{\ol{\y}},\Z/n)$ for $i\not=2$.
Therefore if we replace $X_\y$ with $V_\y$, the assumptions
in Theorem B hold. Thus we may replace $X_F$ with $V_F$.
We put
$$
C^\et_{\Z/p^\nu}(j):=\Coker(H^2_{\et,Y}(V_{\ol{\y}},\Z/p^\nu(j))
\lra H^2_\et(V_{\ol{\y}},\Z/p^\nu(j))),
$$
$$
I^\et_{\Z/p^\nu}(j):=\Image(H^2_{\et,Y}(V_{\ol{\y}},\Z/p^\nu(j))\lra 
H^2_\et(V_{\ol{\y}},\Z/p^\nu(j)))
$$
and $C^\et_{\Z_p}:=\varprojlim C^\et_{\Z/p^\nu}$,
$I^\et_{\Z_p}:=\varprojlim I^\et_{\Z/p^\nu}$.
Thanks to a theorem of M. Artin,
there are the natural isomorphisms
$$
H_\et^\bullet(W_{\ol{\y}},\Z/n(j))\cong
H_B^\bullet(W_{\C},\Z/n(j))
$$
for any separated scheme $W$ of finite type over $\ol{\y}$
with an embedding $\ol{\y}\hra\C$.
Note that $H_\et^\bullet(W_{\ol{\y}},\Z_p)$ is torsion free if and only if
$p$ is prime to the cardinality of $H_B^\bullet(W_{\C},\Z)_{\tor}$.
Therefore we can apply Proposition \ref{cohprop} for the \'etale
cohomology groups of $V_{\ol{\y}}$.
Thus we have the exact sequence
\begin{equation}\label{420et}
0\lra C_{\Z_p}^\et(j)\lra H^2_\et(U^0_{\ol{\y}},\Z_p(j))
\os{\partial_\et}{\lra} 
\Z_p(j-2)^{\op s}\lra0
\end{equation}
of $G_\y$-modules
from Proposition \ref{cohprop} \eqref{cohprop2}.
$C_{\Z_p}^\et(j)$ is a direct summand of 
$H^2_\et(V_{\ol{\y}},\Z_p(j))$ by
Proposition \ref{cohprop} \eqref{cohprop1}.
Moreover $H^2_\et(V_{\ol{\y}},\Z_p(j))$ is torsion free 
and there is an exact sequence
\begin{equation}\label{420et10}
0\lra H^2_\et(V_{\ol{\y}},\Z_p(j))
\os{p^\nu}{\lra} H^2_\et(V_{\ol{\y}},\Z_p(j))
\lra H^2_\et(V_{\ol{\y}},\Z/p^\nu(j))\lra0
\end{equation}
by Proposition \ref{cohprop} \eqref{cohprop3}.

\subsubsection{Proof of Theorem B (1)}
For a scheme $U'$ with $U^0_{\ol{F}}\to U'$
we put
$$
\vg(U')_{\et}:=\Image(\vg(U',\K_2)\ot\Z_p\lra 
H_\et^2(U^0_{\ol{F}},\Z_p(2))).
$$
Then we show
\begin{equation}\label{vnl1}
\vg(U^0_F)_{\et}=\vg(U^0_{\ol{F}})_{\et}.
\end{equation}
\begin{claim}\label{c1}
$\vg(U^0_F)_{\et}\ot\Q_p=
\vg(U^0_{\ol{F}})_{\et}\ot\Q_p$.
\end{claim}
Note that $\vg(U^0_{F'})_{\et}\cap C_{\Z_p}^\et(2)$ is torsion
so that we have $\vg(U^0_{F'})_{\et}\ot\Q_p\cong
\partial_\et(\vg(U^0_{F'})_{\et})\ot\Q_p=
\partial\vg(U^0_{F'},\K_2)\ot\Q_p$ for $F'=F$ or $\ol{F}$.
Therefore the assertion follows from Lemma \ref{d=bmore} (2).
\begin{claim}\label{c11}
$\vg(U^0_{\ol{F}})_{\et}\subset 
H_\et^2(U^0_{\ol{F}},\Z_p(2))^{G_F}$.
\end{claim}
In fact, it follows from Claim \ref{c1} that we have
$\vg(U^0_{\ol{F}})_{\et}\ot\Q_p=\vg(U^0_F)_{\et}\ot\Q_p
\subset H_\et^2(U^0_{\ol{F}},\Q_p(2))^{G_F}$.
Then the assertion follows from the fact that 
$H_\et^2(U^0_{\ol{F}},\Z_p(2))$ is torsion free.
\begin{claim}\label{c2}
$H^2_\et(U^0_F,\Z_p(2))/\vg(U^0_F)_\et$ is torsion free.
\end{claim}
Put $M:=H^2_\et(U^0_F,\Z_p(2))/\vg(U^0_F)_\et$.
Due to the finiteness of the Galois cohomology groups of $F$,
$M$ is finitely generated over $\Z_p$.
Therefore it is enough to show that
$$
M/p^\nu\lra M/p^{\nu+1},\quad x\longmapsto p\cdot x
$$
is injective. Due to Suslin's exact sequence (Theorem \ref{universalc0})
we have an injective map $M/p^\nu \hra H^1_\Zar(U^0_F,\K_2)[p^\nu]
$ and a commutative diagram
$$
\begin{CD}
M/p^\nu@>>>H^1_\Zar(U^0_F,\K_2)[p^\nu]\\
@V{p}VV@VV{\bigcap}V\\
M/p^{\nu+1}@>>>H^1_\Zar(U^0_F,\K_2)[p^{\nu+1}].
\end{CD}
$$
Thus Claim \ref{c2} follows.
\begin{claim}\label{c3}
There is an exact sequence
$$
0\lra H^2_\et(S^0_F,\Z_p(2))
\lra H^2_\et(U^0_F,\Z_p(2))
\lra H^2_\et(U^0_{\ol{F}},\Z_p(2))^{G_F}
\lra 0.
$$
This is split by the section $e$.
\end{claim}
Consider the Hochschild-Serre spectral sequences
$$
E_{2,U^0}^{ij}=H^i(F,H^j_\et(U^0_{\ol{F}},\Z/p^\nu(2)))
\Longrightarrow
H^{i+j}_\et(U^0_F,\Z/p^\nu(2))
$$
$$
E_{2,S^0}^{ij}=H^i(F,H^j_\et(S^0_{\ol{F}},\Z/p^\nu(2)))
\Longrightarrow
H^{i+j}_\et(S^0_F,\Z/p^\nu(2)).
$$
By Proposition \ref{equiporp} \eqref{lemc3},
$H^1_\et(S^0_{\ol{F}},\Z/p^\nu(2))
\cong
H^1_\et(U^0_{\ol{F}},\Z/p^\nu(2))$.
Therefore we have $E_{2,S^0}^{ij}\cong E_{2,U^0}^{ij}$ for
$j=0,1$.
Thus we have
$E_{\infty,S^0}^{11}\cong E_{\infty,U^0}^{11}$ and
$E_{\infty,S^0}^{20}\cong E_{\infty,U^0}^{20}$, and hence
an exact sequence
$$
0\lra H^2_\et(S^0_F,\Z/p^\nu(2))
\lra H^2_\et(U^0_F,\Z/p^\nu(2))
\os{v}{\lra} H^2_\et(U^0_{\ol{F}},\Z/p^\nu(2))^{G_F}.
$$
The rest of the proof is to show that the right arrow 
$v$ is surjective, which is equivalent to
$E_{2,U^0}^{02}=E_{\infty,U^0}^{02}$. To do this it is enough
to show that the map $E_{2,U^0}^{02}\to E_{2,U^0}^{21}$ and
$E_{3,U^0}^{02}\to E_{3,U^0}^{30}$ are zero.
The latter follows from the injectivity of $E_{2,U^0}^{i0}
\to E_{U^0}^{i}$. We show the former.
There is a commutative diagram
$$
\begin{CD}
E_{2,U^0}^{02}@>{a}>> E_{2,U^0}^{21}@>{b}>> E_{U^0}^3\\
@AAA@A{\cong}AA@AA{\bigcup}A\\
E_{2,S^0}^{02}@>{c}>> E_{2,S^0}^{21}@>{d}>> E_{S^0}^3
\end{CD}
$$
with exact rows.
Therefore $a=0$ $\Longleftrightarrow$ $b$ is injective 
$\Longleftrightarrow$ $d$ is injective
$\Longleftrightarrow$ $c=0$. 
However this is clear as $E_{2,S^0}^{i2}=0$.

\medskip

Now we show \eqref{vnl1}.
By Claim \ref{c3} we have a commutative diagram

$$
\xymatrix{
0\ar[r]&
 H^2_\et(S^0_F,\Z_p(2))\ar[r]&
H^2_\et(U^0_F,\Z_p(2))\ar[r]&
 H^2_\et(U^0_{\ol{F}},\Z_p(2))^{G_F}\ar[r]&0\\
0\ar[r]&
 \vg(S^0_F,\K_2)\ot\Z_p \ar[r]^{\pi_F^*}\ar[u]&
 \vg(U^0_F,\K_2)\ot\Z_p \ar[r]\ar[u]&
 \Coker~\pi_F^*\ar[r]\ar[u]&0\\
}
$$
with split exact rows by the section $e_F$. 
Together with Claim \ref{c2}
we have that
\begin{equation}\label{torsu}
H^2_\et(U^0_{\ol{F}},\Z_p(2))^{G_F}/\vg(U^0_F)_\et
\end{equation}
is torsion free.
Therefore,
in the following commutative diagram
$$
\begin{CD}
\vg(U^0_F)_\et@>{\subset}>>
H^2_\et(U^0_{\ol{F}},\Z_p(2))^{G_F}\\
@V{a}VV@|\\
\vg(U^0_{\ol{F}})_\et@>{\subset}>>
H^2_\et(U^0_{\ol{F}},\Z_p(2))^{G_F},
\end{CD}
$$
$\Coker~a$ is torsion free.
On the other hand $\Coker~a$ is torsion by Claim \ref{c1}.
This means that $a$ is bijective.
This completes the proof of \eqref{vnl1} and hence Theorem B (1).

\subsubsection{Proof of Theorem B (2)}
In the same way as Claim \ref{c11} we have
\begin{equation}\label{vnl1claim}
\Image
(\vg(U^0_{\ol{F}},\K_2)\ot\Z_p\to H^2_\et(U^0_{\ol{F}},\Z_p(2)))
\subset H^2_\et(U^0_{\ol{F}},\Z_p(2))^{G_F}.
\end{equation}
Applying ${\mathbb R}\Gamma(G_F,-)$ on \eqref{420et},
we have
$$
\begin{CD}
0@>>>
H^2_\et(U^0_{\ol{F}},\Z/p^\nu(2))^{G_F}
@>>>
(\Z/p^\nu)^{\op s}@>>>
H^1(G_F,C^\et_{\Z/p^\nu}(2))\\
@.@AAA\\
@.\vg(U^0_{\ol{F}},\K_2)/p^\nu
\end{CD}
$$
with an exact row.
Passing to the projective limit, we have
$$
\begin{CD}
0@>>>
H^2_\et(U^0_{\ol{F}},\Z_p(2))^{G_F}
@>{\partial_\et}>>
\Z_p^{\op s}@>>>
\plim{\nu}H^1(G_F,C^\et_{\Z/p^\nu}(2))\\
@.@A{\dlog}AA@AA{=}A\\
@.\vg(U^0_{\ol{F}},\K_2)\ot\Z_p@>{\partial}>>\Z_p^{\op s}.
\end{CD}
$$
We want to show that the cokernel of $\partial$ is torsion free.
In the same way as the proof of Claim \ref{c2}
we can show that
the cokernel of $\dlog$ is torsion free.
Since $C^\et_{\Z/p^\nu}(2)$ is a direct summand of 
$H^2_\et(V_{\ol{F}},\Z/p^{\nu}(2))$,
it is enough to show that 
\begin{equation}\label{plimcohh}
\plim{\nu}H^1(G_F,H^2_\et(V_{\ol{F}},\Z/p^{\nu}(2)))
\end{equation}
is torsion free. 
We have an exact sequence
$$
0\lra H^2_\et(V_{\ol{F}},\Z/p^{\nu-1}(2))
\lra H^2_\et(V_{\ol{F}},\Z/p^{\nu}(2))
\lra H^2_\et(V_{\ol{F}},\Z/p(2))\lra 0
$$
from \eqref{420et10}.
Applying $\R\Gamma(G_F,-)$ to the above,
we see that the following map
$$
H^1(G_F,H^2_\et(V_{\ol{F}},\Z/p^{\nu-1}(2)))
\lra 
H^1(G_F,H^2_\et(V_{\ol{F}},\Z/p^{\nu}(2)))
$$
is injective for all $\nu\geq1$ due to 
the vanishing $H^2_\et(V_{\ol{F}},\Z/p(2))^{G_F}=0$.
Passing to the projective limit,
we have the injectivity of
$$
\plim{\nu}H^1(G_F,H^2_\et(V_{\ol{F}},\Z/p^{\nu}(2)))
\lra 
\plim{\nu}H^1(G_F,H^2_\et(V_{\ol{F}},\Z/p^{\nu}(2))),\quad
x\mapsto p\cdot x.
$$
This means that \eqref{plimcohh} is torsion free.
This completes the proof of Theorem B (2).

\section{Modular elliptic surface}
\label{modularsect}

\def\Kpq{\Phi(X_R,D_R)_{\Z_p}}
\def\kpq{\Phi_{\Z_p}^\prime}

The purpose of this section is to prove Conjecture \ref{conj1}
for a universal elliptic curve $\E_N\to X(N)$ ($N\geq 3$)
over a modular
curve (Theorem \ref{modularthm}).
\subsection{Preliminaries and notations on modular curves}
For a congruence subgroup
$\vg\subset\SL_2(\Z)$
(i.e. $\vg\supset \vg(N)$ for some $N\geq1$),
we denote by $M_k(\vg)$ (resp. $S_k(\vg)$) 
the $\C$-vector space of modular forms (resp. cusp forms)
of weight $k$
with respect to $\vg$.
See \cite{diamond} Def. 1.2.3, or
\cite{shimura} \S 2.1 for the definition
($M_k(\vg)$ is denoted by $A_k(\vg)$ in \cite{shimura}).
\subsubsection{Hecke operators}
We focus on the special case $\vg=\vg(N)$ with $N\geq 3$.
Let
$$
\vg_0=\left\{ 
\begin{pmatrix}
a&b\\
c&d
\end{pmatrix}\in \SL_2(\Z)\left|
\begin{pmatrix}
a&b\\
c&d
\end{pmatrix}\equiv 
\begin{pmatrix}
*&0\\
0&*
\end{pmatrix}\mod N\right.
\right\}.
$$
The group $\vg_0/\vg\cong(\Z/N)^*$ acts on $S_k(\vg)$ by
\begin{equation}\label{cusp1}
f|[\gamma]_k:=(c\tau+d)^{-k}f(\frac{a\tau+b}{c\tau+d}),
\quad
\gamma=
\begin{pmatrix}
a&b\\
c&d
\end{pmatrix}\in \vg_0,
~\tau\in{\frak H}
\end{equation}
where ${\frak H}$ denotes the complex upper half plane.
Since $f|[\gamma]_k$ depends only on $d$ mod $N$, it is sometimes
written as $\langle d\rangle f$, and
$\langle d\rangle$ is called a {\it diamond operator}
(\cite{diamond} 5.2).
Let
$\chi:(\Z/N\Z)^*\to \C^*$ be a Dirichlet character.
We put $\chi(a)=0$ if $(a,N)\not=1$.
We define $S_k(\vg_0,\chi)\subset S_k(\vg)$ as the set of 
all $f\in S_k(\vg)$ satisfying
\begin{equation}\label{cusp1b}
\langle d\rangle f=\chi(d)f\quad
\text{for all }d\in(\Z/N)^*.
\end{equation}
Then we have
\begin{equation}\label{cusp0}
S_k(\vg)=\bigoplus_\chi S_k(\vg_0,\chi)
\end{equation}
where $\chi$ runs over all Dirichlet characters of $(\Z/N\Z)^*$.
The decomposition \eqref{cusp0}
is mutually orthogonal with respect to
the Petersson inner product
(\cite{shimura} (3.5.4)).

There are the {\it Hecke operators} $T_m$ on $S_k(\vg)$.
See \cite{shimura} \S 3.2 for the definition
(where $T_m$ is denoted by $T(m)$).
They satisfy $T_nT_m=T_{nm}$ for $(n,m)=1$ and
\begin{equation}\label{diamond}
T_{p^{r+1}}=
\begin{cases}
T_pT_{p^r}-p^{k-1}\langle p \rangle T_{p^{r-1}}&p\not|N\\
T_pT_{p^r}&p|N
\end{cases}
\end{equation}
on $S_k(\vg)$
for $r\geq 1$ and any prime number $p$ 
(loc.cit. Thm. 3.24 and (3.5.8)).
The Hecke operators and diamond operators are 
mutually commutative and
normal with respect to
the Petersson inner product (loc.cit. Thm. 3.41).
Therefore they
are simultaneously
diagonalizable. 
A common eigen function for all $T_n$
is called a {\it Hecke eigenform}.
If $f\in S_k(\vg)$ is a Hecke eigenform, then there is a unique
Dirichlet character $\chi$ such that $f\in S_k(\vg_0,\chi)$
(loc.cit. Prop. 3.53).

Put $q_N=\exp(2\pi i\tau/N)$.
Let $f=\sum_{n=1}^\infty c_nq^n_N$ be the Fourier expansion
at $\infty$.
Suppose that $f\in S_k(\vg_0,\chi)$. Then the Fourier expansion
of $T_mf$ is given as follows (loc.cit. (3.5.12))
\begin{equation}\label{hecke0000b}
T_mf=\sum_{n=1}^\infty c^*_nq^n_N,\quad
c^*_n=\sum_{a|(n,m)}\chi(a)a^{k-1}c_{mn/a^2}.
\end{equation}
In particular, if $\langle p\rangle f=\sum_{n=1}c'_nq_N^n$
for a prime number $p$, then we have
\begin{equation}\label{hecke0000}
T_pf=\begin{cases}
\sum_{n=1}^\infty (c_{np}+p^{k-1}c'_{n/p})q^n_N&p\not|N\\
\sum_{n=1}^\infty c_{np}q^n_N&p|N
\end{cases}
\end{equation}
where we put $c'_{n/p}=0$ unless $p|n$.
Suppose that $f\in S_k(\vg_0,\chi)$ is a
Hecke eigenform.
If $c_1=1$, then we say that $f$ is {\it normalized}.
In this case we have $T_mf=c_mf$ and
\begin{equation}\label{hecke000}
c_{nm}=c_nc_m\quad((n,m)=1),
\quad c_{p^{r+1}}=c_{p}c_{p^r}-
\chi(p)p^{k-1}c_{p^{r-1}} \quad (r\geq 1)
\end{equation}
(loc.cit. Thm. 3.43).
It is well-known that the characteristic polynomial of $T_m$
has rational integer coefficients for all $m\geq 1$
(loc.cit. (3.5.20)).
Hence all Fourier coefficients of a normalized Hecke eigenform
are algebraic integers in a number field.
\subsubsection{Modular curves :
Deligne and Rapoport \cite{DeRa}}
Let $\vg=\vg(N)$ with $N\geq 3$.
Put $\O_N:=\Z[\zeta_N,1/N]$
where $\zeta_N$ is a primitive
$N$-th root of unity.
The main result of \cite{DeRa} tells that
there are the
algebraic curve
\begin{equation}\label{DeRamodcurve}
X(N)_{\O_N}\lra \Spec\O_N
\end{equation}
and the open subscheme
$Y(N)_{\O_N}\subset X(N)_{\O_N}$
such that
$$X(N)_{\O_N}\times_{\O_N}\C\cong ({\frak H}\cup\Q\cup\{\infty\})
/\vg(N)
$$
$$Y(N)_{\O_N}\times_{\O_N}\Spec\C\cong {\frak H}/\vg(N).
$$
They are called the {\it modular curves}.
The morphism \eqref{DeRamodcurve} is projective and smooth
(loc.cit. IV. Cor. 2.9), and the geometric fiber
is connected (loc.cit. IV. Cor. 5.5).
The complement $\Sigma_{\O_N}
:=X(N)_{\O_N}-Y(N)_{\O_N}$ is the disjoint union
of copies of $\Spec\O_N$ (loc.cit. VII. Cor. 2.5), which we 
put $\Sigma_{\O_N}=\{P_1,\cdots,P_s\}$.
We call the points $P_i$ the {\it cusps}.

Since
$X(N)_{\O_N}$ is the fine moduli of generalized elliptic curves
with level structure $N$ (loc.cit. IV. Def. 2.4),
we have the {\it universal elliptic curve}
$$\pi:\E_N\lra X(N)_{\O_N}.$$ The morphism $\pi$ is projective
over $X(N)_{\O_N}$, and
smooth over $Y(N)_{\O_N}$.
The fiber
$D_{i,\O_N}:=\pi^{-1}(P_i)$ over a cusp
$P_i\in\Sigma_{\O_N}$ is a standard N\'eron $N$-gon.
More precisely, let
$\Delta_{P_i}:\Spec\O_N[[q_N]]\to X(N)_{\O_N}$ be the formal 
neighborhood of $P_i$.
Then the fiber
$\pi^{-1}(\Delta_{P_i})$ is isomorphic to the regular model
$\E_{q,\O_N[[q_N]]}$ in \S \ref{modelsect}
(loc.cit. VII. Cor. 2.4).
In particular, the elliptic surface $\pi:\E_N\to X(N)_{\O_N}$
has only singular fibers of type $I_N$ over the cusps
and satisfies {\bf (Rat)}.

Let $D_{\O_N}:=\sum_i D_{i,\O_N}=\pi^{-1}(\Sigma_{\O_N})$.
Define an invertible sheaf 
$\Omega$
by the exact sequence
$$
0\lra \pi^*\Omega^1_{X(N)/\O_N}(\Sigma_{\O_N})
\lra \Omega^1_{\E_N/\O_N}(\log D_{\O_N})
\lra \Omega\lra 0.
$$
The sheaf $\Omega$ is isomorphic to the dualizing sheaf
${\cal H}^{-1}(R\pi^!\O)$.
Put $
\omega=\pi_*\Omega
$.
This is also an invertible sheaf
which is generated by $du/u$ locally at the cusps
where $du/u$ is the canonical invariant 1-form of the Tate curve
(cf. \eqref{canonicalinvariant}).
Then there is the natural isomorphism
\begin{equation}\label{dera1}
\omega^{\ot 2}\os{\cong}{\lra}
\Omega^1_{X(N)/\O_N}(
\Sigma_{\O_N}),
\quad \frac{du}{u}^{\ot 2}\longmapsto\frac{dq_N}{q_N}
\end{equation}
(loc.cit. VI. 4.5).
Moreover we have the natural isomorphism
\begin{equation}\label{dera0}
M_k(\vg)\os{\cong}{\lra} \vg(X(N)_{\O_N},\omega^{\ot k})
\ot_{\O_N}\C,
\quad f\longmapsto f \frac{du}{u}^{\ot k}
\end{equation}
(loc.cit. VII. 4.6). 
We identify $M_k(\vg)$ with the set of sections
of $\omega^{\ot k}$ by \eqref{dera0}.
Then $S_k(\vg)$ can be identified with
$$
\vg(X(N)_{\O_N},\Omega^1_{X(N)/\O_N}\ot \omega^{\ot k-2})
\ot_{\O_N}\C.
$$
For an $\O_N$-algebra $R$, we put
\begin{equation}\label{dera2}
M_k(\vg)_R\os{\mathrm{def}}{=}\vg(X(N)_{\O_N},\omega^{\ot k})
\ot_{\O_N}R,
\end{equation}
\begin{equation}\label{dera3}
S_k(\vg)_R\os{\mathrm{def}}{=}
\vg(X(N)_{\O_N},\Omega^1_{X(N)/\O_N}\ot \omega^{\ot k-2})
\ot_{\O_N}R.
\end{equation}
Let $f_{i,\chi}\in S_k(\vg,\chi)$ be the normalized
Hecke eigenforms and $f_{i,\chi}=\sum_{n\geq1}c_{i,\chi}(n)q_N^n$
the Fourier expansions. 
Let $F_N=\Q(\zeta_N,c_{i,\chi}(n))_{i,\chi,n}$ which is a finite
extension of $\Q$.
Then $f_{i,\chi}$ is contained in 
$S_k(\vg)_{F_N}=S_k(\vg)_{\O_N}\ot_{\O_N}F_N$ and
form a basis of it:
\begin{equation}\label{basis000}
S_k(\vg)_{F_N}=\bigoplus_{i,\chi} F_N\cdot f_{i,\chi}.
\end{equation}

\subsection{Dlog image of $K_2$ of modular elliptic
surface}
Let $N\geq 3$ be an integer
and $p$ a prime number such that $p\not|2N$.
Let $K$ be a finite unramified extension of $\Q_p$ of
degree $d$. Let $R$ be the ring of integers in $K$.
Suppose $\zeta_N\in R$.
Let $X(N)_R:=X(N)_{\O_N}\times_{\O_N} R$
be the modular curve
over $R$ and 
$$\pi_R: X_R=\E_N\times_{\O_N}R\lra X(N)_R$$ the universal 
elliptic curve.
Let $\Sigma_R=\{P_1,\cdots,P_s\}\subset X(N)_R$ be the cusps. We put
$D_R=\sum_i\pi^{-1}_R(P_i)$ and $U_R=X_R-D_R$.
\begin{thm}\label{modularthm}
$\dlog\vg(U_{\ol{K}},\K_2)\ot_\Z\Q_p=\Phi(X_R,D_R)_{\Z_p}\ot_{\Z_p}\Q_p.$
\end{thm}
\subsection{Proof of Theorem \ref{modularthm}}
We use the following theorem of A. Beilinson.
\begin{thm}[\cite{bei-modular} Cor. 3.1.8.]
The rank of $\dlog\vg(U_{\ol{K}},\K_2)$ is equal to $s$(=the 
number of cusps of $X(N)$). \label{beithm}
\end{thm}
More precisely, Beilinson constructed the {\it Eisenstein symbols}
in $\vg(U_{\ol{K}},\K_2)$, and showed that the boundary map
$\partial$ induces the bijection on the space of
the Eisenstein symbols.
See \cite{scholl} \S 7 for another proof.

\medskip

In order to prove Theorem \ref{modularthm}
it is enough to show 
$\rank_{\Z_p}\Kpq\leq s$.
Recall the definition \eqref{dera2} and \eqref{dera3}.
We have from \eqref{dera1} the natural isomorphisms
\begin{equation}\label{dera5}
M_3(\vg)_R\os{\cong}{\lra}
\vg(X_R,\Omega^2_{X_R/R}(\log D_R)),\quad
S_3(\vg)_R\os{\cong}{\lra}
\vg(X_R,\Omega^2_{X_R/R}).
\end{equation}
Therefore in order to show 
$\rank_{\Z_p}\Kpq\leq s
$,
it is enough to show
\begin{equation}\label{preclaim}
(S_3(\vg)_R\cap \Kpq)\ot\Q_p=0.
\end{equation}

\subsubsection{Lemmas on the Fourier coefficients}
Fix
a cyclotomic basis $\underline{\mu}=\{\zeta_1,\cdots,\zeta_d\}$
of $R$.
Let
$
\kpq
\subset S_3(\vg)_R
$
be the $\Z_p$-submodule generated by
$f\in S_3(\vg)_R$ such that if we express
$$
f=\sum_{k=1}^\infty \sum_{i=1}^d
a^{(i)}_k\frac{\zeta_iq_N^k}{1-\zeta_iq_N^k},
\quad
a^{(i)}_k\in \Z_p$$
then 
$
 a^{(i)}_k\equiv0 \mod p^{2k}\Z_p$
for all $k$ and $i$.
Since $\gamma\in\SL_2(\Z)$ acts on the cusps transitively,
we see
$$
S_3(\vg)_R\cap \Kpq=\bigcap_{\gamma\in \SL_2(\Z)}\gamma\kpq.
$$
\begin{lem}\label{st1}
Let $f\in S_3(\vg)_R$ be a cusp form of weight 3 and level $N$.
If $f\in \Kpq$ then 
$\langle d \rangle f\in \Kpq$ for all $d\in (\Z/N)^*$.
Moreover if $f\in \kpq$ then 
$T_mf\in \kpq$ for all $m\geq1$.
\end{lem}
\begin{pf}
An element $\gamma\in \vg_0/\vg\cong(\Z/N\Z)^*$ induces an
automorphism
of the universal 
elliptic curve $X_R\to X(N)_R$ in a natural way, which
we denote by $\sigma_\gamma$.
By the definition \eqref{cusp1} we have
\begin{equation}\label{auo}
\sigma_\gamma^*(f \frac{dq_N}{q_N}\frac{du}{u})=
\langle d \rangle f
\frac{dq_N}{q_N}\frac{du}{u},
\quad
\gamma=
\begin{pmatrix}
a&b\\
c&d
\end{pmatrix}\in \vg_0.
\end{equation}
Therefore if $f\in\Kpq$ then \eqref{auo} is also
contained in $\Kpq$, namely
$\langle d \rangle f\in\Kpq$ for all $d\in (\Z/N)^*$.

Next suppose that $f\in \kpq$.
It is enough to show $T_\l f\in\kpq$ for all prime $\l$.
Express
$$
f=\sum_{k=1}^\infty \sum_{i=1}^d
a^{(i)}_k\frac{\zeta_iq_N^k}{1-\zeta_iq_N^k},
\quad
\langle \l \rangle f=\sum_{k=1}^\infty \sum_{i=1}^d
b^{(i)}_k\frac{\zeta_iq_N^k}{1-\zeta_iq_N^k}
$$
with $a^{(i)}_k,b^{(i)}_k\in \Z_p$.
Here we put $\langle \l \rangle f=0$ if $\l|N$.
Since $f$ and $\langle \l \rangle f$ are contained in $\kpq$,
we have
\begin{equation}\label{modprrel}
a^{(i)}_k\equiv b^{(i)}_k\equiv 0
\mod k^2\Z_p \quad \text{for all }k~, i
\end{equation}
by definition.
Using \eqref{hecke0000},
a direct calculation yields
$$
T_\l(f)=\sum_{i=1}^d\left(
\sum_{(\l,k)=1}
a^{(i)}_k
\frac{\zeta_i^\l q_N^k}{1-\zeta^\l_iq_N^k}+
+\sum_{\l| k}
a^{(i)}_k
\frac{\zeta_i q_N^k}{1-\zeta_iq_N^k}+
\sum_{k=1}^\infty
\l^2b^{(i)}_k\frac{\zeta_iq_N^{k\l}}{1-\zeta_iq_N^{k\l}}
\right).
$$
In order to show $T_\l(f)\in\kpq$, the only non trivial
part is to show
$$
a^{(i)}_k
\frac{\zeta_i^\l q_N^k}{1-\zeta^\l_iq_N^k}\in \kpq.
$$However this can be proved
by the same argument as in the proof of Lemma \ref{dependrem}.
\end{pf}

\begin{lem}\label{modplemma}
Let $A$ be a discrete valuation ring
and $\pi$ a uniformizer.
Let $\alpha\in A$ be a non-zero element such that 
$e:=\ord_{\pi}(\alpha)\geq 1$.
Let $\{a_n\}_{n\geq 0}$ be a sequence in $A$ satisfying
\begin{equation}\label{heckea000}
a_0=1,\quad a_{n+2}=a_1a_{n+1}-\alpha a_{n}, \quad n\geq 0.
\end{equation}
Suppose that there are integers $d\geq1$ and $s\geq0$
such that
\begin{equation}\label{hecke5}
a_{n+d}\equiv a_n \mod \pi^{en-s}\quad\text{for all }n\geq0.
\end{equation}
Then there is a $d$-th root of unity $\zeta$ such that
\begin{equation}\label{hecke10}
a_n=\frac{\zeta^{n+1}-(\zeta^{-1}\alpha)^{n+1}}{\zeta-\zeta^{-1}\alpha},
\quad n\geq0.
\end{equation}
\end{lem}
\begin{pf}
We have from \eqref{heckea000} that
\begin{equation}\label{hecke00}
a_{n+k+2}=a_{k+1}a_{n+1}-
\alpha a_{k}a_{n} \quad (n,k\geq 0).
\end{equation}
Put $k=d-1,d-2$ in \eqref{hecke00} and apply \eqref{hecke5}:
\begin{equation}\label{hecke6}
-\alpha
a_{d-1}a_{n}\equiv
(1-a_d)a_{n+1}
\mod \pi^{en+e-s}, \quad n\geq 0,
\end{equation}
\begin{equation}\label{hecke7}
(1+\alpha a_{d-2})a_{n}\equiv  a_{d-1}a_{n+1}
\mod \pi^{en-s}, \quad n\geq 0.
\end{equation}
Repeating them, we have
\begin{equation}\label{heckeA}
(-\alpha
a_{d-1})^da_{n}\equiv
(1-a_d)^da_{n+d}\os{\eqref{hecke5}}{\equiv}
(1-a_d)^da_{n}
\mod \pi^{en-s}, \quad n\geq 0,
\end{equation}
\begin{equation}\label{heckeB}
(1+\alpha a_{d-2})^da_{n}\equiv  a^d_{d-1}
a_{n+d}\os{\eqref{hecke5}}{\equiv}  a^d_{d-1}
a_{n}
\mod \pi^{en-s}, \quad n\geq 0.
\end{equation}

We first claim that $a_1\not\equiv 0$ mod $\pi$.
Suppose $a_1\equiv 0$ mod $\pi$.
Then we have $a_n\equiv a_1^n\equiv 0$ 
mod $\pi$ for all $n\geq1$
by \eqref{heckea000}. \eqref{heckeA} yields
$$
a_n\equiv
\left(\frac{-\alpha a_{d-1}}{1-a_d}\right)^da_n
\equiv\cdots\equiv
\left(\frac{-\alpha a_{d-1}}{1-a_d}\right)^{dk}a_n
\equiv 0
\mod \pi^{en-s},\quad k\gg1
$$
as $\ord_\pi(\alpha)=e>0$.
Suppose that there is an integer $m\geq 0$
such that $a_n\equiv0$ mod $\pi^{en+m-s}$ for all $n\geq1$.
Then we have
$$
\alpha a_{n}=a_1a_{n+1}-a_{n+2}\equiv
0 \mod \pi^{en+e+1+m-s}
$$
and hence $a_n\equiv0$ mod $\pi^{en+m+1-s}$.
Therefore the induction yields $a_n\equiv0$ mod $\pi^{en+m-s}$
for all $m\geq 0$ and $n\geq 1$. This means $a_n=0$ for all $n\geq 1$,
which contradicts \eqref{heckea000}.

We now have $a_1\not\equiv 0$ mod $\pi$.
Then $a_n\equiv a_1^n\not\equiv 0$ 
mod $\pi$ for all $n\geq0$
by \eqref{heckea000}.
Therefore \eqref{heckeA} and \eqref{heckeB} imply
$(-\alpha
a_{d-1})^d\equiv(1-a_d)^d$ and
$(1+\alpha a_{d-2})^d\equiv a^d_{d-1}$ mod $\pi^{en-s}$ for all
$n\geq 0$, which means
$(-\alpha
a_{d-1})^d=(1-a_d)^d$ and
$(1+\alpha a_{d-2})^d= a^d_{d-1}$.
Put
\begin{equation}\label{hecke8}
-\alpha
a_{d-1}=\zeta(1-a_d),\quad
1+\alpha a_{d-2}=\zeta'a_{d-1}
\end{equation}
where $\zeta,\zeta'$ are $d$-th roots of unity.
We claim $\zeta=\zeta'$. In fact
due to \eqref{hecke6} and \eqref{hecke7},
we have
$$
-\alpha
a_{d-1}a_n\equiv
-\zeta^{-1}\alpha
a_{d-1}a_{n+1}
\mod \pi^{en+e-s}, \quad n\geq 0,
$$
$$
a_{d-1}a_n\equiv  \zeta^{\prime -1}a_{d-1}a_{n+1}
\mod \pi^{en-s}, \quad n\geq 0.
$$
The above implies $\zeta\equiv\zeta'$ mod $\pi^{en-s}$
for all $n\geq0$, hence $\zeta=\zeta'$.
Thus we have
\begin{align*}
0\os{\eqref{hecke8}}{=}
1-a_d+\zeta^{-1}\alpha a_{d-1}&
\os{\eqref{heckea000}}{=}
1-(a_1a_{d-1}-\alpha a_{d-2})
+\zeta^{-1}\alpha a_{d-1}\\
&\os{\eqref{hecke8}}{=}
\zeta a_{d-1}-a_1a_{d-1}+\zeta^{-1}\alpha a_{d-1}.
\end{align*}
Since $a_{d-1}\not=0$, we have 
$
a_1=\zeta+\zeta^{-1}\alpha.
$
This yields
\eqref{hecke10} by the induction on $n$.
\end{pf}

\subsubsection{End of the proof}
We finish the proof of Theorem \ref{modularthm}.
It is enough to show \eqref{preclaim}.
Suppose that there is a non-zero $f=\sum_{n=1}^\infty
c(n)q_N^n\in S_3(\vg)_R$
such that 
$f\in \Kpq$.
Let $f$ be expressed
$$
f=\sum_{k=1}^\infty \sum_{i=1}^d
a_i(k)\frac{\zeta_iq_N^k}{1-\zeta_iq_N^k}
$$with $a_i(k)\in \Z_p$.
Let $\sigma:R\to R$ be the Frobenius automorphism.
Then
$$ c(p^r)=\sum_{k=0}^r\sum_{i=1}^d
a_i(p^k)\zeta_i^{p^{r-k}}
=\sigma(c(p^{r-1}))+\sum_{i=1}^d
a_i(p^r)\zeta_i, \quad r\geq 1
$$
as $\zeta_i^\sigma=\zeta_i^p$.
Since $f\in \Kpq\subset \kpq$ we have
$$
c(p^{r+1})\equiv\sigma(c(p^{r}))
\mod p^{2r+2}R, \quad r\geq 0.
$$
Repeating it, we have
\begin{equation}\label{diamondexp}
c(p^{r+d})\equiv c(p^{r})
\mod p^{2r+2}R, \quad r\geq 0
\end{equation}
as $\sigma^d=1$.
Let $f_{i,\chi}=\sum_{n=1}^\infty
c_{i,\chi}(n)q_N^n$ be the normalized Hecke eigenforms, and 
express
$f=\sum_{i,\chi}\alpha_{i,\chi}f_{i,\chi}$ with $\alpha_{i,\chi}
\in K(c_{i,\chi}(n))_{i,\chi,n}$ (cf. \eqref{basis000}).
Then \eqref{diamondexp} is written as
\begin{equation}\label{diamondexp1}
\sum_{i,\chi}\alpha_{i,\chi}c_{i,\chi}(p^{r+d})
\equiv 
\sum_{i,\chi}\alpha_{i,\chi}c_{i,\chi}(p^{r})
\mod p^{2r+2}R, \quad r\geq 0.
\end{equation}
On the other hand
By Lemma \ref{st1}, we have $T_mf=
\sum_{i,\chi}c_{i,\chi}(m)\alpha_{i,\chi}f_{i,\chi}\in\kpq$ 
for all $m\geq 1$.
Similarly to \eqref{diamondexp1} we have
\begin{equation}\label{diamondexp2}
\sum_{i,\chi}c_{i,\chi}(m)\alpha_{i,\chi}c_{i,\chi}(p^{r+d})
\equiv 
\sum_{i,\chi}c_{i,\chi}(m)\alpha_{i,\chi}c_{i,\chi}(p^{r})
\mod p^{2r+2}R, \quad r\geq 0
\end{equation}
for all $m\geq1$.
Since $f_{i,\chi}$ are linearly independent, 
\eqref{diamondexp2} yields that there is an integer $s\geq 0$
which does not depend on $r$ such that
\begin{equation}\label{diamondexp3}
\alpha_{i,\chi}c_{i,\chi}(p^{r+d})
\equiv 
\alpha_{i,\chi}c_{i,\chi}(p^{r})
\mod p^{2r+2-s}R', \quad r\geq 0
\end{equation}
for all $i$ and $\chi$ where $R'$ is the ring of integers in 
$K(c_{i,\chi}(n))_{i,\chi,n}$.
Take $f_{i,\chi}$ such that $\alpha_{i,\chi}\not=0$.
Due to \eqref{hecke000} and
\eqref{diamondexp3}, one can
apply Lemma \ref{modplemma} for $a_n=c_{i,\chi}(p^{n})$.
Then we have
\begin{equation}\label{diamondexp4}
c_{i,\chi}(p^{r})
=\frac{\zeta^{r+1}-(\zeta^{-1}\chi(p)p^2)^{r+1}}
{\zeta-\zeta^{-1}\chi(p)p^2}
, \quad r\geq 0.
\end{equation}
However, since $f_{i,\chi}$ is a cusp form of weight 3, we 
have an estimate
$\vert c_{i,\chi}(n)\vert\leq C n^{3/2}$
where $\vert c_{i,\chi}(n)\vert$ denotes the complex absolute
value (\cite{shimura} Lem. 3.62, \cite{diamond} Prop. 5.9.1). 
This is the contradiction.
This completes the proof of \eqref{preclaim} and hence
Theorem \ref{modularthm}.


\section{Example : elliptic surfaces $Y^2=X^3+X^2+t^n$}
\label{expmsect}
Let $\pi:X\to C=\P^1_{\C}$ 
be the minimal elliptic surface over $\C$
such that the general fiber $\pi^{-1}(t)$ is the elliptic curve
defined by
\begin{equation}\label{elsurfeq}
Y^2=X^3+X^2+t^n, \quad n\geq 1.
\end{equation}
The functional $j$-invariant is $-64/(t^n(1+27t^n/4))$.
Put $n=6k+l$ where $k\geq0$ and $1\leq l \leq 6$.
Then the canonical bundle $K_{X}$ is $\pi^{ *}\O(k-1)$
and the Hodge numbers are as follows:
$$
\begin{matrix}
&&1\\
&0&&0\\
k&&10+10k&&k\\
&0&&0\\
&&1
\end{matrix}
$$
In particular $X$ is rational if $1\leq n \leq 6$,
a K3 surface if $7\leq n \leq 12$ and $\kappa(X)=1$ if $n\geq 13$.
There are $(n+2)$ singular fibers, $D_{0}=\pi^{-1}(0)$,
$D_{i}=\pi^{-1}(\sqrt[n]{-4/27}\zeta_n^i)$ ($1\leq i\leq n$)
and
$Y_{\infty}=\pi^{-1}(\infty)$. 
The multiplicative fibers are $D=D_{0}+\cdots+D_{n}$.
$D_{0}$ is of type $I_n$ and $D_{i}$
is of type $I_1$ for $1\leq i\leq n$.
The type of $Y_{\infty}$ is as follows:
\begin{center}
\begin{tabular}{c|cccccc}
$l$&1&2&3&4&5&6\\
\hline
$Y_{\infty}$&II*&IV*&${\mathrm I}_0^*$&IV&II&(smooth)
\end{tabular}
\end{center}
The number of multiplicative fibers is $s=n+1$.
We put
$
\Sigma=\{0,\sqrt[n]{-4/27}\zeta_n^i~(
1\leq i\leq n)\}\subset C$
and $U=X-D$.

\medskip

The purpose of this section is to prove the following:
\begin{thm}\label{two}
Suppose that $n$ is a prime number with
$2\leq n\leq29$. 
Then we have
$$
\rank~\dlog
\vg(U,\K_2)=2.
$$
The dlog image is generated by
$$
\dlog\left\{
\frac{Y-X}{Y+X},-\frac{t^n}{X^3}
\right\},\quad
\dlog\left\{
\frac{iY-(X+2/3)}{iY+(X+2/3)},\frac{t^n+4/27}{(X+2/3)^3}
\right\}
$$
\end{thm}
Note $\dlog
\vg(U,\K_2)=
\dlog
\vg(U-Y_{\infty},\K_2)$ (Lemmas \ref{d=b} and \ref{d=bmore}).
\subsection{Proof of Theorem \ref{two}}
Let $\partial:\vg(U,\K_2)\to \Z^{\op 1+n}$ be the boundary map
where the base $(0,\cdots,\os{b}{1},\cdots,0)\in\Z^{\op 1+n}$ 
corresponds to
the fiber $D_b$ for $0\leq b \leq n$.
We want to show that the rank of $\partial\vg(U,\K_2)$ is 2.
A direct calculation yields
$$
\partial\left\{
\frac{Y-X}{Y+X},-\frac{t^n}{X^3}
\right\}=(n,0,\cdots,0),
$$
$$
\partial\left\{
\frac{iY-(X+2/3)}{iY+(X+2/3)},-\frac{t^n+4/27}{(X+2/3)^3}
\right\}=(0,1,\cdots,1).
$$
Therefore we have $\rank~\partial\vg(U,\K_2)\geq2$ and hence
it is enough to show
$\rank~\partial\vg(U,\K_2)\leq2$.
To do this we use Theorem \ref{key3}.
\subsubsection{Step 1 : Choice of $p\geq 5$ and $\pi_R:X_R\to C_R$}
Let $p\not\vert6n$ be a prime number.
Let $
K=\Q_p(\sqrt{-1},\sqrt[n]{-4/27},\zeta_n)
$ be an unramified extension over $\Q_p$
and $R$ the ring of integers in $K$.
Let $$\pi_R:X_R\to C_R$$ be the elliptic surface over $R$ obtained
from the defining equation $Y^2=X^3+X^2+t^n$ in a natural way
such that $X_K:=X_R\times_RK$ is minimal.
One can easily check that it is an elliptic surface in the sense
of \S \ref{setupsect} and
satisfies {\bf (Rat)}.
Since $6\not\vert n$, the fiber $Y_{\infty}$ is not
a multiplicative type. Therefore
$H^3_\et(X_{\ol{K}},\Z_p)$ is torsion free
by Lemma \ref{simplyc} and
Proposition \ref{cohprop} \eqref{cohprop3}.
Thus 
the condition in Theorem \ref{key3} (1) is satisfied.
Let us see the condition in Theorem \ref{key3} (2).
Let $F=\Q(\sqrt{-1},\sqrt[n]{-4/27},\zeta_n)$.
One can check that $\pi_K:X_K\to C_K$ is defined over $F$,
which we write $\pi_F:X_F\to C_F$.
It is easy to see that $\pi_F:X_F\to C_F$ satisfies {\bf(Rat)}.
We discuss the condition
\begin{equation}\label{U2}
H^2(X_{\ol{F}},\Z/p(2))^{G_F}=0.
\end{equation}
To do this, we see the Frobenius action.
Suppose that $X_F$ has a good reduction at
a prime $\frak l$ of $F$ which is prime to $p$.
Let $\l$ be a rational prime such that ${\frak l}\vert \l$.
Let $\kappa$ be the residue field at $\frak l$
and put $t=[\kappa:\F_\l]$.
Let $\Frob_{\frak l}\in\Gal(\ol{\kappa}/\kappa)$ be the Frobenius
automorphism.
By the proper smooth base change theorem we have
$$
H^2_\et(X_{\ol{F}},\Z/p(2))\cong
H^2_\et(X_{\ol{\kappa}},\Z/p(2)).
$$
By Proposition \ref{cohprop} \eqref{cohprop3}
we have
$$
H^2_\et(X_{\ol{\kappa}},\Z/p(2))\cong
H^2_\et(X_{\ol{\kappa}},\Z_p(2))\ot\Z/p\Z
$$
and $H^2_\et(X_{\ol{\kappa}},\Z_p(2))$ is torsion free.
Denote by $\Frob_{\frak l}^{\text{arith}}$ the action of 
$\Frob_{\frak l}$ on 
$H^2_\et(X_{\ol{\kappa}},\Z_p(j))$. 
If $p$ satisfies
\begin{equation}\label{ph20}
p\not\vert
\det(1-\Frob_{\frak l}^{\text{arith}}:
H^2_\et(X_{\ol{\kappa}},\Z_p(2))),
\end{equation}
then \eqref{U2} holds.
There is the geometric Frobenius endomorphism on $X_{\ol{\kappa}}$.
Denote by $\Frob_{\frak l}^{\text{geo}}$ the action 
of the geometric Frobenius on 
$H^2_\et(X_{\ol{\kappa}},\Z_p(j))$.
Then $\Frob_{\frak l}^{\text{geo}}\cdot
\Frob_{\frak l}^{\text{arith}}=
\Frob_{\frak l}^{\text{arith}}\cdot
\Frob_{\frak l}^{\text{geo}}$ 
is the identity.
By the Poincare duality theorem, the dual of
$H^2_\et(X_{\ol{\kappa}},\Z_p(2))$ is 
$H^2_\et(X_{\ol{\kappa}},\Z_p)$. Thus we have
\begin{equation}\label{primeonh21}
\det(1-\Frob_{\frak l}^{\text{arith}}:
H^2_\et(X_{\ol{\kappa}},\Z_p(2)))=
\det(1-\Frob_{\frak l}^{\text{geo}}:
H^2_\et(X_{\ol{\kappa}},\Z_p)).
\end{equation}
$X_F$ is defined over $\Q$ and 
$X_\kappa$ is defined over the prime field $\F_\l$, which
we denote by $X_\Q$ and $X_{\F_\l}$ respectively.
Let $\Frob_{\l}^{\text{geo}}$ be the geometric Frobenius
action on $H^2_\et(X_{\ol{\F}_\l},\Z_p)
=H^2_\et(X_{\ol{\kappa}},\Z_p)$ and
$\alpha_i$ its eigenvalues. 
Since $\Frob_{\frak l}^{\text{geo}}=(\Frob_{\l}^{\text{geo}})^t$
we have
\begin{equation}\label{primeonh21s}
\det(1-\Frob_{\frak l}^{\text{geo}}:
H^2_\et(X_{\ol{\kappa}},\Z_p))=
\prod_{i=1}^{10+12k}(1-\alpha_i^{t}).
\end{equation}
The eigenvalues $\alpha_i$ are computed from the zeta function
of $X_{\F_\l}$.
Namely letting $\nu_m(X_{\F_\l})$ be the cardinality of
the $\F_{\l^m}$-rational points of $X_{\F_\l}$, 
the zeta function $Z(X_{\F_\l},T)\in \Z[[T]]$ is defined
as follows:
$$
Z(X_{\F_\l},T)\os{\mathrm{def}}{=}\exp\sum_{m=1}^\infty
\frac{\nu_m(X_{\F_\l})}{m}T^m.
$$
Due to the Lefschetz trace formula
we have
$$
Z(X_{\F_\l},T)=
\frac{1}{(1-T)(1-\l^{2}T)
\det(1-\Frob_{\l}^{\text{geo}}T:
H^2_\et(X_{\ol{\F}_\l},\Z_p))}
$$
Hence we have
\begin{equation}\label{primeonh23}
\nu_m(X_{\F_\l})=1+\l^{2m}+\sum_{i=1}^{10+12k}\alpha_i^m.
\end{equation}
Moreover it follows from the Poincare duality that we have
\begin{equation}\label{primeonh22}
\prod_{i=1}^{10+12k}\alpha_i=\pm \l^{10+12k}.
\end{equation}
Let us compute $\nu_n(X_{\F_\l})$.
Put
$$
X^o:=\Spec\F_\l[X,Y,t]/(Y^2-X^3-X^2-t^n)\hra X_{\F_\l}.
$$
Since $X_{\F_\l}=(X^o-\{t=0\})\cup D_0 \cup Y_\infty \cup e(\P^1)$, 
we have
\begin{equation}\label{primeonh25}
\nu_m(X_{\F_\l})=1+(12k-n+11)\l^m+\nu_m(X^o)
\quad \text{unless }6\vert n.
\end{equation}
Suppose that $n$ is a prime number and
$(\l$ mod $n)$ is a generator
of $(\Z/n)^*$.  
Then $\F_{\l^m}\to\F_{\l^m}$ $x\mapsto x^n$ is bijective
unless $n-1\vert m$ and hence we have
\begin{equation}\label{primeonh26f}
\nu_m(X^o)=\l^{2m}
\quad\text{unless }n-1\vert m.
\end{equation}
By \eqref{primeonh23}, \eqref{primeonh25}
and \eqref{primeonh26f} we have
\begin{equation}\label{primeonh27}
\sum_{i=1}^{10+12k}\alpha_i^m=
(12k-n+11)\l^m
\quad\text{unless }n-1\vert m.
\end{equation}
\begin{claim}\label{weil0}
Suppose $n\geq11$. Then
\eqref{primeonh22} and
\eqref{primeonh27} yield
$$
\sum_{i=1}^{10+12k}\alpha_i^{n-1}=
\begin{cases}
(10+12k)\l^{n-1} & \text{if }
\prod_i\alpha_i=-\l^{10+12k}\\
10 & \text{if }
\prod_i\alpha_i=\l^{10+12k}\\
\end{cases}
$$
$$
\sum_{i=1}^{10+12k}\alpha_i^{2n-2}=
(10+12k)\l^{n-1}  \quad\text{if }
\prod_i\alpha_i=\l^{10+12k}.
$$
\end{claim}
\begin{pf}
Exercise on symmetric polynomials.
\end{pf}
By the Weil-Riemann conjecture we have
\begin{equation}\label{prim29}
\vert\alpha_i^\sigma\vert=\l\quad \text{for any }
\sigma:\ol{\Q}\to\C.
\end{equation}
Therefore Claim \ref{weil0} implies $\alpha_i^{2n-2}=\l^{2n-2}$ 
for all $1\leq i \leq 10+12k$.
Since $2n-2$ is divided by
$t=[\kappa:\F_\l]$, \eqref{ph20} holds if $p\not|(1-\l^{2n-2})$ for 
$n\geq11$.
There exist infinitely many $\l$ such that
$(\l$ mod $p)$ is a generator
of $(\Z/p)^*$ and  
$(\l$ mod $n)$ is a generator
of $(\Z/n)^*$.  
For such $\l$,
if $p-1\not\vert 2n-2$, then $p\not\not\vert (1-\l^{2n-2})$
and hence \eqref{ph20}.
Thus \eqref{U2}
holds if $p-1\not\vert 2n-2$ for $n\geq 11$.

In the cases $n=2,3,5,7$, we check \eqref{ph20}
as the case may be.

\smallskip

\noindent\underline{Case n=2}.
We put $\l=13$. Then $\kappa=\F_{13}$ ($t=1$).
By a direct calculation we have
$$
\nu_1(X_{\F_{13}})=13^2+10\cdot 13+1\quad(\Longrightarrow
\sum_{i=1}^{10}\alpha_i=10\cdot 13).
$$
Then it follows from the Weil-Riemann conjecture \eqref{prim29}
that we have $\alpha_i=13$ for all $1\leq i\leq 10$.
Thus if $p\not\vert13-1$ (i.e. $p\geq5$) 
then \eqref{ph20} and hence \eqref{U2} hold.

\noindent\underline{Case n=3}.
We put $\l=7$. Then $\kappa=\F_{7^6}$ ($t=6$).
$$
\nu_1(X_{\F_7})=7^2+10\cdot 7+1.
$$
Therefore we have $\alpha_i=7$ for all $1\leq i\leq 10$.
Thus if $p\not\vert7^6-1$ ($\Longleftrightarrow$ $p\not=2,3,19,43$) 
then \eqref{ph20} holds.
On the other hand let $\l=13$. Then $\kappa=\F_{13^3}$ ($t=3$)
and
$$
\nu_1(X_{\F_{13}})=13^2+10\cdot 13+1.
$$
Therefore we have $\alpha_i=13$ for all $1\leq i\leq 10$
and if $p\not\vert13^3-1$ ($\Longleftrightarrow$ $p\not=2,3,61$)
then \eqref{ph20} holds.
As a result, if $p\geq 5$ then \eqref{U2} holds.

\noindent\underline{Case n=5}.
We put $\l=19$. Then $\kappa=\F_{19^2}$ ($t=2$).
$$
\nu_m(X_{\F_5})=\begin{cases}
19^{2m}+6\cdot 19^m+1&2\not\vert m\\
19^{4}+10\cdot 19^2+1&m=2.
\end{cases}
$$
Therefore we have $\alpha_i^2=19^2$ for all $1\leq i\leq 10$.
Thus if $p\geq7$ then \eqref{ph20} holds.

\noindent\underline{Case n=7}.
We put $\l=13$. Then $\kappa=\F_{13^2}$ ($t=2$).
$$
\nu_m(X_{\F_{13}})=\begin{cases}
13^{2m}+16\cdot 13^m+1&2\not\vert m\\
13^{4}+22\cdot 13^2+1&m=2.
\end{cases}
$$
Therefore we have $\alpha_i^2=13^2$ for all $1\leq i\leq 22$.
Thus if $p\not\vert2\cdot3\cdot7\cdot13$ then \eqref{ph20} holds.

\medskip

Summarizing the above, 
the elliptic surface $\pi_R:X_R\to C_R$ over $R$
satisfies
the conditions in Theorem \ref{key3}
in the following cases.
\begin{center}
\begin{tabular}{c|c|c|c|c|c}
$n$&2&3&5&7&$n\geq11$\\
\hline
$p$&$p\geq5$&$p\geq5$&$p\geq7,~p\not=19$
&$p\not\vert2\cdot3\cdot7\cdot13$
&$p-1\not\vert 2(n-1)$
\end{tabular}
\end{center}

\subsubsection{Step 2 : Computation of $\dim_{\F_p}
\Phi(X_R,D_R)_{\F_p}$}
Let $p$,
$K=\Q_p(\sqrt{-1},\sqrt[n]{-4/27},\zeta_n)$
and $\pi_R:X_R\to C_R$ be as above.
Put
$$
Y'=\sqrt{-1}Y,\quad X'=-X-\frac{2}{3},\quad 
t'=\left(\sqrt[n]{\frac{-4}{27}}\right)^{-1}t,
\quad
t_i=\zeta_n^{-i}t'-1 \quad (1\leq i \leq n).
$$
Then 
\begin{equation}\label{elsurfeqc1}
Y^{\prime2}=X^{\prime3}+X^{\prime2}+\frac{4}{27}
(t^{\prime n}-1)
=X^{\prime3}+X^{\prime2}+\frac{4}{27}
((t_i+1)^n-1).
\end{equation}
Let $q_i\in t_iR[[t_i]]$ be the formal power series such that
\begin{equation}\label{qexp1}
j=\frac{1}{q_i}+744+196884q_i+\cdots=
\frac{-432}{(t_i+1)^n((t_i+1)^n-1)}.
\end{equation}
$q_i$ 
is the period of the Tate curve in the neighborhood of 
the singular fiber $D_i$ ($1\leq i \leq n$).
Conversely the $q_i$-expansion of $t_i$ is as follows:
\begin{equation}\label{qexp2}
t_i=-\frac{432}{n}q_i+\frac{41472n+93312}{n^2}q_i^2+\cdots.
\end{equation}
Then $R((q_i))= R((t_i))$ and the isomorphism
between the Tate curve $E_{q_i,R((q_i))}$ \eqref{tate0}
and $U_R\times_{S_R}\Spec R((q_i))$ 
is given by
$$
X'=4\frac{x+(1-E_4^{1/4})/12}{E_4^{1/2}},
\quad
Y'=4\frac{2y+x}{E_4^{3/4}}
$$
where $E_4$ is the Eisenstein series:
$$
E_4=1+240\sum_{n=1}^\infty\frac{n^3q_i^n}{1-q_i^n},
\quad E_4^{1/4}=1+60q_i-4860q_i^2+\cdots.
$$
In particular we have
\begin{equation}\label{qexp3}
\frac{dX'}{Y'}=E_4^{1/4}\frac{dx}{2y+x}=E_4^{1/4}\frac{du}{u}.
\end{equation}
$\vg(X_R,\Omega^2_{X_R/R}(\log(D_R)))$ is generated by
$$
t^{\prime a}dt'\frac{dX'}{Y'},\quad
\frac{dt'}{t'}\frac{dX'}{Y'},\quad
\frac{dt'}{t'-\zeta_n^b}\frac{dX'}{Y'},
\quad (0\leq a \leq k-1,~  1\leq b \leq n)
$$
where $n=6k+l$ with $k\geq 0$ and $1\leq l \leq 6$.
Since
$$
\partial_{\DR}\left(
t^{\prime a}dt'\frac{dX'}{Y'}
\right)=(0,\cdots,0),\quad
\partial_{\DR}\left(
\frac{dt'}{t'}\frac{dX'}{Y'}
\right)=(\sqrt{-1},0,\cdots,0),
$$
$$
\partial_{\DR}\left(
\frac{dt'}{t'-\zeta_n^b}\frac{dX'}{Y'}
\right)=(0,\cdots,\os{b}{1},\cdots,0),
$$
we have
\begin{align*}
\vg(X_R,\Omega^2_{X_R/R}(\log D_R))_{\Z_p}
&=\sum_{a=0}^{k-1}
Rt^{\prime a}dt'\frac{dX'}{Y'}+
\Z_p\sqrt{-1}\frac{dt'}{t'}\frac{dX'}{Y'}+
\sum_{b=1}^n\Z_p\frac{dt'}{t'-\zeta_n^b}\frac{dX'}{Y'}\\
&=\sum_{a=0}^{k-1} Rt^{a}_idt_i\frac{dX'}{Y'}+
\Z_p\sqrt{-1}\frac{dt_i}{t_i+1}\frac{dX'}{Y'}+
\sum_{b=1}^n\Z_p\frac{dt_i}{t_i+1-\zeta_n^{b-i}}\frac{dX'}{Y'}.
\end{align*}

\medskip

We show $\dim_{\F_p}\D(X_R,D_R)_{\F_p}\leq2$
as the case may be.
We check it only in case $n=7$.
The other cases are similar (left to the reader).

\medskip

Let $n=7$.
We take $p=11$. Then $K=\Q_{11}(\sqrt{-1},\zeta_7)
=\Q_{11}(\zeta)$ $(\zeta:=\sqrt{-1}\zeta_7)$, 
$R=\Z_{11}[\zeta]$ and
$[K:\Q_{11}]=6$.
Fix a cyclotomic basis $\underline{\mu}=\{1,\zeta,\zeta^2,\cdots,
\zeta^5\}$.
Put
\begin{align*}
a&=\zeta_7+\zeta_7^2+\zeta^4_7=4+11+11^2+7\cdot 11^3+\cdots\\
b&=\zeta_7^3+\zeta_7^5+\zeta^6_7
=6+9\cdot 11+9\cdot 11^2+3\cdot11^3+\cdots.
\end{align*}
The minimal polynomial of $\zeta_7$ over $\Q_{11}$ is 
$x^3-ax^2+bx-1$. 
$\vg(X_R,\Omega^2_{X_R/R}(\log D_R))_{\Z_{11}}$
is generated by
$$
Rdt'\frac{dX'}{Y'},\quad
\Z_{11}\sqrt{-1}\frac{dt'}{t'}\frac{dX'}{Y'},\quad
\Z_{11}\frac{dt'}{t'-\zeta_7^b}\frac{dX'}{Y'} \quad(1\leq b \leq 7).
$$
Let us compute
$$
\phi_i:\vg(X_R,\Omega^2_{X_R/R}(\log D_R))_{\Z_{11}}
\lra \prod_{k\geq1}(\Z_{11}/k^2\Z_{11})^{\op 6}
$$
for $1\leq i\leq n$.
Since 
\begin{align*}
&\dlog\left\{
\frac{Y-X}{Y+X},-\frac{t^7}{X^3}
\right\}=7\sqrt{-1}\frac{dt'}{t'}\frac{dX'}{Y'},\\
&\dlog\left\{
\frac{iY-(X+2/3)}{iY+(X+2/3)},-\frac{t^7+4/27}{(X+2/3)^3}
\right\}=\frac{d(t^{\prime 7})}{t^{\prime 7}-1}\frac{dX'}{Y'}
\end{align*}
we have $$\phi_i\left(\sqrt{-1}\frac{dt'}{t'}\frac{dX'}{Y'}
\right)=0,
\quad
\phi_i\left(\frac{d(t^{\prime 7})}{t^{\prime 7}-1}
\frac{dX'}{Y'}\right)=0.$$
Let
$$
\frac{dt'}{t'-\zeta_7^b}\frac{dX'}{Y'}=
f_b(q_i)\frac{dq_i}{q_i}\frac{du}{u}.
$$
Express
\begin{align*}
f_b(q_i)&=E_4^{1/4}\frac{1}{t_i+1-\zeta_7^{b-i}}q_i\frac{dt_i}{dq_i}\\
&=E_4^{1/4}\frac{1}{1-\zeta_7^{b-i}}\sum_{m\geq0}
\left(\frac{-t_i}{1-\zeta_7^{b-i}}\right)^m\cdot
q_i\frac{dt_i}{dq_i}\\
&=
\sum_{k\geq1}\frac{a_kq_i^k}{1-q_i^k}
+\frac{b_k\zeta q_i^k}{1-\zeta q_i^k}
+\cdots+\frac{f_k\zeta^5 q_i^k}{1-\zeta^5 q_i^k}
\end{align*}
for $b\not\equiv i\mod 7$ and
\begin{align*}
f_b(q_i)&=E_4^{1/4}\frac{1}{t_i}q_i\frac{dt_i}{dq_i}\\
&=
1+\sum_{k\geq1}\frac{a_kq_i^k}{1-q_i^k}
\end{align*}
for $b=i$. By definition,
$\phi_i(\frac{dt'}{t'-\zeta_7^b}\frac{dX'}{Y'})=
(\bar{a}_k,\cdots,\bar{f}_k)_{k\geq1}$
where $\bar{a}_k=(a_k$ mod $k^2\Z_{11})$. The
following is the table for $k=11,22,33,44$
(don't check it by hand).
\begin{center}
\begin{tabular}{c|cccccc|cccccc}
$b-i$
&$\bar{a}_{11}$&$\bar{b}_{11}$&$\bar{c}_{11}$
&$\bar{d}_{11}$&$\bar{e}_{11}$&$\bar{f}_{11}$
&$\bar{a}_{22}$&$\bar{b}_{22}$&$\bar{c}_{22}$
&$\bar{d}_{22}$&$\bar{e}_{22}$&$\bar{f}_{22}$\\
\hline
1&
$99$&$0$&$99$&$0$&$55$&$0$&$11$&$0$&$99
$&$0$&$99$&$0$\\
2&
$77$&$0$&$55$&$0$&$88$&$0$&$110$&$0$&$55
$&$0$&$99$&$0$\\
3&
$11$&$0$&$55$&$0$&$33$&$0$&$99$&$0$&$88
$&$0$&$11$&$0$\\
4&
$66$&$0$&$88$&$0$&$99$&$0$&$99$&$0$&$88
$&$0$&$44$&$0$\\
5&
$33$&$0$&$33$&$0$&$55$&$0$&$110$&$0$&$33
$&$0$&$22$&$0$\\
6&
$22$&$0$&$33$&$0$&$33$&$0$&$99$&$0$&$0
$&$0$&$88$&$0$\\
0&
$55$&$0$&$0$&$0$&$0$&$0$&$77$&$0$&$0
$&$0$&$0$&$0$\\
\end{tabular}
\end{center}
\begin{center}
\begin{tabular}{c|cccccc|cccccc}
$b-i$
&$\bar{a}_{33}$&$\bar{b}_{33}$&$\bar{c}_{33}$
&$\bar{d}_{33}$&$\bar{e}_{33}$&$\bar{f}_{33}$
&$\bar{a}_{44}$&$\bar{b}_{44}$&$\bar{c}_{44}$
&$\bar{d}_{44}$&$\bar{e}_{44}$&$\bar{f}_{44}$\\
\hline
1&
$99$&$0$&$88$&$0$&$99$&$0$&$33$&$0$&$55
$&$0$&$55$&$0$\\
2&
$0$&$0$&$99$&$0$&$55$&$0$&$77$&$0$&$11
$&$0$&$66$&$0$\\
3&
$88$&$0$&$11$&$0$&$44$&$0$&$11$&$0$&$0
$&$0$&$11$&$0$\\
4&
$33$&$0$&$55$&$0$&$88$&$0$&$55$&$0$&$55
$&$0$&$0$&$0$\\
5&
$55$&$0$&$66$&$0$&$11$&$0$&$11$&$0$&$55
$&$0$&$55$&$0$\\
6&
$88$&$0$&$44$&$0$&$66$&$0$&$33$&$0$&$66
$&$0$&$55$&$0$\\
0&
$0$&$0$&$0$&$0$&$0$&$0$&$22$&$0$&$0
$&$0$&$0$&$0$\\
\end{tabular}
\end{center}
Similarly, we put
$$
\zeta^bdt'\frac{dX'}{Y'}=
\zeta^{b+8i}dt_i\frac{dX'}{Y'}=
g_b(q_i)\frac{dq_i}{q_i}\frac{du}{u}
$$
$$
g_b(q_i)=\zeta^{b+8i}E_4^{1/4}q_i\frac{dt_i}{dq_i}
=\sum_{k\geq1}\frac{a_kq_i^k}{1-q_i^k}+
\frac{b_k\zeta q_i^k}{1-\zeta q_i^k}
+\cdots+\frac{f_k\zeta^5 q^k_i}{1-\zeta^5 q_i^k}.
$$
\begin{center}
\begin{tabular}{c|cccccc|cccccc}
$b+8i$
&$\bar{a}_{11}$&$\bar{b}_{11}$&$\bar{c}_{11}$
&$\bar{d}_{11}$&$\bar{e}_{11}$&$\bar{f}_{11}$
&$\bar{a}_{22}$&$\bar{b}_{22}$&$\bar{c}_{22}$
&$\bar{d}_{22}$&$\bar{e}_{22}$&$\bar{f}_{22}$\\
\hline
0&
$2$&$0$&$0$&$0$&$0$&$0$&$93$&$0$&$0
$&$0$&$0$&$0$\\
1&
$0$&$2$&$0$&$42$&$0$&$54$&$109$&$95$&$106
$&$92$&$8$&$101$\\
2&
$25$&$0$&$86$&$0$&$79$&$0$&$49$&$0$&$79
$&$0$&$22$&$0$\\
3&
$0$&$0$&$0$&$44$&$0$&$79$&$86$&$0$&$101
$&$66$&$1$&$29$\\
4&
$0$&$0$&$42$&$0$&$44$&$0$&$66$&$0$&$48
$&$0$&$68$&$0$\\
5&
$0$&$54$&$0$&$25$&$0$&$86$&$60$&$101$&$116
$&$49$&$88$&$37$\\
\end{tabular}
\end{center}
\begin{center}
\begin{tabular}{c|cccccc|cccccc}
$b+8i$
&$\bar{a}_{33}$&$\bar{b}_{33}$&$\bar{c}_{33}$
&$\bar{d}_{33}$&$\bar{e}_{33}$&$\bar{f}_{33}$
&$\bar{a}_{44}$&$\bar{b}_{44}$&$\bar{c}_{44}$
&$\bar{d}_{44}$&$\bar{e}_{44}$&$\bar{f}_{44}$\\
\hline
0&
$114$&$0$&$0$&$0$&$0$&$0$&$53$&$0$&$0
$&$0$&$0$&$0$\\
1&
$0$&$24$&$0$&$93$&$0$&$28$&$17$&$27$&$91
$&$83$&$22$&$3$\\
2&
$113$&$0$&$7$&$0$&$106$&$0$&$35$&$0$&$43
$&$0$&$67$&$0$\\
3&
$0$&$24$&$0$&$88$&$0$&$28$&$77$&$0$&$60
$&$110$&$44$&$38$\\
4&
$20$&$0$&$15$&$0$&$113$&$0$&$99$&$0$&$17
$&$0$&$51$&$0$\\
5&
$0$&$60$&$0$&$40$&$0$&$52$&$35$&$3$&$65
$&$35$&$77$&$72$\\
\end{tabular}
\end{center}
The above tables yields that
the kernel of $\phi_1$
is contained in
$$ 
\Z_{11}\sqrt{-1}\frac{dt'}{t'}\frac{dX'}{Y'}
+\Z_{11}\frac{dt^{\prime 7}}{t^{\prime 7}-1}\frac{dX'}{Y'}
+{11} \cdot\vg(X_R,\Omega^2_{X_R/R}(\log D_R))_{\Z_{11}}.
$$
Hence we have $\dim_{\F_{11}}\D(X_R,D_R)_{\F_{11}}\leq2$.
This completes the proof.

\subsection{Indecomposable parts of $K_1(X)$ with higher rank}
Let $X$ be a nonsingular variety over $\C$.
Let $Y\subset X$ be a closed subscheme.
The Adams operations on $K_i^Y(X)_\Q=K_i^Y(X)\ot\Q
\cong K'_i(Y)\ot\Q$ give rise to the decomposition
$\bigoplus_j K_i^Y(X)^{(j)}$ (\cite{soulecanada}). 
In case $i=1$ we have the decomposition
$K_1(X)_\Q=\bigoplus_{j=1}^{\dim X+1} K_1(X)^{(j)}$.
Of particular interest to us here is $K_1(X)^{(2)}$, which is
isomorphic to
Bloch's higher Chow group $\CH^2(X,1)\ot\Q$.
Let $\NS(X)$ be the N\'eron-Severi group of $X$.
Then there is the natural map $\C^*\ot \NS(X)_\Q\to K_1(X)^{(2)}$.
We denote its cokernel by $K_1^{\mathrm{ind}}(X)^{(2)}$ and 
call the {\it indecomposable $K_1$} of $X$.
An element $x\in K_1(X)^{(2)}$ is called {\it indecomposable}
if it is non trivial in $K_1^{\mathrm{ind}}(X)^{(2)}$.
\begin{lem}\label{dec0lem}
Let $\pi:X\to C$ be a minimal elliptic surface
over $\C$ and $U^0$ the complement of the all singular fibers.
Let $\NS(X)'\subset \NS(X)$ be the subgroup generated by
the irreducible components of singular fibers and the section $e(C)$.
Then there is an exact sequence
\begin{equation}\label{dec0}
K_2(U^0)_\Q\os{\partial}{\lra} \Q^{\op s}
\lra K_1(X)^{(2)}/(\C^*\ot\NS(X)').
\end{equation}
\end{lem}
\begin{pf}
Let $D_i$ be the multiplicative fibers
and $Y_i$ the singular fibers of
other types.
Let $D_i=\sum_k D_i^{(k)}$ and $Y_i=\sum_k Y_i^{(k)}$ be
the irreducible decompositions.
Quillen's localization exact sequence
(\cite{Sr} Prop.(5.15))
yields an exact sequence
\begin{equation}\label{dec1}
K_2(U^0)^{(2)}\lra
\bigoplus_{i} K'_1(D_i)^{(1)}\op \bigoplus_{j}K'_1(Y_j)^{(1)}
\lra K_1(X)^{(2)}
\lra K_2(U^0)^{(1)}.
\end{equation}
On the other hand, there are exact sequences
\begin{equation}\label{dec2}
\bigoplus_k K_1(D_i^{(k)})
\lra
K'_1(D_i) \lra \Z \lra 0,
\end{equation}
\begin{equation}\label{dec2i}
\bigoplus_k K_1(Y_j^{(k)})
\lra
K'_1(Y_j)\lra 0,
\end{equation}
\begin{equation}\label{dec3}
\bigoplus_k \C^*\cdot [D_i^{(k)}]
\lra
K'_1(D_i)^{(1)} \lra \Q\lra 0,
\end{equation}
\begin{equation}\label{dec3i}
\bigoplus_k \C^*\cdot [Y_j^{(k)}]
\lra
K'_1(Y_j)^{(1)}\lra 0.
\end{equation}
(As for $D_i$, the above is shown in \S \ref{modelsect}.
However the same argument also works for $Y_j$).
Let $\NS(X)^{\prime\prime}\subset \NS(X)'$ be the subgroup
generated by all irreducible components of the singular fibers.
Then
$
\NS(X)'=\NS(X)^{\prime\prime}\op \Z[e(C)].
$
We claim that the map 
$$
\C^*\ot \NS(X)'_\Q\lra K_1(X)^{(2)}
$$
is injective. In fact let $K_1(X)^{(2)}\to
K_1(D_i^{(k)})^{(2)\op}\op
K_1(Y_j^{(k)})^{(2)\op}\cong \C^{*\op}$
be the pull-back.
Then the composition $\C^*\ot \NS(X)'\lra K_1(X)^{(2)}\to 
\C^{*\op}$
is given by the intersection matrix on $NS(X)'$.
Since the intersection pairing on $\NS(X)'$ is
non-degenerate (cf. proof of Lemma \ref{intlem}),
it is bijective.
Now a commutative diagram
$$
\begin{CD}
@.0@.0\\
@.@AAA@AAA\\
@.
\Q^{\op n+1}@>>>
K_1(X)^{(2)}/\C^*\ot\NS(X)^{\prime\prime}\\
@.@A{\partial_1}AA@AAA\\
K_2(U^0)^{(2)}@>{\partial_2}>>
\bigoplus_{i} K'_1(D_i)^{(1)}
\op \bigoplus_{j}K'_1(Y_j)^{(1)}@>>> K_1(X)^{(2)}\\
@.@AAA@AAA\\
@. \bigoplus_{i,k} \C^*\ot[D_i^{(k)}]\op
\bigoplus_{j,k} \C^*\ot[Y_j^{(k)}]@>>>
\C^*\ot\NS(X)^{\prime\prime}\\
@.@.@AAA\\
@.@.0
\end{CD}
$$
with $\partial=\partial_1\partial_2$ yields an exact sequence
\begin{equation}\label{dec4}
K_2(U^0)_\Q\os{\partial}{\lra} \Q^{\op s}\lra
K_1(X)^{(2)}/\C^*\ot\NS(X)^{\prime\prime}.
\end{equation}
The rest of the proof is to show 
\begin{equation}\label{dec6}
\Image\left(
\bigoplus_{i} K'_1(D_i)^{(1)}_\Q
\to K_1(X)^{(2)}\right)
\bigcap \left(\C^*\ot[e(C)]\right)=0.
\end{equation}
Let $X_t$ be a smooth fiber. Let 
$K_1(X)^{(2)}\to K_1(X_t)^{(2)}$ be the pull-back and 
$K_1(X_t)^{(2)}\to \C^*$
the norm map (also called the transfer map) 
for $X\to\Spec \C$. Then the composition
$\C^*\ot[e(C)]\to K_1(X)^{(2)}\to K_1(X_t)^{(2)}\to \C^*$ 
is bijective.
On the other hand the composition 
$K'_1(D_i)\to K_1(X)^{(2)}\to K_1(X_t)^{(2)}$ is
clearly zero. This shows that the image of $\op_i K'_1(D_i)^{(1)}$
does not intersect with $\C^*\ot [e(C)]$.
Thus we have \eqref{dec6}.
\end{pf}

\begin{cor}\label{declem}
Let $\pi:X\to\P^1$ be the minimal elliptic surface 
over $\C$
defined by 
$Y^2=X^3+X^2+t^n$.
Suppose that $n$ is prime to $30$.
Then there is an exact sequence
\begin{equation}\label{dec7}
K_2(U^0)_\Q\os{\partial}{\lra} \Q^{\op n+1}
\lra K_1^{\mathrm{ind}}(X)^{(2)}.
\end{equation}
\end{cor}
\begin{pf}
In this case $\NS(X)'_\Q=\NS(X)_\Q$ (\cite{stiller1} p.185 Example 7,
\cite{stiller2} Example 4).
\end{pf}

By Theorem \ref{two} and Corollary \ref{declem} we have
\begin{thm}\label{dectwo}
Let $\pi:X\to\P^1$ be the minimal elliptic surface 
over $\C$
defined by 
$Y^2=X^3+X^2+t^n$ and $D_i$ the multiplicative fibers.
Suppose that $n$ is a prime number with
$7\leq n\leq 29$. Then we have
$$
\dim~
\Image\left(
\bigoplus_{i=0}^n K'_1(D_i)_\Q
\lra K_1^{\mathrm{ind}}(X)^{(2)}\right)=n-1.
$$
\end{thm}



\noindent
Graduate School of Mathematics, Kyushu University,
Fukuoka 812-8581,
JAPAN

\medskip

\noindent
E-mail address : asakura@math.kyushu-u.ac.jp

\end{document}